%% file: main.tex
\pdfoutput=1

\documentclass{article}
\usepackage{arxiv}
\usepackage[utf8]{inputenc}
\usepackage[T1]{fontenc}    
\usepackage{url}            
\usepackage{booktabs}       
\usepackage{amsfonts}       
\usepackage{nicefrac}      
\usepackage{microtype} 

\usepackage{enumitem}

\usepackage{graphicx}
\usepackage{amsmath,amsthm, amssymb}
\numberwithin{equation}{section}
\usepackage[dvipsnames]{xcolor}
\usepackage[style=numeric, doi=true, url=true, eprint=true]{biblatex}
\usepackage{esint}
\usepackage{xfrac}
\usepackage[colorlinks=true, linkcolor=blue, citecolor=blue, urlcolor=magenta, breaklinks]{hyperref}
\usepackage[nameinlink]{cleveref}
\crefname{lemma}{Lemma}{Lemmata}
\crefname{thm}{Theorem}{Theorems}
\crefname{bmk}{Remark}{Remarks}
\crefname{assump}{Assumption}{Assumptions}
\crefname{cor}{Corollary}{Corollaries}
\crefname{defn}{Definition}{Definitions}
\crefname{equation}{Equation}{Equations}
\crefname{alg}{Algorithm}{Algorithms}

\crefname{prop}{Proposition}{Propositions}
\crefname{exmp}{Example}{Examples}
\crefname{appendix}{Appendix}{Appendices}
\AddToHook{env/lemma/begin}{\crefalias{thm}{lemma}}
\AddToHook{env/prop/begin}{\crefalias{thm}{prop}}
\AddToHook{env/cor/begin}{\crefalias{thm}{cor}}
\AddToHook{env/assump/begin}{\crefalias{thm}{assump}}
\AddToHook{env/defn/begin}{\crefalias{thm}{defn}}
\AddToHook{env/exmp/begin}{\crefalias{thm}{exmp}}
\AddToHook{env/bmk/begin}{\crefalias{thm}{bmk}}
\AddToHook{env/alg/begin}{\crefalias{thm}{alg}}

\crefname{section}{Section}{Sections}
\newtheoremstyle{break}
  {\topsep}   
  {\topsep}   
  {\itshape} 
  {}         
  {\bfseries}
  {.}         
  {\newline} 
  {}         
\theoremstyle{break}
\newtheorem{thm}{Theorem}[section]
\newtheorem{lemma}[thm]{Lemma}

\newtheorem{cor}[thm]{Corollary}
\newtheorem{assump}[thm]{Assumption}
\newtheorem{defn}[thm]{Definition} 
 
\newtheorem{bmk}[thm]{Remark}
\newtheorem{alg}[thm]{Algorithm}
\usepackage{mathtools} 
\DeclareMathOperator*{\esssup}{ess\,sup}
\title{Quantitative approximation for Fourier neural operators in controlled ordinary differential equations}
\renewcommand{\shorttitle}{}
\author{
  Dean Zenner \\
  Department of Mathematics\\
  Saarland University\\
  Saarbrücken, Germany \\
  \texttt{dean.zenner@num.uni-sb.de} \\
       \And
  Andreas Rupp \\
  Department of Mathematics\\
  Saarland University\\
  Saarbrücken, Germany \\
  \texttt{andreas.rupp@uni-saarland.de} \\
        \And
  Thomas Schuster \\
  Department of Mathematics\\
  Saarland University\\
  Saarbrücken, Germany \\
  \texttt{thomas.schuster@num.uni-sb.de} \\
}
\begin{document}
\maketitle
\begin{abstract}
Fourier neural operators (FNOs) avoid the curse of parametric complexity when approximating the solution operator of certain partial differential equations (PDEs). While recent studies primarily focus on PDEs, in this article we derive upper error bounds for approximating the solution operator of controlled ordinary differential equations (CODEs) $\dot{y}=f(y,u)$ using a FNO. Under a Carathéodory setting and suitable regularity assumptions, including a smallness condition, we prove that the size of the FNO grows at most algebraically in \(1/\varepsilon\), namely like \(\varepsilon^{-(1+d/(2(s-1)))}\) up to logarithmic factors for GELU networks, and is thus sub-quadratic for sufficiently smooth \(f\), where \(d=1+n_y+n_u\) and \(s\) denotes the Sobolev regularity of \(f\). Our proof is based on emulating the Picard iteration through the layers of the FNO while exploiting the differentiation formula of the discrete Fourier transform. This result contributes to the quantitative approximation theory of FNOs for physically motivated problems and provides a theoretical foundation for leveraging a FNO as efficient forward surrogate for the solution operator of CODEs, especially in inverse parameter estimation where repeated evaluations of the forward operator are needed.
\end{abstract}

\keywords{Fourier neural operators \and ordinary differential equations \and quantitative approximation \and inverse parameter estimation}

\smallskip
\section{Introduction}
Neural operators, acting between infinite-dimensional function spaces, constitute a natural extension of classical finite neural networks. In recent years, they have attracted considerable attention in the field of machine learning due to their ability to enabling surrogate modeling of high-dimensional and complex physical systems. Such systems can be modeled by partial differential equations (PDEs) or controlled ordinary differential equations (CODEs), whose numerical solving may become time consuming, particularly when repeated evaluations are required. Using neural operators to replace the solution operator of a PDE or CODE, where the latter also acts as an operator between function spaces due to the presence of control functions, may work out this issue. 
Once a neural operator is learned, it allows for a faster computational interference time while also maintaining an architecture that enables the processing of batched data. This also makes them an encouraging tool in solving (Bayesian) inverse problems via the Ensemble Kalman Inversion, especially in cases where a large data set of multiple experiments is present.

Various different architectures to realize a neural operator were introduced over the last years.
One of the most promising ones is the so-called Fourier neural operator (FNO) which employs the (inverse) Fourier transform inside the network's layered structure. While Universal Approximation Theorems of the FNO are well studied and suggest that any continuous operator can be approximated, 
quantitative approximation studies show that FNO faces exponential growth in model complexity, the so-called \textit{curse of parametric complexity}~\cite{lanthaler2025parametriccomplexityoperatorlearning}, even when approximating $C^r$- or Lipschitz-regular operators. The curse of parametric complexity is an infinite-dimensional analogue to the well-known curse of dimensionality common in high-dimensional approximation problems.

Recent research suggests that this curse of parametric complexity can be overcome when approximating the solution operator of certain PDEs by emulating numerical methods like the Picard iteration. While current studies therefore expand the class of PDEs which allow for breaking the curse using a FNO, CODEs have not been studied in this context, although they play an important role in modeling physics while also being solvable through the Picard iteration. Therefore in this article, we focus on CODEs of form $\dot{y}=f(t,y, u)$, where $f$ is a non-linear function and $u$ are the Lebesgue-integrable control functions.

This article starts with introducing a system of CODEs under a Carathéodory setting, that is that the controls $u$ do not need to be continuous, to determine appropriate mathematical assumptions to show the well-posedness of such a system while also deriving a continuous solution operator. This results in reformulating the CODE as a fixed point integral formula and thus, we follow the idea of emulating the Picard iteration by relying on the differentiation formula of the (discrete) Fourier transform to realize the integration of the controlled ordinary differential equation. Thus, the FNO constitutes a natural candidate for this approach. By sharpening our assumptions, which requires to assume more regularity on $u$ and $f$, while for the latter also a smallness condition is needed, we break the curse by showing that the size of the FNO scales at most algebraically, and under a sufficiently smooth $f$ even sub-quadratically up to logarithmic factors.

Our article is organized as follows:
\begin{itemize}[leftmargin=*]
    \item In~\cref{Sec:WPCODE} we introduce the controlled ordinary differential equation $\dot{y}=f(t,y, u)$ under suitable assumptions, such that we obtain the existence of a continuous control-to-solution operator. As we are interested in deriving error bounds for approximating this solution operator with a Fourier neural operator, we start to sharpening the assumptions on the CODE such that we can apply tools from Fourier analysis.
    \item The Fourier neural operator is being formally introduced and defined in~\cref{Sec:FNOCODE}. One essential key-point of the FNO is the usage of the (inverse) Fourier transform inside the network's architecture. Since deriving error bounds requires to quantify the network's size, we define the $\Psi$-FNO which relates to the numerical implementable version of a FNO. Afterwards, we state known quantitative approximation results regarding finite neural networks, as we need to approximate $f$ by such networks. At the end of~\cref{Sec:FNOCODE}, we propose a Picard type fixed-point iteration scheme which produces an approximate output to the true solution of the CODE while also allowing to be emulated through the layers of $\Psi$-FNO.
    \item In~\cref{Sec:QAT} we establish the main result. Under the sharpened assumptions from~\cref{Sec:WPCODE}, which include a smallness condition on $f$, we show that the size of the $\Psi$-FNO scales at most algebraically when approximating the control-to-solution operator. Thus, the size of the $\Psi$-FNO avoids the curse of parametric complexity.
    \item To conclude the article, we present in~\cref{Sec:AIPE} the application of a learned FNO in an inverse parameter estimation problem motivated by a physical example. Due to the nature of the given example, it is the case that a large data set of multiple physical experiments is present. We therefore propose to utilize the learned FNO in an Ensemble Kalman Inversion which is extended to a Kaczmarz-like structure.
\end{itemize}
\subsection{Related works}
Neural operators~\cite{kovachki2021neural} were first introduced in~\cite{li2020neuraloperatorgraphkernel} by realizing them as graph kernel networks. By doing so, they offer a new approach to operator learning compared to other methods such as DeepONet~\cite{Lu_2021} and PCA-Net~\cite{lanthaler2023operatorlearningpcanetupper}. Numerous neural operator architectures, such as the Fourier neural operator~\cite{li2021fourierneuraloperatorparametric}, the convolutional neural operator~\cite{raonic2023convolutionalneuraloperatorsrobust}, the Laplace neural operator~\cite{cao2023lnolaplaceneuraloperator}, the orthogonal polynomial neural operator~\cite{liu2023renderuntonumericsorthogonal} and the Wavelet neural operators~\cite{tripura2022waveletneuraloperatorneural}, have been developed in recent years. Owing to their strong empirical approximation capabilities and fast inference times, they have established themselves as state-of-the-art tools for surrogate modeling physical simulations~\cite{wagner2024neural}, which also makes them interesting objects to study in the context of inverse problems~\cite{Yin_2022,behroozi2025sensitivityconstrainedfourierneuraloperators,Yin_2023,Cao_2023,scherzer2026neuraloperatorssolvingnonlinear,wang2024latentneuraloperatorsolving,molinaro2023neuralinverseoperatorssolving}, due to the reliance on repeatedly evaluating the forward operator to numerically solve inverse problems. We refer the reader to~\cite{azizzadenesheli2024neuraloperatorsacceleratingscientific,kovachki2024operatorlearningalgorithmsanalysis} for an overview of neural operators and to~\cite{Nganyu_Tanyu_2023} for neural operators in the context of inverse problems.

One of the most important theoretical foundations in the field of machine learning are the Universal Approximation Theorems, which state that for a given continuous operator and any accuracy level, there exists a neural operator capable of approximating the target operator over a compact set to any desired degree of accuracy. These theorems further justify the application of neural operators. Multiple Universal Approximation Theorems have been established for various kinds of neural operators~\cite{kovachki2021neural}, including the FNO~\cite{kovachki2021universalapproximationerrorbounds, lanthaler2024nonlocalitynonlinearityimpliesuniversality}, as well as for the DeepONet~\cite{lanthaler2022errorestimatesdeeponetsdeep} and PCA-Net~\cite{lanthaler2023operatorlearningpcanetupper}.

Analogous to the case of finite-dimensional neural networks, neural operators suffer from the curse of parametric complexity, the infinite-dimensional analogue to the curse of dimensionality. This curse is an exponential growth in size with respect to an increasing accuracy, which is also present when approximating $C^r$- or Lipschitz-regular operators as discussed in~\cite{lanthaler2025parametriccomplexityoperatorlearning}. In~\cite{kovachki2021universalapproximationerrorbounds}, they overcome this curse when approximating the solution operator of the Darcy-type elliptic equation and the incompressible Navier-Stokes equations. The idea of the proof relies on recasting the underlying PDEs as fixed point problems which then allows for utilizing the Picard iteration. The latter is then realized by the layers of the investigated neural operator. The authors of 
~\cite{furuya2024quantitativeapproximationneuraloperators,furuya2026approximationtheorylaplacianbasedneural, kow2026quantitativeapproximationtheoremneural} follow the same idea, leading to multiple neural operators that approximate solution operators of various PDEs that do not suffer from the curse of parametric complexity.

Although much encouraging work has been done to couple neural operators with ODEs~\cite{Lu_2021,jin2022mionetlearningmultipleinputoperators,behroozi2025sensitivityconstrainedfourierneuraloperators, liu2026inverseneuraloperatorode,zhou2025parameterestimationstructuraldynamics}, no effort has yet been made to derive quantitative approximation rates for approximating the solution operator of a controlled ODE with a Fourier neural operator. Since neural operators usually act on function spaces, the object of interest in the case of ODEs is the setting when so-called control functions are present, which then act as the input functions for the neural operator. We therefore refer to such ODEs as controlled ordinary differential equations (CODEs) and highlight~\cite{book_sontag, decarvalhoneto2025generalversioncaratheodorysexistence} for further reading on this topic. In this article, we derive error bounds for the Fourier neural operator when approximating the solution operator of CODEs. This is motivated by the fact that the differentiation formula of the Fourier transform suggests that the emulation of a CODE, when recast as a fixed point integral formula, is realizable through the layers of a FNO.

Such controlled systems can arise in practice, e.g., when controlling and modeling the process of electric arc furnaces~\cite{identification_electric_arc_furnace, MeierThomas}. To obtain meaningful predictions, simulations of these furnaces require inverse parameter estimation~\cite{RATHABA200389,MarulandaDurango,Lpez2014CALIBRATIONOP}. As electric arc furnaces are run under different control settings, their operation generates large data sets that can be used for parameter estimation. However, the long inference times of complex physical simulations can hinder a successful inversion. We therefore propose to apply the Ensemble Kalman Inversion (EKI) where one additionally replaces the simulation by a FNO surrogate. Interestingly the authors of~\cite{kow2026quantitativeapproximationtheoremneural,iglesias2025deeponetacceleratedbayesianinversionmoving} share the same idea of deploying the EKI with a learned neural operator. However, we extend this approach to a Kaczmarz structure, which is commonly encountered in deterministic inverse problems~\cite{Leit_o_2016,Kowar2002ConvergenceAO,NanYu2026InverseProblemsandImaging,Li_2018}. For further reading on EKI we refer the reader to~\cite{Iglesias_2013,Iglesias_2016,iglesias2020adaptiveregularisationensemblekalman} and for inverse problems we suggest~\cite{clason2021regularizationinverseproblems,Tarantola2004InversePT,Rieder2003TheAI,rothermel2019parameterestimationtemperaturedependent,Rothermel2020SolvingAI,klein2019sequentialsubspaceoptimizationrecovering}.
\subsection{Our results and contributions}
After introducing the CODE system in a Carathéodory setting, we give appropriate assumptions that lead to a well-defined operator acting between the spaces of Lebesgue integrable and continuous functions. The corresponding result is summarized in~\cref{my_caratheodory}. We then strengthen these assumptions to establish the setting for deriving error bounds when deploying the Fourier neural operator as a surrogate of the solution operator. Our main result~\cref{main theorem} states that the size of the FNO scales at most algebraically and thus, exponential growth in the size of the FNO with respect to an increasing accuracy can be avoided. Under a sufficiently smooth $f$ the size in fact grows at most sub-quadratically up to logarithmic factors. Proving this relies on a reformulation of the CODE as a fixed point integral formula, while utilizing numerous results from Fourier analysis. The non-linearity $f$ in the CODE system can be efficiently approximated by a finite neural network. Consequently, the idea of our proof is similar to that of~\cite{furuya2024quantitativeapproximationneuraloperators}, implying that other problems with a comparable structure can be resolved in the same manner. At the end of the paper, we introduce an inverse parameter estimation problem in which a set of measurement data from multiple experiments is present. We propose solving this problem through surrogate modeling and applying the Ensemble Kalman Inversion in a Kaczmarz-type fashion. 
\subsection{Notation}
Let us briefly introduce the notation used throughout this paper.
\begin{itemize}
    \item We set $\mathbb{N}=\{1,2,\dots\}$ and $\mathbb{N}_0=\mathbb{N}\cup\{0\}$. Further we set $\mathbb{R}^+=(0,\infty)$.
    \item Let $\mathbb{T}^d=[0,2\pi]^d, d\in\mathbb{N}$, be the periodic domain which we identify with the standard torus.
    \item By $\Vert\cdot\Vert_{\mathbb{R}^n}, n\in\mathbb{N}$, we denote the euclidean norm.
    \item For $q\in[1,\infty]$ the Lebesgue space $L^q(I;\mathbb{R}^{n})$ consists of all measurable functions $u:I\rightarrow\mathbb{R}^n$ such that
    \begin{equation*}
        \Vert u\Vert_{L^q(I;\mathbb{R}^n)} := \left\{ \begin{aligned}
            &{\big( \int_I\Vert u(x)\Vert_{\mathbb{R}^{n}}^qdx\big)}^{1/q}, \ \ \text{if} \ \ q<\infty,\\
                                            &\esssup_{x\in I}\Vert u(x)\Vert_{\mathbb{R}^n}, \ \ \text{if} \ \ q=\infty.    
        \end{aligned}
         \right.
    \end{equation*}
    \item By $C_\textnormal{abs}(I;\mathbb{R}^n)$ we denote the space of absolutely continuous functions endowed with the sup-norm. For $0\leq s\leq\infty$ we denote by $C^s(I;\mathbb{R}^n)$ the space of $s$-times continuously differentiable functions.
    \item For $s\in\mathbb{N}_0, p\in [1, \infty]$ the Sobolev space $W^{s,p}(I;\mathbb{R}^n)$ is defined by
    \begin{align*}
        W^{s,p}(I;\mathbb{R}^n):=\{u\in L^p(I;\mathbb{R}^n) \ | \ D^\alpha u\in L^p(I:\mathbb{R}^n) \ \ \text{for all} \ \ \alpha\in\mathbb{N}^d_0 \ \ \text{with} \ \ \vert\alpha\vert\leq s \},
    \end{align*}
    where $D^\alpha u$ denotes the weak derivative of order $s$ and $\alpha$ is a multi-index such that 
    \begin{align*}
    \alpha&=(\alpha_1,\dots,\alpha_d),\
    D^\alpha u = (\frac{\partial^{\alpha_1}}{\partial^{\alpha_1}x_1}\dots\frac{\partial^{\alpha_d}}{\partial^{\alpha_d}x_d})u,\
        \vert\alpha\vert=\sum_{i=1}^d\alpha_i.
    \end{align*}
    The $W^{s,p}(I;\mathbb{R}^n)$-norm is given by
    \begin{equation*}
        \Vert u\Vert_{W^{s,p}(I;\mathbb{R}^n)} := \left\{ \begin{aligned}
            &{\big( \sum_{0\leq\vert\alpha\vert\leq s}\Vert D^\alpha u\Vert^p_{L^p(I;\mathbb{R}^n)}\big)}^{1/p}, \ \ \text{if} \ \ p<\infty,\\
                                            &\max_{0\leq\vert\alpha\vert\leq s}\Vert D^\alpha u\Vert_{L^\infty(I;\mathbb{R}^n)}, \ \ \text{if} \ \ p=\infty.    
        \end{aligned}
         \right.
        \end{equation*}
    We define the notation $H^s(I;\mathbb{R}^n):=W^{s,2}(I;\mathbb{R}^n)$.
    \item In general we write $\Vert f\Vert_{\mathcal{X}(I)}$ for functions $f\in\mathcal{X(I;\mathbb{R})}$. Thus, we also just write $\mathcal{X}(I)$ instead of $\mathcal{X}(I;\mathbb{R})$.
\end{itemize}
All proofs of our results are found in the appendix of the respective chapters.
\section{Well-posedness for controlled ordinary differential equations (CODE)}
\label{Sec:WPCODE}
We are interested in systems of controlled ordinary differential equations, that is
\begin{equation}
\label{eq: ODE constant controls}
\begin{aligned}
\dot{y}(t) &= f(t, y(t), u(t)), \ t\in I:=(t_1, t_2)\subset\mathbb{R},\ -\infty< t_1 < t_2 <\infty,\\
y(t_1) &= y^0,
\end{aligned}
\end{equation}
where $f:I\times\mathbb{R}^{n_y}\times\mathbb{R}^{n_u}\rightarrow\mathbb{R}^{n_y}, y:I\rightarrow\mathbb{R}^{n_y}, u:I\rightarrow\mathbb{R}^{n_u},\ n_y, n_u\in\mathbb{N}$. Further we assume that the \textit{control} function $u$ does not necessarily need to be a continuous but a Lebesgue integrable function. Such setting is usually referred to as Carathéodory setting in literature. See~\cite{book_sontag},~\cite{decarvalhoneto2025generalversioncaratheodorysexistence},~\cite{roubicek},~\cite{o2013existence} and~\cite{Beldzinski} for the study of controlled ODEs. We introduce standard assumptions, under which the system is \textit{well posed}. By well-posedness we mean that for every control function $u$ and initial value $y^0$ there exists a unique solution $y$ of~\cref{eq: ODE constant controls}. By slightly extending the assumptions we obtain continuity of the solution operator, i.e. the mapping
\begin{align*}
    (u,y^0)\mapsto y
\end{align*}
is continuous with respect to appropriately chosen function spaces and their respective norms.
Note that the derivative and integration are understood component-wise throughout this paper.

\begin{thm}[Well-posedness of CODE]\label{my_caratheodory}
    Let $n_y,n_u\in\mathbb{N}$ and $I\subset\mathbb{R}$ bounded and consider $f=f(t, y(t), u(t)):I\times\mathbb{R}^{n_y}\times\mathbb{R}^{n_u}\rightarrow\mathbb{R}^{n_y}$. Assume the following properties, which we refer to as Carathéodory properties, 
    \begin{enumerate}
        \item $u\in L^q(I;\mathbb{R}^{n_u}), 1\leq q\leq\infty$,
        \item $f(t,\cdot,\cdot)$ is continuous for almost all $t\in I$,
        \item $f(\cdot, y, u)$ is measurable for all $y\in\mathbb{R}^{n_y}, u\in\mathbb{R}^{n_u}$,
        \item there exists $\alpha\in L^p(I;[0,\infty)), 1\leq p\leq\infty$ with $\Vert f(t,y,u)\Vert_{\mathbb{R}^{n_y}}\leq\alpha(t)$ for almost all $t\in I$, all $y\in\mathbb{R}^{n_y}$ and $u\in\mathbb{R}^{n_u}$,
        \item there exists $m\in L^1(I;[0,\infty))$ with $\Vert f(t, y, u)-f(t, z, u)\Vert_{\mathbb{R}^{n_y}}\leq m(t)\Vert y-z\Vert_{\mathbb{R}^{n_y}}$ for almost every $t\in I$ and all $y,z\in\mathbb{R}^{n_y}, u\in\mathbb{R}^{n_u}$,
        \item for every fixed $t_1\in I$ and $y\in C(I;\mathbb{R}^{n_y})$ the mapping $u\mapsto\big(t\mapsto\int^t_{t_1}f(s, y(s), u(s)) ds\big)$ is continuous from $L^q(I;\mathbb{R}^{n_u})$ to $C(I;\mathbb{R}^{n_y})$.
    \end{enumerate}
Then, for every $t_1\in I, y^0\in\mathbb{R}^{n_y}$ and $u\in L^q(I;\mathbb{R}^{n_u})$, there exists a unique function $y\in C_\text{abs}(I;\mathbb{R}^{n_y})$ such that
    \begin{align*}
        \dot{y}(t)&=f(t,y(t),u(t)) \ \text{almost everywhere in} \ I,\\
        y(t_1) &= y^0,
    \end{align*}
so that $y$ is a solution of the controlled ODE. Further the solution $y$ depends continuously on $y^0$ and $u$, that is, continuous in the sense of
\begin{align*}
    \mathbb{R}^{n_y}\times L^q(I;\mathbb{R}^{n_u})&\rightarrow C(I;\mathbb{R}^{n_y}),\ (y^0, u)\mapsto y.
\end{align*}
To highlight the dependence on initial value $y^0$ and control $u$ we write $y(t;u,y^0)$.
\end{thm}

To obtain error rates for approximating above mapping with a Fourier neural operator, we need to introduce more structure regarding the regularity of admissible controls. We define following \textit{control set}
\begin{align*}
    \mathcal{U}^{l,I,A}:=\{u\in H^{l}(I;\mathbb{R}^{n_u}) \ | \ \Vert u\Vert_{H^l(I;\mathbb{R}^{n_u})}\leq A\in\mathbb{R}^+\}, 
\end{align*}
where we for now assume $l\geq 1$. For initial values consider
\begin{align*}
    y^0\in Y^0\subset\mathbb{R}^{n_y},
\end{align*}
where $Y^0\subset\mathbb{R}^{n_y}$ is bounded. These assumptions immediately lead to bounded CODE solutions under the Carathéodory properties.
\begin{lemma}\label{bound of original ode solution}
    Consider $f$ from~\cref{my_caratheodory} with the assumptions (2)-(6) and $u\in \mathcal{U}^{l,I,A}, l\geq 1, y^0\in Y^0\subset\mathbb{R}^{n_y}$ additionally. Then we find a $Y\in\mathbb{R}^+$ such
    that for the solution $y(t; u, y^0)\in C_\textnormal{abs}(I;\mathbb{R}^{n_y})$ of the CODE system it holds
    \begin{align*}
        y(t;y^0,u) \in {(-Y,Y)}^{n_y}
    \end{align*}
    for all $t\in I$.
\end{lemma}
When working with FNOs, we usually study periodic functions on the domain $\mathbb{T}^d=[0,2\pi]^d, d\in\mathbb{N}$. Therefore, to make Fourier analysis in this setting possible, we focus on solutions that at least live on a smaller time domain 
\begin{align*}
    J:=(\hat{t}_1,\hat{t}_2)\subset (0,2\pi)=\mathbb{T},\ \hat{t}_1<\hat{t}_2.
\end{align*}
However, for our CODE system, we still assume the system to be defined for an interval $I=(t_1,t_2)$ with $\mathbb{T}\subset I$, while the Carathéodory properties are assumed to hold on $I$. The reason for such interval $I$ is discussed at the end of~\cref{sec: Proof of well-posedness of CODE}.

By definition, the control functions do not need to be periodic but the following lemma shows that we are able to consider extended control functions such that they are periodic. The lemma can be found in the appendix of~\cite{kovachki2021universalapproximationerrorbounds,lanthaler2024nonlocalitynonlinearityimpliesuniversality} where it is stated that the lemma relies on a result by~\cite{Stein}.
\begin{lemma}[\cite{kovachki2021universalapproximationerrorbounds,lanthaler2024nonlocalitynonlinearityimpliesuniversality,Stein}]\label{periodic extension operator}
    Let $D\subset\mathbb{R}^d, d\in\mathbb{N}$ be a bounded, Lipschitz domain~\cite{ern2016mollificationstronglylipschitzdomains}. There exists a continuous, linear operator $\mathcal{E}:W^{s,p}(D)\rightarrow W^{s,p}_\text{per}(F)$ for any $s\geq 0$ and $1\leq p\leq\infty$
    where $F\subset\mathbb{R}^d$ is a bounded hypercube with $D\Subset F$ such that, for any $u\in W^{s,p}(D)$,
    \begin{itemize}
        \item $\mathcal{E}(u)|_D = u$,
        \item $\mathcal{E}(u)\in W^{s,p}_\text{per}(F)$ is periodic on $F$ (including its derivatives).
    \end{itemize}
Furthermore, $\mathcal{E}$ maps continuously differentiable functions to continuously differentiable functions, so
\begin{align*}
    \mathcal{E}(C^s(\overline{D}))\subset C^s_\text{per}(F),
\end{align*}
and hence defines a continuous mapping $\mathcal{E}:C^s(\overline{D})\rightarrow C^s_\text{per}(F)$.
\end{lemma}
We fix this extension operator $\mathcal{E}$ and introduce following control set to keep notation simple
\begin{align}\label{set of periodically extended controls}
    \mathcal{U}^{l,\mathbb{T},A}=\{\mathcal{E}{u}\in H^{l}(\mathbb{T};\mathbb{R}^{n_u}) \ | \ u\in\mathcal{U}^{l,J,A}\},
\end{align}
where we now assume to have at least $l\geq 2$ regarding the regularity. The controls $u\in\mathcal{U}^{l,\mathbb{T},A}$ are the periodic versions of the control functions from $\mathcal{U}^{l,J,A}$.
We sometime omit to explicitly write $\mathcal{E}u$ and by abuse of notation just write $u\in\mathcal{U}^{l,\mathbb{T},A}$ when it is obvious. Further we define the periodically extended time variable as 
\begin{align*}
    \tau:\mathbb{T}&\rightarrow\mathbb{R},\ t\mapsto\mathcal{E}\text{id}(t).
\end{align*}

One of the reasons we introduce the periodic time variation $\tau$ is displayed in the following. Ignoring the inner function $y$, the function
\begin{align*}
    t\mapsto f(t,y(t),u(t)),\ t\in\mathbb{T}
\end{align*}
cannot be periodic alone due to the $t$-dependence. Thus we need to introduce such periodic time variable. 
As $\tau$ is dependent on $t$ and constructed the same way as the controls $u\in\mathcal{U}^{l,\mathbb{T},A}$, we may consider $\tau$ as an additional control function. 
Therefore to simplify analysis and notation, we study for now the autonomous system
\begin{equation}
\label{eq: autonomous ODE}
\left.
\begin{aligned}
    \dot{y}(t) &= f(y(t),u(t)),\ t\in\mathbb{T},\\
    y(\hat{t}_1) &= y^0,
\end{aligned}
\right\}
\end{equation}
where $u\in\mathcal{U}^{l,\mathbb{T}, A}, l\geq 2$ but with $\textnormal{Im}(u)\subset\mathbb{R}^{1+n_u}$. Looking at the Carathéodory assumptions, this means we require
\begin{align*}
    t\mapsto f(t,y,u)
\end{align*}
to be continuous instead of measurable for every $y\in\mathbb{R}^{n_y}, u\in\mathbb{R}^{n_u}$.

The functions $u$ from the introduced control set $\mathcal{U}^{l,\mathbb{T},A}$ share the same bound with respect to the sup-norm. Since we discussed that the solutions are also bounded, this will allow us to study the non-linearity $f$ on a bounded domain, which is essential if one is interested in approximating $f$ by finite neural networks.
\begin{lemma}\label{bound of periodicially extended controls}
For the control functions $u\in\mathcal{U}^{l,\mathbb{T},A}, l\geq 2$, with $\textnormal{Im}(u)\subset\mathbb{R}^{1+n_u}$ there exists a $U^\prime\in\mathbb{R}^+$ such that 
\begin{align*}
    u(t)\in(-U^\prime, U^\prime)^{1+n_u}
\end{align*}
for all $t\in\mathbb{T}$.
\end{lemma}
Later, when reformulating CODE as Picard iteration, we have to recast it in a way such that it is applicable for
emulation through the layers of a FNO. This leads us to replace $f$ with an approximation by finite neural networks. An approximation like that eventually introduces an accuracy error inside the Picard iteration. Thus we have to enlarge the domain of interest.

Choose a $\gamma\in\mathbb{R}^+$ such that
\begin{equation}\label{choosage of gamma}
\begin{aligned}
    Y':=Y+\gamma>1
\end{aligned}
\end{equation}
and define the domain
\begin{align*}
D^\prime:={(-Y^\prime-1,Y^\prime+1)}^{n_y}\times{(-U^\prime-1,U^\prime+1)}^{1+n_u}
\end{align*}
to study $f$ on. The margin $1$ ensures that the closed boxes ${[-Y^\prime,Y^\prime]}^{n_y}$ and ${[-U^\prime,U^\prime]}^{1+n_u}$, which contain the values of all functions inserted into $f$ below, lie in the interior of $D^\prime$.
Before introducing Fourier neural operators we summarize all final assumptions on the CODE system.
\begin{assump}[Final assumptions of CODE system]\label{Final assumptions of ODE system}
    \begin{equation*}
\left\{
    \begin{aligned}
        &J := (\hat{t}_1,\hat{t}_2)\subset(0,2\pi),\\
        &y^0\in Y^0,\\
        &u\in\mathcal{U}^{l,J,A}, l\geq 2,\\
        &f \ \text{as in} \ \cref{my_caratheodory} \ \text{with}\
        \alpha\in L^\infty(J;[0,\infty)),\\
        &s \geq 2,\\
        &D^\beta f{|_{D^\prime}} \in L^\infty(D^\prime;\mathbb{R}^{n_y}),\
        \forall\beta\in\mathbb{N}^d_0\ (d:=1+n_y+n_u),
        \vert\beta\vert\leq s,\\
        &\Vert f\Vert_{W^{s,\infty}(D^\prime;\mathbb{R}^{n_y})}\leq E\in\mathbb{R}^+.
    \end{aligned}
\right.
\end{equation*}
\end{assump}
Note that we assume that $f\in W^{s,\infty}, s\geq 2$ with $\Vert f\Vert_{W^{s,\infty}(D^\prime;\mathbb{R}^{n_y})}\leq E\in\mathbb{R}^+$ which is needed for the efficient approximation by a finite neural network, see~\cref{finite neural network approx for my case,quantitative approximation with GELU}.
The assumption $\Vert f\Vert_{W^{s,\infty}(D^\prime;\mathbb{R}^{n_y})}\leq E\in\mathbb{R}^+$ is a \textit{smallness condition} where we need to choose $E$ as discussed in~\cref{choosage of E} to secure convergence of our proposed Picard type iteration.
\section{Fourier neural operator for CODE}
\label{Sec:FNOCODE}
Fourier neural operators, first introduced in~\cite{li2021fourierneuraloperatorparametric}, have strong approximation
capabilities~\cite{kovachki2021universalapproximationerrorbounds, lanthaler2024nonlocalitynonlinearityimpliesuniversality} while offering to apply methods from Fourier analysis including the differentiation formula. The latter formula allows us to rewrite an integral formula using the (inverse) Fourier transform, which makes FNOs an encouraging candidate for the emulation of the Picard iteration. Thus we begin by introducing all Fourier analysis related notations and definitions needed for understanding our main result~\cref{main theorem}. Most of those notations and definitions are from the paper~\cite{kovachki2021universalapproximationerrorbounds} as our article is intended to build upon.

\subsection{Definition}
We start by introducing neural operator architecture in a broader sense.
Let $d\in\mathbb{N}$ and denote by $D$ a domain in $\mathbb{R}^d$. 
Consider arbitrary Banach spaces of form $\mathcal{A}(D;\mathbb{R}^{n_a})$ and $\mathcal{U}(D;\mathbb{R}^{n_u})$ with $d_a,d_u\in\mathbb{N}$.
A general neural operator is a mapping between such spaces denoted by
\begin{align*}
    \mathcal{N}:\mathcal{A}(D;\mathbb{R}^{n_a})&\rightarrow\mathcal{U}(D;\mathbb{R}^{n_u}),\
    a \mapsto \mathcal{N}(a) = \mathcal{Q}\circ\mathcal{L}_L\circ\mathcal{L}_{L-1}\circ\dots\circ\mathcal{L}_1\circ\mathcal{R}(a)
\end{align*} 
for a given \textit{depth} $L\in\mathbb{N}$ where $\mathcal{R}:\mathcal{A}(D;\mathbb{R}^{n_a})\rightarrow\mathcal{U}(D;\mathbb{R}^{n_v}), d_v\geq d_u$ is a \textit{lifting} operator (acting locally) defined by
\begin{equation}
    \label{lifting_layer_no_positional_encoding}
\begin{aligned}
\mathcal{R}(a)(x)=Ra(x), \ \ R\in\mathbb{R}^{d_v\times d_a}, 
\end{aligned}
\end{equation}
and $\mathcal{Q}:\mathcal{U}(D;\mathbb{R}^{d_v})\rightarrow\mathcal{U}(D;\mathbb{R}^{d_u})$ is a local \textit{projection} operator of form
\begin{equation}
    \label{projecting_layer_no_positional_encoding}
\begin{aligned}
    \mathcal{Q}(v)(x)=Qv(x), \ \ Q\in\mathbb{R}^{d_u\times d_v}.
\end{aligned}
\end{equation}
The non-linear operator layers $\mathcal{L}_l:\mathcal{U}(D;\mathbb{R}^{n_v})\rightarrow \mathcal{U}(D;\mathbb{R}^{n_v}), v\mapsto \mathcal{L}_l(v)$ have following structure
\begin{align*}
    \mathcal{L}_l(v)(x)=\sigma(W_l v(x) + b_l(x) + (\mathcal{K}(a; \theta_l) v)(x)), \ \ \forall x\in D.
\end{align*}
The \textit{weight} matrix $W_l\in\mathbb{R}^{n_v\times n_v}$ and the \textit{bias} $b_l(x)\in\mathcal{U}(D;\mathbb{R}^{n_v})$ define affine point-wise mappings $W_lv(x)+b_l(x), l=1,\dots,L$, while the linear operators $\mathcal{K}(a,\theta)$, acting non-local, are integral operators of form
\begin{align*}
(\mathcal{K}(a;\theta)v)(x) = \int_D \kappa_\theta (x,y;a(x),a(y)) v(y)dy, \ \ \forall x\in D,
\end{align*}
where the integral kernel $\kappa_\theta:\mathbb{R}^{2(d+d_a)}\rightarrow\mathbb{R}^{d_v\times d_v}$ is a neural network parametrized by $\theta\in\Theta$. The non-linear \textit{activation function} $\sigma$ acts component-wise.
Fourier neural operators are a special case of these general neural operators where the kernel is of form $\kappa_\theta =\kappa_\theta(x-y)$. In this case our linear operators are convolutions
\begin{align*}
    (\mathcal{K}(\theta)v)(x) = \int_D\kappa_\theta(x-y) v(y)dy, \ \ \forall x\in D.
\end{align*}
The domain is now considered as $D={[0,2\pi]}^d$. Concluding, FNOs are concatenations of
\begin{align*}
    \mathcal{N}(a) = \mathcal{Q}\circ\mathcal{L}_L\circ\mathcal{L}_{L-1}\circ\dots\circ\mathcal{L}_1\circ\mathcal{R}(a),
\end{align*}
where 
\begin{equation}
     \label{nonlinear FNO layer}
\begin{aligned}
    \mathcal{L}_l(v)(x) = \sigma\biggl(W_l v(x) + b_l(x) + \mathcal{F}^{-1}\big(P_l(k)\cdot \mathcal{F}(v)(k)\big)(x)\biggr),
\end{aligned}
\end{equation}
while $\mathcal{F}\ (\mathcal{F}^{-1})$ denotes the (inverse) Fourier transform. The \textit{Fourier multiplier} matrix $P_l:\mathbb{Z}^d\rightarrow\mathbb{C}^{d_v\times d_v}$ defines the coefficients of a non-local, linear mapping via the Fourier transform. 
It stems from the integral kernel but we consider it directly parametrized in Fourier space as we have $P_\theta(k)=\mathcal{F}(\kappa_\theta)(k)$.
We impose that $P_\theta(-k)=\overline{P_\theta(k)}$, where the bar denotes entry-wise complex conjugation, to ensure that $(\mathcal{K}(\theta)v)(x)$ is a real-valued function for a real-valued input function $v$.

As we are interested in deriving approximation rates for the FNO, we need to quantify the size, i.e. the amount of trainable parameters of such architecture. Note that the layers of a FNO require to compute an integral exactly to perform the (inverse) Fourier transform. Thus, an implementable version of the FNO, the so-called $\Psi$-FNO, is introduced by utilizing the discrete (inverse) Fourier transform. For this, we need the Fourier wavenumbers first.
\begin{defn}[Fourier wavenumbers]
A set of Fourier wavenumbers $k\in\mathcal{K}_N$ is given by
\begin{align*}
    \mathcal{K}_N:=\{k\in\mathbb{Z}^d \ | \ \vert k\vert_\infty\leq N\},
\end{align*}
where $k=(k_1,\dots,k_d)\in\mathbb{Z}^d$ and
    \begin{align*}
        \vert k\vert_\infty:=\max_{i=1,\dots,d}\vert k_i\vert\leq N.
    \end{align*}
\end{defn}
In the following we introduce a uniform discretized grid of $\mathbb{T}^d=[0,2\pi]^d$.
Set the index set as
\begin{align*}
\mathcal{J}_N:={\{0,\dots,2N\}}^d,
\end{align*}
and the grid points with the corresponding function evaluations as
\begin{align*}
    t_j &:= \frac{2\pi j}{2N+1},\
    v_j := v(t_j), \ j\in\mathcal{J_N}.
\end{align*}
\begin{defn}[Discrete Fourier transform]
The discrete Fourier transform $\mathcal{F}_N:\mathbb{R}^{\mathcal{J}_N}\rightarrow\mathbb{C}^{\mathcal{K}_N}$ is defined by
\begin{align*}
    \mathcal{F}_N(v)(k):=\frac{1}{{(2N+1)}^d}\sum_{j\in\mathcal{J}_N}v_j e^{-2\pi i\langle j,k\rangle/{(2N+1)}}
\end{align*}
and the inverse $\mathcal{F}_N^{-1}:\mathbb{C}^{\mathcal{K}_N}\rightarrow\mathbb{R}^{\mathcal{J}_N}$,
\begin{align*}
    \mathcal{F}_N^{-1}(\hat{v})(j):=\sum_{k\in\mathcal{K}_N}\hat{v}_k e^{2\pi i\langle j,k\rangle /{(2N+1)}}.
\end{align*}
\end{defn}
When working with the discrete (inverse) Fourier transform, a natural function space to study are the trigonometric polynomials of corresponding degree $N\in\mathbb{N}$. Thus, we introduce the precise definition.
\begin{defn}[Trigonometric polynomials]
    Let $N\in\mathbb{N}$, then denote by $L_N^2(\mathbb{T}^d)$ the space of trigonometric polynomials $v_N:\mathbb{T}^d\rightarrow\mathbb{R}$ defined by
    \begin{align*}
        v_N(t)=\sum_{\vert k\vert_\infty\leq N} c_k e^{i\langle t, k\rangle},
    \end{align*}
    with $c_k\in\mathbb{C}$.
\end{defn}
We recall that a trigonometric polynomial is uniquely defined by its values on the grid $\mathcal{J}_N.$
With the following projection, the so-called pseudo-spectral projection, we are able to to project continuous function, and thus Sobolev function under certain assumptions, on the space of trigonometric polynomials.
\begin{defn}[Pseudo-spectral projection]
    The pseudo-spectral projection is given by 
    \begin{align*}
    \mathcal{I}_N: C(\mathbb{T}^d)\rightarrow L^2_N(\mathbb{T}^d),
\end{align*}
where $\mathcal{I}_N v$ of a continuous function $v$ is defined as the unique trigonometric polynomial in $L^2_N(\mathbb{T}^d)$ such that
\begin{align*}
    \mathcal{I}_N v(t_j) = v(t_j), \ \forall j\in\mathcal{J}_N.
\end{align*}
\end{defn}
This projection allows us, to give an analogue definition of the $\Psi$-FNO compared to the continuous FNO.
\begin{defn}[$\Psi$-FNO]
    A $\Psi$-FNO is a mapping 
    \begin{align*}
        \mathcal{N}:\mathcal{A}(\mathbb{T}^d;\mathbb{R}^{d_a})\rightarrow\mathcal{U}(\mathbb{T}^d;\mathbb{R}^{d_u}), \ a\mapsto\mathcal{N}(a)
    \end{align*}
    of form
    \begin{align*}
        \mathcal{N}(a)=\mathcal{Q}\circ\mathcal{I}_N\circ\mathcal{L}_L\circ\mathcal{I}_N\circ\dots\circ\mathcal{L}_1\circ\mathcal{I}_N\circ\mathcal{R}(a),
    \end{align*}
    where $\mathcal{I}_N$ is the pseudo-spectral Fourier projection. The lifting and projection layer are given again as before, see~\cref{lifting_layer_no_positional_encoding} and~\cref{projecting_layer_no_positional_encoding}.
    The non-linear layers $\mathcal{L}_l, l=1,\dots,L$ are given by~\cref{nonlinear FNO layer}, where we allow the Lipschitz continuous activation function $\sigma=\sigma_l$ to depend on the layer $l$. Lemma D.4 and Lemma D.5 of~\cite{kovachki2021universalapproximationerrorbounds} remain valid in this case, since their proofs only use the Lipschitz continuity of the activation functions.
\end{defn}
A $\Psi$-FNO is uniquely defined as by its restriction on  the finite-dimensional subspace $L_N^2(\mathbb{T}^d;\mathbb{R}^{d_a})\subset\mathcal{A}(\mathbb{T}^d;\mathbb{R}^{d_a})$. Further a $\Psi$-FNO can be naturally identified with a finite-dimensional mapping
\begin{align*}
    \widehat{\mathcal{N}}:\mathbb{R}^{d_a\times \mathcal{J}_N}\rightarrow\mathbb{R}^{d_u\times \mathcal{J}_N}, \ a\mapsto\widehat{\mathcal{N}}(a),
\end{align*}
with input $a={\{a_j\}}_{j\in\mathcal{J}_N}\in\mathbb{R}^{d_a\times\mathcal{J}_N}$ corresponding to the point-values $a_j=a(t_j)$ on the grid ${\{t_j\}}_{j\in\mathcal{J}_N}$ and $\mathcal{J}_N:={\{0,\dots,2N\}}^d$. The non-linear layers $\widehat{\mathcal{L}}_l$ for $l=1,\dots,L$ are
\begin{align*}
    \widehat{\mathcal{L}}{(v)}_j = \sigma\biggl(W_l v_j + b_{l,j} + \mathcal{F}_N^{-1}{\big(P_l(k)\cdot\mathcal{F}_N(v)(k)\big)}_j\biggr)
\end{align*}
for $j\in\mathcal{J}_N$. 
To quantify the approximation error with respect to the network's size, we give following definition.
\begin{defn}[Size of $\Psi$-FNO]
    The depth and width of a $\Psi$-FNO $\mathcal{N}$ are defined by 
    \begin{align*}
        \textnormal{depth}(\mathcal{N}):=L,\ \textnormal{width}(\mathcal{N}):=d_v\vert\mathcal{J}_N\vert=d_v\vert\mathcal{K}_N\vert = {(2N+1)}^{d}d_v,\
        \textnormal{lift}(\mathcal{N}):=d_v.
    \end{align*}
    The precise size of $\Psi$-FNO is then given by 
    \begin{align*}
        \textnormal{size}(\mathcal{N})=\underbrace{d_v d_u}_{\textnormal{size}(\mathcal{Q})} + L\underbrace{(d_v^2 + d_v\vert\mathcal{J}_N\vert + d_v^2\vert\mathcal{J}_N\vert)}_{\textnormal{size}(\mathcal{L}_l)} + \underbrace{d_a d_v}_{\textnormal{size}(\mathcal{R})}.
    \end{align*}
\end{defn}

We use the estimate
\begin{align}\label{estimate for size of psi fno}
    \textnormal{size}(\mathcal{N})\lesssim\textnormal{depth}(\mathcal{N})\textnormal{width}(\mathcal{N})\textnormal{lift}(\mathcal{N}),
\end{align}
where we assume that $\max(d_a,d_u)\leq d_v$ which implies that we have $\textnormal{size}(\mathcal{N})\sim Ld_v^2\vert\mathcal{J}_N\vert$.

Besides the pseudo-spectral projection, we also use the $L^2$-orthogonal projection to project the space of square-integrable functions onto the space of trigonometric polynomials. 
\begin{defn} [$L^2$-orthogonal projection]
    The $L^2$-orthogonal projection
\begin{align*}
    P_N: L^2(\mathbb{T}^d)\rightarrow L^2_N(\mathbb{T}^d)
\end{align*}
is defined by 
\begin{align*}
    P_N(\sum_{k\in \mathbb{Z}^d}c_k e^{i\langle t,k\rangle}) = \sum_{|k|_\infty\leq N} c_k e^{i\langle t,k\rangle}.
\end{align*} 
\end{defn}
Recall, that we need to reformulate the Picard iteration of the CODE in such a way that emulation through the FNO layers, or to be more precise through the $\Psi$-FNO layers, is feasible. Then, to secure well-definedness regarding the reformulated Picard iteration, we need following additional projection.
\begin{defn}[Mean reducing $L^2$-orthogonal projection]
    Define the projection 
    \begin{align*}
    \dot{P}_N:L^2(\mathbb{T}^d)\rightarrow\dot{L}^2_N(\mathbb{T}^d)
    \end{align*}
    by removing the mean, i.e. $\dot{P}_N v=P_N v - \fint_{\mathbb{T}^d}v(t)dt$ or equivalently
    \begin{align*}
        \dot{P}_N(\sum_{k\in\mathbb{Z}^d}c_k e^{i\langle t,k\rangle})=\sum_{0<\vert k\vert_\infty\leq N}c_k e^{i\langle t,k\rangle}, \ \forall {(c_k)}_{k\in\mathbb{Z}^d}\in\ell^2(\mathbb{Z}^d).
    \end{align*}
\end{defn}
See also~\cite{ShenSpectral} and~\cite{LiJiyong} for further error estimates of the $L^2$-orthogonal and pseudo-spectral projection.
We conclude this subsection by briefly defining the Dirichlet kernel as it is needed for the formulation of our main result~\cref{main theorem}.
\begin{defn}[Dirichlet kernel]\label{Dirichlet kernel}
    For $N\in\mathbb{N}$ and $t\in\mathbb{R}$ we denote by
    \begin{align*}
        D_N(t):=\sum_{\vert k\vert\leq N} e^{itk}
    \end{align*}
    the Dirichlet kernel.
\end{defn}
Further additional definitions and statements of $\Psi$-FNOs needed solely for proving~\cref{main theorem} are found in~\cref{sec: Fourier analysis}.
\subsection{Finite neural network approximation}
As mentioned earlier, we need to find a way of dealing with the non-linearity $f$ in the CODE system to derive error bounds. Fortunately, there are useful results which imply under certain assumptions the existence of finite neural networks capable of avoiding the curse of dimensionality when approximating other functions. These results mostly rely on which type of function is approximated and which kind of activation function is assumed inside the neural network's architecture.

Following definition, presenting one of the most known activation functions, is from~\cite{guehring2019errorboundsapproximationsdeep}. We also recommend said reference for a further reading regarding the precise definition
of a neural network as well as regarding the quantities \textit{depth} and \textit{width} describing a finite network's architecture.
\begin{defn}[ReLU activation function]
    The function 
    \begin{align*}
        \sigma:\mathbb{R}&\rightarrow\mathbb{R},\ x\mapsto\max(0, x)
    \end{align*}
    is called ReLU (Rectified Linear Unit) activation function.
\end{defn}

The next lemma is an extended version to the $n$-dimensional output case of Theorem 4.1 found in~\cite{guehring2019errorboundsapproximationsdeep} by applying their Lemma 2.9 and Lemma 2.11, such that the obtained network is a standard neural network, which is a network without having skip connections. They show in their Theorem 4.1 the existence of a constant $c$ and a finite neural network $\varphi^f_\varepsilon$ with one-dimensional output such that the network's layers $L(\varphi^f_\varepsilon)$ are bounded by $L(\varphi^f_\varepsilon)\leq c\cdot\log(\varepsilon^{-s/(s-r)})$, while the amount of total neurons $N(\varphi^f_\varepsilon)$ is bounded by $N(\varphi^f_\varepsilon)\leq c\cdot \varepsilon^{-d/(s-r)}\log(\varepsilon^{-s/(s-r)})$. Such network may allow for skip connections which is not aligning with our setting as we want to apply~\cref{Psi FNO representation of I_N NN} and the $\Psi$-FNO does not allow for such connections. Further, we want to quantify the width of the finite neural network and not the total amount of neurons due to~\cref{estimate for size of psi fno} and~\cref{Psi FNO representation of I_N NN}. In  the proof of their Lemma 2.11 they construct a standard neural network $\varphi^f_{\textnormal{st},\varepsilon}$ such that the output aligns with the output from the network $\varphi^f_\varepsilon$ having skip connections. They show, while having the same amount of layers, that the standard network's neurons per layer $N_l^\textnormal{st}, l=1,\dots L$ are bounded by $N_l^\textnormal{st}\leq 2N(\varphi^f_\varepsilon), l=1,...,L$. The $n$-dimensional output case follows by placing $n$ such networks $\varphi^f_{i,\varepsilon}, i=1,\dots,n$ with possible skip connections in parallel (cf.\ their Lemma 2.9) implying the existence of a single network $\phi^f_{\varepsilon}$ with $n$-dimensional output and 
\begin{itemize}
    \item[(i)] $L(\phi^f_{\varepsilon})\leq \max_{i=1,\dots,n}{L(\varphi^f_{i,\varepsilon})}\leq\max_{i=1,\dots,n}\{c_i\}\log(\varepsilon^{-s/(s-r)})$,
        \item[(ii)] $N(\phi^f_{\varepsilon})\leq \sum_{i=1}^n N(\varphi^f_{i,\varepsilon})-(n-1)d\leq n\cdot\max_{i=1,\dots,n}\{c_i\}\varepsilon^{-d/(s-r)}\log(\varepsilon^{-s/(s-r)})$. 
\end{itemize}
Using their Lemma 2.11 on $\phi^f_{\varepsilon}$ we obtain a standard neural network $\phi^f_{\textnormal{st},\varepsilon}$ with $n$-dimensional output such that
\begin{itemize}
    \item $L(\phi^f_{\textnormal{st},\varepsilon})\leq \max_{i=1,\dots,n}\{c_i\}\log(\varepsilon^{-s/(s-r)})$,
    \item $N_l(\phi^f_{\textnormal{st},\varepsilon})\leq 2N(\phi^f_{\varepsilon})\leq 2\cdot n\cdot\max_{i=1,\dots,n}\{c_i\}\varepsilon^{-d/(s-r)}\log(\varepsilon^{-s/(s-r)})$.
\end{itemize}
Denoting $L(\phi^f_{\textnormal{st},\varepsilon})$ by $\textnormal{depth}(\phi^f_{\textnormal{st},\varepsilon})$ and the maximum amount of neurons in a layer by $\textnormal{width}(\phi^f_{\textnormal{st},\varepsilon})$, we end up with following result.

\begin{lemma}[Quantitative approximation of ReLU neural network]\label{finite neural network approx for my case}
    Let $d,n,s\in\mathbb{N}$ where $s\geq 2$ with $1\leq p\leq\infty$, $E>0$, and $0\leq r\leq 1$. Then, there exists a constant $c_\textnormal{ReLU}=c_\textnormal{ReLU}(d,n,s,p,E,r)>0$ with the following properties:

    For any $\varepsilon\in(0,\sfrac{1}{2})$ there is a standard neural network architecture 
    with $d$-dimensional input and 
    $n$-dimensional output such that for any $f\in\{g\in W^{s,p}({(0,1)}^d):\Vert g\Vert_{W^{s,p}({(0,1)}^d)}\leq E\}$, there is a standard neural network $\phi^f_{\textnormal{st},\varepsilon}$ that has said architecture such that
    \begin{align*}
        \Vert \phi^f_{\textnormal{st},\varepsilon} - f\Vert_{W^{r,p}({(0,1)}^d;\mathbb{R}^n)}\leq\varepsilon,
    \end{align*}
    where $\sigma$ is the ReLU activation function and
    \begin{itemize}
        \item[(i)] $\textnormal{depth}(\phi^f_{\textnormal{st},\varepsilon})\leq c_\textnormal{ReLU}\cdot\log(\varepsilon^{-s/(s-r)})$,
        \item[(ii)] $\textnormal{width}(\phi^f_{\textnormal{st},\varepsilon}) \leq c_\textnormal{ReLU}\cdot\varepsilon^{-d/(s-r)}\cdot\log(\varepsilon^{-s/(s-r)})$. 
    \end{itemize}
\end{lemma}
To conclude, finite neural networks using the ReLU activation function have under certain assumptions quantitative approximation results, which avoid the curse of dimensionality. However, in practice one commonly uses ``smoother'' versions of the ReLU activation. Thus, we consider an additional activation function to study how the approximation rates are
affected. The activation function that we also consider is the Gaussian Error Linear Unit, or GELU for short, which was first introduced in~\cite{hendrycks2023gaussianerrorlinearunits} and is given here in its standard form as it can also be found in~\cite{yang2026deepneuralnetworksgeneral}.
\begin{defn}[GELU activation function]\label{GELU activation function}
    We define 
    \begin{align*}
        \textnormal{GELU}(x):=\frac{x}{\sqrt{2\pi}}\int^x_{-\infty} e^{-{t^2}/{2}}dt.
    \end{align*}
\end{defn}
The GELU activation function is smooth, meaning it is a $C^\infty$-function and therefore a standard neural network using the GELU function is also smooth.
Analogue to the quantitative approximation statement using the ReLU function, we can find in the literature a similar result. Following theorem is from~\cite{yang2026deepneuralnetworksgeneral} where it is found as Theorem 3\@.
\begin{lemma}[Quantitative approximation of GELU neural networks~\cite{yang2026deepneuralnetworksgeneral}]\label{quantitative approximation with GELU}
    Let $d\in\mathbb{N}$. For any $f\in W^{s,\infty}({[0,1]}^d)$ with $0\leq r<s$, and for any $M,L\in\mathbb{N}$, there
    exists a $GELU$-neural network $\phi$ with depth $c_{1,\textnormal{GELU}}L\log L$ and width $c_{2,\textnormal{GELU}}M \log M$ such that
    \begin{align*}
        \Vert f - \phi\Vert_{W^{r,\infty}({[0,1]}^d)} \leq c_{0,\textnormal{GELU}} \Vert f\Vert_{W^{s,\infty}({[0,1]}^d)}M^{-2(s-r)/d}L^{-2(s-r)/d},
    \end{align*}
where the constants $c_{0,\textnormal{GELU}},c_{1,\textnormal{GELU}},c_{2,\textnormal{GELU}} > 0$ are independent of $M$ and $L$.
\end{lemma}
To apply the last lemmata on our case of interest, consider $f$ from~\cref{Final assumptions of ODE system} with $d=1+n_y+n_u$ and specify $r=1$ in~\cref{finite neural network approx for my case} while choosing a $\varepsilon_\textnormal{ReLU}>0$ according to~\cref{choosage of epsilon ReLU}.
Thus, we obtain the existence of a standard ReLU neural network $f_\textnormal{ReLU}:\mathbb{R}^d\rightarrow\mathbb{R}^{n_y}$ with
\begin{align*}
    \Vert f_\textnormal{ReLU}-f\Vert_{L^\infty(D^\prime;\mathbb{R}^{n_y})}\leq\Vert f_\textnormal{ReLU}-f\Vert_{W^{1,\infty}(D^\prime;\mathbb{R}^{n_y})}\leq\varepsilon_\textnormal{ReLU},
\end{align*}
and 
\begin{itemize}
    \item[(i)] $\textnormal{depth}(f_\textnormal{ReLU})\leq c_\textnormal{ReLU}\cdot\log(\varepsilon_\textnormal{ReLU}^{-s/(s-1)})$,
    \item[(ii)] $\textnormal{width}(f_\textnormal{ReLU}) \leq c_\textnormal{ReLU}\cdot\varepsilon_\textnormal{ReLU}^{-d/(s-1)}\cdot\log(\varepsilon_\textnormal{ReLU}^{-s/(s-1)})$,
\end{itemize}
where the constant $c_\textnormal{ReLU}>0$ is now further dependent on the domain $D^\prime$.
Now for the GELU case, we again choose $r=1$ and apply~\cref{quantitative approximation with GELU}.
To simplify the numerical analysis, we set the width of the GELU network as fixed with value $M_0\in\mathbb{N}$ and study the approximation rates under the influence regarding the amount of layers. We observe that there then exists a constant $c_{0,\textnormal{GELU}}=c_{0,\textnormal{GELU}}(M_0, E)>0$ such that
    \begin{align*}
        \Vert f_\textnormal{GELU}-f\Vert_{L^\infty(D^\prime;\mathbb{R}^{n_y})}\leq\Vert f_\textnormal{GELU} - f\Vert_{W^{1,\infty}(D^\prime;\mathbb{R}^{n_y})} \leq c_{0,\textnormal{GELU}} L^{-2(s-1)/d},
    \end{align*}
where we set $\varepsilon_\textnormal{GELU}:=c_{0,\textnormal{GELU}}L^{-2(s-1)/d}$ and $L\in\mathbb{N}$ chosen as in~\cref{choosage of epsilon GELU} for fixed $s\geq 2$ and $d=1+n_y+n_u$.
In the following we will start to use the notation of 
\begin{align}\label{choosage of fnet}
    f_\textnormal{net}&\in\{f_\textnormal{ReLU}, f_\textnormal{GELU}\},\
    \varepsilon_\textnormal{net}\in\{\varepsilon_\textnormal{ReLU},\varepsilon_\textnormal{GELU}\},
\end{align}
which refers to above obtained ReLU, as well as GELU neural network, for the approximation of $f$. 
\subsection{Approximate Picard iteration of~\texorpdfstring{\cref{eq: autonomous ODE}}{the CODE}}
As earlier mentioned, and as one can observe in the proof of~\cref{my_caratheodory}, the CODE can be solved via Picard iteration by introducing the so-called Picard operator. Therefore, let us formally introduce said operator.

Note that the CODE is well defined on the time interval $\mathbb{T}=[0, 2\pi]$. This allows us to define the Picard operator also on $\mathbb{T}$, as we see in the following. Define
\begin{equation}
\label{definition phi}
\begin{aligned}
    \Phi: C(\mathbb{T};\mathbb{R}^{n_y})\rightarrow C(\mathbb{T};\mathbb{R}^{n_y}),\    [\Phi(y;u, y^0)](t):= y^0 + \int^t_{\hat{t}_1} f(y,u) ds,
\end{aligned}
\end{equation}
where $u\in\mathcal{U}^{l,\mathbb{T},A}, l\geq 2$ and $\hat{t}_1\in\mathbb{T}$ with $J=(\hat{t}_1,\hat{t}_2)$. To keep notation simple we sometimes write $v:=(u, y^0)$ and start to use either $\Phi(y;v)$ or $\Phi(y;u,y^0)$. Further in the appendix, the~\cref{phi kappa contraction} establishes that $\Phi$ defines a $\kappa$-contraction with $\kappa<1/2$.

To conclude this section, we state the aforementioned Picard type fixed point iteration scheme that produces an approximate solution to the true solution of the CODE. This iteration scheme has the advantage to be emulated through the $\Psi$-FNO layers introduced earlier. To study the convergence of the proposed algorithm, we define the corresponding Picard operator.
For $N\in\mathbb{N}$ we define the operator $\Phi_{N,\textnormal{net}}: L^2_N(\mathbb{T};\mathbb{R}^{n_y})\rightarrow L^2_N(\mathbb{T};\mathbb{R}^{n_y})$ by
\begin{equation}
\label{definition phi n net}
\begin{aligned}
    [\Phi_{N,\textnormal{net}}(y_{N,\textnormal{net}};v_N)](t)&:={y}^0 + \int^t_{\hat{t}_1}\dot{P}_N\mathcal{I}_N f_\textnormal{net}(y_{N,\textnormal{net}}, u_N)ds + c_0(\tau_N(t)-\hat{t}_1)
\end{aligned}
\end{equation}
for $u\in \mathcal{U}^{l,\mathbb{T},A}$ and $v_N:=(u_N, y^0, \tau_N-\hat{t}_1),\ u_N:=P_N u,\ \tau_N:=P_N\tau,\ c_0 := \fint^{}_{\mathbb{T}} \mathcal{I}_N f_\textnormal{net} (y_{N,\textnormal{net}}, u_N)ds.$ We use following notation to state the $K$-times composition of the above operators
\begin{align*}
    \Phi^{[K]} &= \underbrace{\Phi\circ\Phi\dots\Phi\circ\Phi}_{K\text{-times}},
\end{align*}
analogue for $\Phi_{N,\textnormal{net}}^{[K]}$.

We introduce, as the domain to study the convergence on, the set
\begin{align*}
\mathfrak{B}_N:=\{z\in L^2_N(J;\mathbb{R}^{n_y}) \mid \| z \|_{W^{1,\infty}(\mathbb{T};\mathbb{R}^{n_y})} \le Y' \},
\end{align*}
meaning that we consider (restrictions to $J$ of) trigonometric polynomials whose $W^{1,\infty}(\mathbb{T};\mathbb{R}^{n_y})$-norm is at most $Y^\prime$. Here $L^2_N(J;\mathbb{R}^{n_y}):=\{z|_J \ | \ z\in L^2_N(\mathbb{T};\mathbb{R}^{n_y})\}$; since trigonometric polynomials are real-analytic, $z$ is uniquely determined by its restriction $z|_J$. We identify both, so that $\mathfrak{B}_N\subset C(\overline{J};\mathbb{R}^{n_y})$ and operators defined on $L^2_N(\mathbb{T};\mathbb{R}^{n_y})$ act on $\mathfrak{B}_N$.

Following algorithm is the Picard reformulation of CODE while allowing for being emulated through the layers of 
a $\Psi$-FNO. 
\begin{alg}[Approximate solution of CODE]\label{our algorithm for CODE}
Let $N\in\mathbb{N}$. 

\textbf{Input:} control $u\in\mathcal{U}^{l,J,A}$, $l\geq 2$, initial value $y^0=y(\hat{t}_1)\in Y^0$, shifted identity approximation $\tau_N-\hat{t}_1$.

\textbf{Output:} ${y}_{N,\textnormal{net}}\in L_N^2(\mathbb{T};\mathbb{R}^{n_y})$ such that ${{{y}}_{{N,\textnormal{net}}}}{|_J}\approx \hat{y}$ where $\hat{y}$ solves $\hat{y}=\Phi(\hat{y};u, y^0)|_J$.

(1) Compute $u_N:=P_N(\tau,\mathcal{E}{u})\in L^2_N(\mathbb{T};\mathbb{R}^{1+n_u})$.

(2) Set
\begin{align*}
    K&:=\lceil\frac{\log(YN)}{\log(1/\kappa)}\rceil.
\end{align*}

(3) Initial value ${y}_{N,\textnormal{net}}^{(0)}(t):=0\in \mathfrak{B}_N \subset L^2_N(\mathbb{T};\mathbb{R}^{n_y})$.

(4) For $k=1,\dots,K$ calculate
\begin{align*}
{y}^{(k)}_{N,\textnormal{net}}(t) = y^0 + \int^t_{\hat{t}_1}\dot{P}_N\mathcal{I}_N f_\textnormal{net}({y}^{(k-1)}_{N,\textnormal{net}}, u_N)ds + c_{0,k}(\tau_N(t)-\hat{t}_1) \subset L^2_N(\mathbb{T};\mathbb{R}^{n_y}),
\end{align*}
where
\begin{align*}
    c_{0,k}:= \fint^{2\pi}_{0} \mathcal{I}_N f_\textnormal{net}(y^{(k-1)}_{N,\textnormal{net}}, u_N)ds.
\end{align*}
(5) Set ${y}_{N,\textnormal{net}}:={y}^{(K)}_{N,\textnormal{net}}$.
\end{alg}
The well-definedness of above algorithm is proven in the beginning of~\cref{sec: well definedness algorithm and main result proof}. Further, note that the algorithm produces an approximate solution of the CODE. That this approximate converges to the true solution is shown in~\cref{approximation by our algorithm}.
\section{Quantitative approximation theorem}
\label{Sec:QAT}
The following quantitative approximation theorem for a $\Psi$-FNO approximating the solution operator of a CODE is our main result.
\begin{thm}\label{main theorem}
    Let~\cref{Final assumptions of ODE system} hold, including the smallness condition~\cref{choosage of E}, and let $d=1+n_y+n_u$. There are constants 
    $\varepsilon_0\in(0,\sfrac{1}{2}]$ and $C>0$ depending on $n_y,n_u,E, \alpha, m, s, l, A, Y^0, \gamma, J$, on the fixed extension operator $\mathcal{E}$ and on the GELU width parameter $M_0$, but not on $\varepsilon$, such that for every $\varepsilon\in(0,\varepsilon_0)$ there exist a $N=N(\varepsilon)\leq C\varepsilon^{-1}$ and two $\Psi$-FNOs $\widehat{\mathcal{N}}\in\{\widehat{\mathcal{N}}_\textnormal{ReLU},\widehat{\mathcal{N}}_\textnormal{GELU}\}, \widehat{\mathcal{N}}:L_N^2(\mathbb{T};\mathbb{R}^{(1+n_u)+n_y+2})\rightarrow L^2_N(\mathbb{T};\mathbb{R}^{n_y})$ with
    \begin{align*}
        \sup_{y^0\in Y^0,u\in\mathcal{U}^{l,J,A}}\Vert \hat{y}(\cdot;u,y^0) - \widehat{\mathcal{N}}\big(u_N,y^0,\tau_N-\hat{t}_1,D_N(\hat{t}_1-\cdot)\big)|_J \Vert_{C(J;\mathbb{R}^{n_y})} \leq \varepsilon,
    \end{align*}
    where $\hat{y}(\cdot;u,y^0)$ solves $\dot{y}=f(y,u)$ on $J$ with $y(\hat{t}_1)=y^0, u_N:=P_N(\tau,\mathcal{E}u)\in L^2_N(\mathbb{T};\mathbb{R}^{1+n_u})$ is the projected periodically extended control with first component $\tau_N=P_N\tau$, $D_N$ is the Dirichlet kernel and $y^0$ enters as a constant function. Moreover,
    \begin{align*}
        \textnormal{size}(\widehat{\mathcal{N}}_\textnormal{ReLU}) \leq C {(\frac{1}{\varepsilon})}^{\frac{2d}{s-1}+1}{\big(\log(\frac{1}{\varepsilon})\big)}^4,\
        \textnormal{size}(\widehat{\mathcal{N}}_\textnormal{GELU}) \leq C {(\frac{1}{\varepsilon})}^{\frac{d}{2(s-1)}+1} {\big(\log(\frac{1}{\varepsilon})\big)}^2,
    \end{align*}
    where subscript indicates whether $f$ is approximated by ReLU or a GELU network. In the case of ReLU, the remaining layers of $\widehat{\mathcal{N}}$ use the GELU activation function, see~\cref{linear approximation lemma explanation}.
\end{thm}
We see that we obtain similar results as in~\cite{kovachki2021universalapproximationerrorbounds} for their approximation regarding the Navier-Stokes equation, however, the log-terms here are taken to the power of four when using a ReLU network to approximate $f$, while they are only squared when using the GELU activation function. Neglecting the log-terms in both cases, we see that the size of the $\Psi$-FNO grows at most algebraically in ${1}/{\varepsilon}$ and that the size of $\mathcal{\widehat{N}_\textnormal{ReLU}}$ scales even sub-quadratically for the case $s> 2d+1$, while the size of $\mathcal{\widehat{N}_\textnormal{GELU}}$ scales at most sub-quadratically already when $s>d/2+1$. The latter requirement is favorable compared to $s> 2d+1$, due to the fact that $d$ is the amount of control and solution functions of a CODE system. In complex physical simulations the amount of control and solution functions can be in the range of hundreds. Nonetheless, the curse of parametric complexity is avoided in~\cref{main theorem}.
\section{Application to inverse parameter estimation}
\label{Sec:AIPE}
Motivated by the efficient approximation of the constructed FNO, we show how one can apply a FNO surrogate model in an inverse parameter estimation problem where a large data set of measurement data is available for inversion.

Consider for an interval $J\subset(0,2\pi)$ the case where one has the solution operator
\begin{align*}
    S:\mathcal{U}^{l, J, A}\times Y^0\times \prod_{i=1}^{n_x} [\hat{x}_{1,i}, \hat{x}_{2,i}]\rightarrow C(J;\mathbb{R}^{n_y})
\end{align*}
coming from a CODE system as described in the previous sections. Here, we introduce an additional parameter dependence on $x\in\prod_{i=1}^{n_x} [\hat{x}_{1,i}, \hat{x}_{2,i}]\subset\mathbb{R}^{n_x}$. One can think of $S$ being
implemented by a complex computer simulation. Then, in the presence of measurement data, it is common to infer an inverse parameter estimation such that the computer simulation provides accurate predictions
describing the real underlying physical process~\cite{Anthony}. One possible bottleneck for the parameter inversion are long inference times. Ensemble methods are derivative-free and only require evaluations of the forward operator, which the surrogate provides at low cost. Therefore, we propose the following method, which combines the EKI introduced in~\cite{Iglesias_2013}, in the form used in~\cite{iglesias2025deeponetacceleratedbayesianinversionmoving}, with the Kaczmarz method from~\cite{Leit_o_2016}.
 
Consider a physical process described by $S$ under a fixed $x^\dagger\in\prod_{i=1}^{n_x} [\hat{x}_{1,i}, \hat{x}_{2,i}]\subset\mathbb{R}^{n_x}, n_x\in\mathbb{N}$.  If said process is run under different control functions and initial values, we have a data set consisting of $N_\textnormal{exp}\in\mathbb{N}$ real experiments
\begin{align*}
    D_\textnormal{exp}:=\{(u_i,y_i^0,y_i^\delta)\}_{i=1,\dots,N_\textnormal{exp}},
\end{align*}
while an observation operator
\begin{align*}
    \mathcal{O}:C(J;\mathbb{R}^{n_y})\rightarrow\mathbb{R}^{n_\textnormal{meas}},\ n_\textnormal{meas}\in\mathbb{N},\ y\mapsto \big(y_{c_r}(t^{\mathrm{m}}_r)\big)_{r=1}^{n_{\mathrm{meas}}}
\end{align*}
maps the continuous simulation output to the discrete measurements per experiment where $y_{c_r}$ denotes the measurement of component $c_r$ at the time point $t^m_r\in J$. Varying the parameters $x$, the controls $u$ and initial values $y^0$, leads to a data set
\begin{align*}
    D_\textnormal{sim}:=\{(u_l,y_l^0,x_l, y_l)\}_{l=1,\dots,N_\textnormal{sim}},\ N_\textnormal{sim}\in\mathbb{N}
\end{align*}
generated by the simulation. This data set is used for a FNO $\widehat{\mathcal{N}}$ that serves as a surrogate for the simulation. For doing the approximation with a FNO, one treats the parameters $x$ as additional but constant controls and thus they serve as input functions to the FNO.
This allows us to introduce a forward operator per experiment, i.e. for $i=1,\dots,N_\textnormal{exp}$,
\begin{align*}
    F_i:\prod_{i=1}^{n_x} [\hat{x}_{1,i}, \hat{x}_{2,i}]\subset\mathbb{R}^{n_x}\rightarrow\mathbb{R}^{n_\textnormal{meas}},\ x\mapsto\mathcal{O}\circ\widehat{\mathcal{N}}\big((u_i, y_i^0, x)\big),
\end{align*}
where the influence of the parameters $x$ is studied under fixed control and initial value. To formulate several multiple forward operators is inspired by~\cite{Leit_o_2016}.
We also can define one single forward operator
\begin{align*}
    \mathbf{F}:\prod_{i=1}^{n_x} [\hat{x}_{1,i}, \hat{x}_{2,i}]\subset\mathbb{R}^{n_x}\rightarrow\mathbb{R}^{N_\textnormal{exp}\times n_\textnormal{meas}},\ x\mapsto\bigg(\mathcal{O}\circ\widehat{\mathcal{N}}\big((u_1, y_1^0, x)\big),\dots,\mathcal{O}\circ\widehat{\mathcal{N}}\big((u_{N_\textnormal{exp}}, y_{N_\textnormal{exp}}^0, x)\big)\bigg)
\end{align*}
for a more concise notation. From a mathematical point of view, we are interested in solving the inverse problem
\begin{align*}
    \mathbf{F}(x)=\mathbf{y}^\delta\in\mathbb{R}^{N_\textnormal{exp}\times n_\textnormal{meas}},
\end{align*}
where $\mathbf{y}^\delta$ denotes the set of all real measurements with
\begin{align*}
    \mathbf{y}^\delta:=(y_1^\delta,\dots,y_{N_\textnormal{exp}}^\delta), y_i^\delta\in\mathbb{R}^{n_\textnormal{meas}},
\end{align*}
where it is already hinted at the presence of noise in the data measurements. Depending on the amount of discrete measurements done per experiment, the number $N_\textnormal{exp}\cdot n_\textnormal{meas}$ already depicts that a parameter inversion may not be feasible when done under a solution operator with long computation time. For each measurement per experiment, we assume following relationship
\begin{align*}y_i^\delta = F_i (x^\dagger) + \varepsilon_{i,\textnormal{approx}} + \varepsilon_{i,\textnormal{sim}} + \varepsilon_{i,\textnormal{meas}},\ i=1,...,N_\textnormal{exp},
\end{align*}
where 
\begin{align*}
\varepsilon_{i,\textnormal{approx}}&\sim{N}(\overline{\varepsilon}_\textnormal{approx}, \Gamma_\textnormal{approx}),\
\varepsilon_{i,\textnormal{sim}}\sim{N}(\overline{\varepsilon}_\textnormal{sim},\Gamma_\textnormal{sim}),\
\varepsilon_{i,\textnormal{meas}}\sim{N}(0,\Gamma_\textnormal{meas})
\end{align*}
are respectively the approximation-, simulation- and measurement error. We assume that the three errors are Gaussian, mutually independent and independent of the inputs, so that their covariances add up as below. Note that we assume a simulation error, which is not uncommon in the theory of computer model calibration~\cite{Kennedy2001BayesianCO}. We can now consider the mean values $\overline{\varepsilon}_\textnormal{approx}, \overline{\varepsilon}_\textnormal{sim}$ and rewrite the relation as
\begin{align*}
y_i^\delta -\overline{\varepsilon}_\textnormal{approx}
-\overline{\varepsilon}_\textnormal{sim}= F_i (x^\dagger) + \zeta_i, i=1,...,N_\textnormal{exp},
\end{align*}
where 
\begin{align*}
\varepsilon_{i,\textnormal{approx}}-\overline{\varepsilon}_\textnormal{approx}+\varepsilon_{i,\textnormal{sim}}-\overline{\varepsilon}_\textnormal{sim}+\varepsilon_{i,\textnormal{meas}}=:\zeta_i\sim{N}(0,\Gamma_y)
\end{align*}
with an inflated noise covariance matrix defined as
\begin{align*}\Gamma_\textnormal{approx}+\Gamma_\textnormal{sim}+\Gamma_\textnormal{meas}:=\Gamma_y\in\mathbb{R}^{n_\textnormal{meas}\times n_\textnormal{meas}}.
\end{align*}
We end up with shifted measurement data set denoted by $\mathbf{y}^\zeta\in\mathbb{R}^{N_\textnormal{exp}\times n_\textnormal{meas}}$. By $y_i^\zeta\in\mathbb{R}^{n_\textnormal{meas}}, i=1,\dots,N_\textnormal{exp}$, we denote the shifted data measured per experiment for which we assume, a common noise level
\begin{align*}
    \Vert\Gamma^{-1/2} (y_i^\zeta - F_i(x^\dagger))\Vert\leq\eta\in\mathbb{R}^+.
\end{align*}
Since $\Vert\Gamma_y^{-1/2}\zeta_i\Vert^2$ is $\chi^2$-distributed with $n_\textnormal{meas}$ degrees of freedom, such a bound can only hold with a prescribed probability. For $\delta\in(0,1)$, the choice $\eta^2:=n_\textnormal{meas}+2\sqrt{n_\textnormal{meas}\log(N_\textnormal{exp}/\delta)}+2\log(N_\textnormal{exp}/\delta)$ guarantees, by the tail bound in~\cite[Lemma 1]{laurent2000adaptive} and a union bound over the experiments, that the noise level holds for all $i=1,\dots,N_\textnormal{exp}$ simultaneously with probability at least $1-\delta$. Note that $\mathbb{E}\Vert\Gamma_y^{-1/2}\zeta_i\Vert^2=n_\textnormal{meas}$, so the naive choice $\eta=\sqrt{n_\textnormal{meas}}$ is exceeded with a probability of roughly $1/2$.

The above problem is being solved using the Ensemble Kalman Inversion in which we take presence of a measurement data set into account. The idea is to iterate through each experiment of the data set. We do so by extending the Ensemble Kalman Inversion in a Kaczmarz-like style. 
Before, we introduce a component-wise transformation on parameters coming from~\cite{iglesias2025deeponetacceleratedbayesianinversionmoving, Causon2024RealtimeBI}
\begin{align*}
\mathcal{T}(x)_i = \log\Big(\frac{\hat{x}_{2,i}-x_i}{x_i-\hat{x}_{1,i}}\Big), \ \ \ \ \mathcal{T}^{-1}(\theta)_i=\frac{\hat{x}_{2,i}+e^{\theta_i}\hat{x}_{1,i}}{1+e^{\theta_i}},\quad i=1,\dots,n_x,
\end{align*}
and generate a prior ensemble $x^{(j)}_0\sim U(\prod_{i=1}^{n_x} [\hat{x}_{1,i}, \hat{x}_{2,i}]\big),\ j=1,\dots,N_\textnormal{ens}$, with $N_\textnormal{ens}\geq 2$
using the uniform distribution. Then for every iteration step we will have an updated ensemble $x^{(j)}_k$ at iteration $k$.

The Ensemble Kalman Inversion is performed in the following way.
\begin{alg}[Kaczmarz-EKI]\label{EKI Kaczmarz}
Set $\rho > 1$ and choose a maximal number of iterations $k_{\max}\in\mathbb{N}$.

\textbf{Input:} forward operators $F_i$, input and shifted measurement data $\{(u_i,y_i^0,y_i^\zeta)\},\ {i=1,\dots,N_\textnormal{exp}}$, prior ensemble $\{x^{(j)}_0\}_{j=1,\dots,N_\textnormal{ens}}$,
symmetric positive definite inflated noise covariance matrix $\Gamma_y\in\mathbb{R}^{n_\textnormal{meas}\times n_\textnormal{meas}}$, noise level $\eta>0$, and $h$ with $0<h<\min_{i}(\hat{x}_{2,i}-\hat{x}_{1,i})/2$.

For $k=0,1,\dots,k_{\max}$ and $[k]:=(k \bmod N_\textnormal{exp})+1$

(1) Apply metric projection onto $\prod_{i=1}^{n_x} [\hat{x}_{1,i}+h, \hat{x}_{2,i}-h]$, for a small $h>0$, denoted by $\mathcal{P}_h$,
\begin{align*}
x^{(j)}_k \leftarrow \mathcal{P}_h x^{(j)}_k
\ \text{for} \ j=1,\dots,N_\textnormal{ens}.
\end{align*}

(2) Evaluate forward maps $F_{i,k}^{(j)}\leftarrow F_i(x_k^{(j)})$ for $j=1,\dots,N_\textnormal{ens}$ and $i=1,\dots,N_\textnormal{exp}$.

(3) Compute
\begin{align*}
\omega_k := \begin{cases}
      1 & \Vert \Gamma_y^{-\sfrac{1}{2}}(\overline{F}_{[k],k} - y^\zeta_{[k]})\Vert > \rho\eta,\\
      0 & \text{otherwise},
      		\end{cases}
\end{align*}
where $\overline{{F}}_{i,k}:=\frac{1}{N_\textnormal{ens}}\sum_{j=1}^{N_\textnormal{ens}}{F}_{i,k}^{(j)}$ for $i=1,\dots,N_\textnormal{exp}$.

(4) Stop and return $\{x^{(j)}_k\}_{j=1,\dots,N_\textnormal{ens}}$ if $\Vert \Gamma_y^{-\sfrac{1}{2}}(\overline{F}_{[k-i+1],k-i+1} - y^\zeta_{[k-i+1]})\Vert \leq \rho\eta$ for all $i=1,\dots,N_\textnormal{exp}$.

(5) Perform component-wise transformation
\begin{align*}
\theta^{(j)}_k \leftarrow \mathcal{T}(x^{(j)}_k) \ \text{for} \ j=1,\dots,N_\textnormal{ens}.
\end{align*}

(6) Compute with
\begin{align*}
C_k^{\theta F} &= \frac{1}{N_\textnormal{ens}-1}\sum^{N_\textnormal{ens}}_{j=1}[(\theta_k^{(j)} - \overline{\theta}_k)\otimes ({F}_{[k],k}^{(j)} - \overline{{F}}_{[k],k})],\ \overline{\theta}_k:=\frac{1}{N_\textnormal{ens}}\sum_{j=1}^{N_\textnormal{ens}} \theta_k^{(j)},\\
C_k^{FF} &= \frac{1}{N_\textnormal{ens}-1}\sum^{N_\textnormal{ens}}_{j=1}[({F}_{[k],k}^{(j)} - \overline{{F}}_{[k],k})\otimes({F}_{[k],k}^{(j)} - \overline{{F}}_{[k],k})].
\end{align*}

(7) Sample $\xi_k^{(j)}\sim N(0,\Gamma_y)$.

(8) Calculate update 
\begin{align*}
\theta^{(j)}_{k+1} = \theta^{(j)}_k + \omega_k C^{\theta F}_k(\Gamma_y + C^{FF}_k)^{-1}\big(y^\zeta_{[k]}+\xi_k^{(j)}-{F}_{[k],k}^{(j)}\big),
\ j=1,\dots,N_\textnormal{ens}.
\end{align*}

(9) Transform back
\begin{align*}
x_{k+1}^{(j)}\leftarrow \mathcal{T}^{-1}(\theta^{(j)}_{k+1}).
\end{align*}

\textbf{Output:} final ensemble $\{x^{(j)}_k\}_{j=1,\dots,N_\textnormal{ens}}$ at the iteration $k$ at which the algorithm stops in step (4), or $\{x^{(j)}_{k_{\max}+1}\}_{j=1,\dots,N_\textnormal{ens}}$ if the stopping criterion is not met for $k\leq k_{\max}$.
\end{alg}
\begin{bmk}
The final ensemble of~\cref{EKI Kaczmarz} should not be interpreted as a sample of the posterior distribution: even for linear forward operators the ensemble collapses when EKI is iterated~\cite{schillings2017analysis}, and here the data of each experiment is reused in every sweep with the full covariance $\Gamma_y$. We therefore use the mean of the final ensemble as a point estimate and the discrepancy principle in step (4) as stopping rule. We do not claim regularization properties for~\cref{EKI Kaczmarz}; for EKI, such properties have been established for modified iterations, see~\cite{Iglesias_2016}. Since the forward operators are non-linear due to the FNO, we refer to~\cite{weissmann2022gradient} for a convergence analysis of EKI for non-linear forward models.
\end{bmk}

As aforementioned the notation using multiple forward operators is inspired by~\cite{Leit_o_2016} where they also use the parameter $\omega_k$ combined with a projection step. Here, the metric projection step helps avoiding numerical issues while the transformation on parameters helps avoiding to leave the parameter space used for training the surrogate model as it is stated in~\cite{iglesias2025deeponetacceleratedbayesianinversionmoving}.
Note that neither the noise level $\eta$ nor the covariance matrix $\Gamma_y$ are dependent on the experiments. This is just due to the fact that we
want to keep the presentation simple but an extension to the experiment dependent case is favorable, especially in such physical settings where the time length of an experiment can vary. Further, recall that in the beginning of this paper it was mentioned, that neural operators maintain an architecture to handle batched data, that is the simultaneous but fast evaluation on a set of data. Thus, we highlight the strength of a FNO in the above algorithm due to the requirement of performing step (2).
\section{Summary and future work}
We derived quantitative rates for approximating the solution operator of 
controlled ordinary differential equations using the implementable version of the 
Fourier neural operator. These results show that we overcome the curse of parametric 
complexity by stating that the size of the FNO grows at most algebraically in $1/\varepsilon$, and even sub-quadratically if $s>d/2+1$ (GELU) or $s>2d+1$ (ReLU), suggesting the Fourier neural operator as a suitable surrogate model for 
physical systems described by CODEs.

Examples of such systems include physical simulations of electric arc furnaces.
Those simulations require inverse parameter estimation based on data sets consisting of multiple experiments data.
This is not feasible when relying on long computational times that come with high dimensional systems.
We suggest mitigating this issue by deploying a learned FNO surrogate within an Ensemble Kalman Inversion. 
Our proposed EKI method takes multiple measurement data into account by allowing to iterate through
the data set, which we refer to as the Kaczmarz structure. Future work will demonstrate the practical effectiveness
of the described method by presenting numerical results on an electric arc furnace simulation.
\subsection*{Acknowledgements}
Dean Zenner was supported by the ``AG der Dillinger Hüttenwerke''.
\printbibliography
\input{appendixKap2}
\input{appendixKap3}
\input{appendixKap4}
\end{document}

%% file: appendixKap2.tex
\appendix
\crefalias{section}{appendix}
\section{Proof of well-posedness of CODE}\label{sec: Proof of well-posedness of CODE}

Let $d, n,n_y,n_u\in\mathbb{N}$.
Before proving~\cref{my_caratheodory} we state the well known fundamental theorem of calculus giving us
a concise recast of ODEs as definite integrals.  
\begin{thm}[Fundamental theorem of calculus]\label{Differential and Integralcalculus}
    Let $I\subset\mathbb{R}$ be open. If $y\in C_\textnormal{abs}(I;\mathbb{R}^n)$, then the derivative $\dot{y}$ exists almost everywhere in $I$, is finite and in $L^1(I;\mathbb{R}^n)$. Additionally we have
    \begin{align*}
        y(t_2)-y(t_1)=\int^{t_2}_{t_1}\dot{y}(t)dt \ \ \forall t_1, t_2\in I.
    \end{align*} 
    Conversely, if $h\in L^1(I;\mathbb{R}^n)$ and $t_1\in I$ fixed, then the function
    \begin{align*}
        y(t):=\int^t_{t_1}h(s)ds
    \end{align*}
    defines an absolutely continuous function $y\in C_\textnormal{abs}(I;\mathbb{R}^n)$ with $\dot{y}(t)=h(t)$ almost everywhere in $I$.
\end{thm}
\begin{proof}
    See E3.6 in~\cite{alt2006lineare}.
\end{proof}
Following extension of the Banach fixed point theorem allows us to derive the continuous dependence on the control functions in~\cref{my_caratheodory}.
\begin{thm}[Parameter-dependent Banach fixed point theorem]\label{Parameter-dependent Banach fixed point Theorem}
Let $X$ be a (non-empty) Banach space and $\Lambda$ a metric space. Let $\Phi:X\times \Lambda\rightarrow X$ such that for every fixed $x\in X$ the mapping $\lambda\mapsto\Phi(x, \lambda)$ is continuous.
Assume there exists a $\kappa\in(0, 1)$ with
\begin{align*}
    \Vert\Phi(x, \lambda)-\Phi(y,\lambda)\Vert_X\leq\kappa\Vert x-y\Vert_X \ \text{for all} \ \lambda\in\Lambda, x,y\in X.
\end{align*}
Then 
\begin{itemize}
    \item for every $\lambda\in\Lambda$ exists a unique fixed point $x(\lambda)\in X$ with $\Phi(x(\lambda), \lambda)=x(\lambda)$. Further for every $y\in X$
    \begin{align*}
        \Vert y-x(\lambda)\Vert_X\leq \frac{1}{1-\kappa}\Vert y-\Phi(y,\lambda)\Vert_X,
    \end{align*}
    \item the mapping $\lambda\mapsto x(\lambda)$ is continuous from $\Lambda$ to $X$.
\end{itemize}
\end{thm}
\begin{proof}
    First point is just applying the Banach fixed point theorem for each fixed $\lambda\in\Lambda$. For the second point consider a sequence $(\lambda_n)_{n\in\mathbb{N}}$ with $\lambda_n\rightarrow\lambda_0\in\Lambda$. For $y=x(\lambda_0)$
    \begin{align*}
        \Vert x(\lambda_0)-x(\lambda_n)\Vert_X\leq \frac{1}{1-\kappa}\Vert \Phi(x(\lambda_0),\lambda_0)-\Phi(x(\lambda_0),\lambda_n)\Vert_X,
    \end{align*}
    where the right hand side vanishes due to the assumed continuity.
\end{proof}
With both established results we can now give the proof of~\cref{my_caratheodory}.
\begin{proof}
    We want to apply the parameterized Banach fixed point theorem from before.
    Define for $X:=C(I;\mathbb{R}^{n_y}),\ \Lambda:=L^q(I;\mathbb{R}^{n_u})$, the operator
    \begin{align*}
        \Phi:X\times\Lambda\times\mathbb{R}^{n_y}\rightarrow X,\ [\Phi(y, u, y^0)](t):=y^0+\int^t_{t_1}f(s, y(s), u(s))ds.
    \end{align*}
    First we have to argue that $\Phi$ is well defined. For this, the integral has to exist. The first three assumptions ensure that $f$ is measurable with respect to $t$ while the fourth assumption shows that the integral is finite, since for $1\leq p<\infty$ it is
    \begin{align*}
        \Vert f(\cdot,y(\cdot),u(\cdot))\Vert_{L^p(I;\mathbb{R}^{n_y})}
        &=(\int_I\Vert f(s,y(s),u(s))\Vert_{\mathbb{R}^{n_y}}^p ds)^{1/p}\\
        &\leq (\int_I {\alpha(s)}^p ds)^{1/p}\\
        &= (\int_I \vert\alpha(s)\vert^p ds)^{1/p}\\
        &= \Vert \alpha\Vert_{L^p(I)}<\infty
    \end{align*}
    and for $p=\infty$ we have
    \begin{align*}
        \Vert f(t,y(t),u(t))\Vert_{\mathbb{R}^{n_y}} \leq \alpha(t) \leq \Vert\alpha\Vert_{L^\infty}<\infty,
    \end{align*}
    for almost all $t\in I$. Thus it is $\Vert f(t,y,u)\Vert_{L^\infty(I;\mathbb{R}^{n_y})} = \esssup_{t\in I}\Vert f(t,y(t),u(t))\Vert_{\mathbb{R}^{n_y}}<\infty$. By that we conclude $\dot{y}(\cdot)=f(\cdot,y(\cdot),u(\cdot))\in L^p(I;\mathbb{R}^{n_y})\subseteq L^1(I;\mathbb{R}^{n_y}), 1\leq p\leq\infty$ for finite domains. In the next step, we need to show that $\Phi$ defines a contraction. We show that $\Phi$ is a contraction with respect to the equivalent Bielecki norm 
    \begin{align*}
    \Vert y\Vert_X:=\sup_{t\in I}\Vert y(t)\Vert_{\mathbb{R}^{n_y}} e^{-\lambda \vert t-t_1\vert},
    \end{align*} 
    where $\lambda>0$ is chosen later.
    Fix $y^0=y(t_1)$ and $u^0\in \Lambda$. For $\hat{y}=\Phi(y, u^0, y^0)$ and $\hat{z}=\Phi(z,u^0, y^0)$ it is 
    \begin{align*}
        \Vert\hat{y}(t)-\hat{z}(t)\Vert&\leq \Big\vert\int^t_{t_1}\Vert f(s, y(s), u^0)-f(s, z(s), u^0)\Vert_{\mathbb{R}^{n_y}} ds\Big\vert \\
&\leq\Big\vert\int^t_{t_1}m(s)\Vert y(s)-z(s)\Vert_{\mathbb{R}^{n_y}} ds\Big\vert\leq\Vert y-z\Vert_X\Big\vert\int^t_{t_1}m(s)e^{\lambda \vert s-t_1\vert}ds\Big\vert.
    \end{align*}
    Hence it is $\Vert\hat{y}-\hat{z}\Vert_X\leq\kappa\Vert y-z\Vert_X$, for 
    $\kappa:=\sup_{t\in I}\big\vert\int^t_{t_1}m(s)e^{-\lambda\vert t-s\vert}ds\big\vert$, where we used $\vert t-t_1\vert-\vert s-t_1\vert=\vert t-s\vert$ for $s$ between $t_1$ and $t$.
    We write $m=m_1+m_2$ with $m_1:=m\chi_{\{ m > M \}}$ and $m_2:=m\chi_{\{ m \leq M \}}$. Choose $M>0$ large enough such that $\int_I m_1(s)ds<\sfrac{1}{4}$. Then for all $t\in I$ it is
    \begin{align*}
    \Big\vert\int^t_{t_1} m_1(s)e^{-\lambda\vert t-s\vert}ds\Big\vert\leq \int_I m_1(s)ds <\frac{1}{4}.
    \end{align*}
    Further we have by substitution
    \begin{align*}
    \Big\vert\int^t_{t_1} m_2(s)e^{-\lambda\vert t-s\vert}ds\Big\vert\leq M\Big\vert\int^t_{t_1} e^{-\lambda\vert t-s\vert}ds\Big\vert \leq M \int^\infty_0 e^{-\lambda s}ds=\frac{M}{\lambda}.
    \end{align*}
    When we choose $\lambda>4M$, we then have $\kappa<\sfrac{1}{2}$ and thus $\Phi$ is contraction measured in the Bielecki norm. We still need continuity of $u\mapsto\Phi(y,u,y^0)$ considering the $L^q$ norm and the sup-norm. 
    Note that we can now endow the space of continuous functions with the sup-norm since it is equivalent to the Bielecki norm. For $u,v\in\Lambda$ it is
    \begin{align*}
        &\sup_{t\in I}\Vert [\Phi(y,u,y^0)](t) - [\Phi(y,v,y^0)](t)\Vert_{\mathbb{R}^{n_y}}\
        =\sup_{t\in I}\Vert \int^t_{t_1}f(s, y(s), u(s))ds - \int^t_{t_1}f(s, y(s), v(s))ds\Vert_{\mathbb{R}^{n_y}} \rightarrow 0 
    \end{align*}
    for $\Vert u-v\Vert_{L^q}\rightarrow 0$, due to the last assumption in the statement.
    Hence the parameterized Banach fixed point states that for $y^0\in\mathbb{R}^{n_y}$ and $u\in\Lambda$ there exists a unique fixed point $y\in C(I;\mathbb{R}^{n_y})$ such that 
    \begin{align*}
        \Phi(y, y^0, u) = y.
    \end{align*}
    This implies
    \begin{align*}
        y(t)=[\Phi(y, u, y^0)](t)=y^0+\int^t_{t_1}f(s, y(s), u(s))ds \ \ \text{for all} \ \ t\in I,
    \end{align*}
    which is being equivalent to $y\in C_\textnormal{abs}(I;\mathbb{R}^{n_y})$ with $\dot{y}(t)=f(t, y(t), u(t))$ for almost all $t\in I$, due to the~\cref{Differential and Integralcalculus}.
    Further the mapping
    \begin{align*}
        \mathbb{R}^{n_y}\times L^q(I;\mathbb{R}^{n_u})\rightarrow C(I;\mathbb{R}^{n_y}),\
        (y^0, u) \mapsto y=\Phi(y,u,y^0)
    \end{align*}
    is continuous.
\end{proof}

Later, we study $\Phi(y;u,y^0)$ parameterized by $(u,y^0)$. Then, looking at the above proof, $\Phi$ is remains a contraction with respect to $y$ which we summarize in following corollary.
\begin{cor}\label{phi kappa contraction}
    Fix $t_1\in I$ and $y(t_1)=y^0\in\mathbb{R}^{n_y}$ and $u\in L^q(I;\mathbb{R}^{n_u}), 1\leq q\leq\infty$. Consider $\Phi(y;u,y^0): C(I;\mathbb{R}^{n_y})\rightarrow C(I;\mathbb{R}^{n_y}),\
            y\mapsto y^0 + \int_{t_1}^{t} f(s,y(s),u(s))$. Then there exists a $\kappa<\sfrac{1}{2}$ such that $\Phi$ is a $\kappa$-contraction in the Bielecki norm, meaning
    \begin{align*}
        \Vert\Phi(y;u,y^0) - \Phi(z;u,y^0)\Vert_X \leq \kappa\Vert y-z\Vert_X
    \end{align*}
    for all $y,z\in C(I;\mathbb{R}^{n_y})$.
\end{cor}
With this, we have proven the well-posedness of our CODE system.
The following proof of~\cref{bound of original ode solution} shows that the assumptions of~\cref{my_caratheodory} together with the bounded initial values $y^0\in Y^0\subset\mathbb{R}^{n_y}$ are enough to show that solutions of the CODE system are also bounded, while not relying on properties of $u\in L^q(I;\mathbb{R}^{n_u})$. This is due to the assumption (4) in the statement as we can observe in the proof of~\cref{bound of original ode solution}.
\begin{proof}
    From $y^0\in Y^0\subset\mathbb{R}^{n_y}$ bounded there exists a $S>0$ such that we have $\Vert y^0\Vert_{\mathbb{R}^{n_y}} \leq S$ for all $y^0\in Y^0\subset\mathbb{R}^{n_y}$.
    From~\cref{Differential and Integralcalculus} and Assumption (4) in~\cref{my_caratheodory} we can show that it holds for the solution that
    \begin{align*}
        \Vert y(t)\Vert_{\mathbb{R}^{n_y}} &\leq \Vert y^0\Vert_{\mathbb{R}^{n_y}} + \int^t_{t_1} \Vert f(s, y(s), u(s)) \Vert_{\mathbb{R}^{n_y}} ds\\
                        &< S + \int^{t_2}_{t_1} \alpha(s) ds\\
                        &\leq S + C\Vert \alpha \Vert_{L^p}=:Y < \infty,
    \end{align*}
    for all $t\in I$ and all $1\leq p\leq\infty$ where $C=C(p, I)>0$ is the constant obtained by estimating the $L^1$-norm upwards with the $L^p$-norm,$1\leq p\leq\infty$.
\end{proof}

Next we want to prove the boundedness of controls $u\in\mathcal{U}^{l, \mathbb{T},A}$. As they are Sobolev functions, we need a well-known Sobolev embedding statement, which will help us several times throughout this paper. 
\begin{defn}
    Consider the Sobolev space $W^{s,p}(I;\mathbb{R}^n), I\subset\mathbb{R}^d$. We call the number
    \begin{align*}
        \gamma:=s-\frac{d}{p}
    \end{align*}
    the Sobolev number.
\end{defn}
\begin{lemma}[Sobolev embedding]\label{Sobolev into continuous functions embedding}
    Let $I\subset\mathbb{R}^d$ be a bounded, open subset, $1\leq p<\infty, s\in\mathbb{N}_0$ and $\gamma$ the Sobolev number to $W^{s,p}(I;\mathbb{R}^n)$. 
    Then for $\gamma>k$ we have the continuous embedding
    \begin{align*}
        W^{s,p}(I)\hookrightarrow C^k(I;\mathbb{R}^n)\subset C(I;\mathbb{R}^n),
    \end{align*}
    where the latter are endowed with the sup norm.
    In particular by setting $p=2, k=0$ and then $s>\sfrac{d}{2}$, we have a continuous embedding
        \begin{align*}
        H^s(I;\mathbb{R}^n)\hookrightarrow C(I;\mathbb{R}^n).    
        \end{align*}
        Thus, there exists a constant $C=C(s, d)>0$ such that
\begin{align*}
    \Vert v\Vert_{L^\infty(I;\mathbb{R}^n)}\leq C\Vert v\Vert_{H^s(I;\mathbb{R}^n)}, \ \ \forall v\in H^s(I;\mathbb{R}^n).
\end{align*}
Further, when $p=2,k=1$ and then $s>\sfrac{d}{2}+1$, we have the continuous embedding
\begin{align*}
    H^s(I;\mathbb{R}^n)\hookrightarrow C^1(I;\mathbb{R}^n)\subset W^{1,\infty}(I;\mathbb{R}^n),
\end{align*}
so there exists a constant $C=C(s,d)>0$ such that
\begin{align*}
    \Vert v\Vert_{W^{1,\infty}(I;\mathbb{R}^n)}\leq C\Vert v\Vert_{H^s(I;\mathbb{R}^n)}, \ \ \forall v\in H^s(I;\mathbb{R}^n).
\end{align*}
\end{lemma}
\begin{proof}
    See Corollary 9.13 in~\cite{Brzis2010FunctionalAS}.
\end{proof}
Before we now prove~\cref{bound of periodicially extended controls}, we first investigate the boundedness of the 
control functions, when they are not periodically extended using the linear periodic extension operator.
\begin{lemma}
For  $l>\sfrac{1}{2}$ we have 
\begin{equation}\label{bound of original control function C1}
\begin{aligned}
    \Vert u\Vert_{L^\infty(I;\mathbb{R}^{n_u})} \leq C\Vert u\Vert_{H^l(I;\mathbb{R}^{n_u})}\leq CA ,\ \forall \ u\in \mathcal{U}^{l,I,A},
\end{aligned}
\end{equation}
where $C=C(H^l(I;\mathbb{R}^{n_u}), L^\infty(I;\mathbb{R}^{n_u}) )>0$ is a constant depending on $l>\sfrac{1}{2}$.
\end{lemma}
\begin{proof}
This is immediate by the Sobolev embedding~\cref{Sobolev into continuous functions embedding} and the definition of $\mathcal{U}^{l,I,A}$.
\end{proof}
This will benefit us at several places but for now it is important to note that we can conclude from it, that there exists a $U\in\mathbb{R}^+$ such that
\begin{align*}
    u(t)\in {(-U,U)}^{n_u}
\end{align*}
for all $t\in I, u\in \mathcal{U}^{l,I,A}$ since $u$ has a continuous representative.

We prove~\cref{bound of periodicially extended controls} next.
\begin{proof}
Due to the Sobolev embedding,~\cref{Sobolev into continuous functions embedding}, and the continuity of $\mathcal{E}$ we have
\begin{align*}
    \Vert\mathcal{E}u\Vert_{C(\mathbb{T})} \leq C_1\Vert \mathcal{E}u\Vert_{H^l(\mathbb{T})}
    \leq C_1C_{\mathcal{E}u}\Vert u\Vert_{H^l(J)}
    \leq C_1C_{\mathcal{E}u}A
    \leq {C}_u C_{\mathcal{E}u}A,
\end{align*}
with $C_1=C_1(H^l(\mathbb{T}),C(\mathbb{T}))>0$ and $C>0$ defined as
\begin{align}\label{Definition of Cu}
    {C}_u :=\max\{C_1,C_2\}.
\end{align}
Here, $C_2=C_2(H^l(\mathbb{T};\mathbb{R}^{n_u}),W^{1,\infty}(\mathbb{T};\mathbb{R}^{n_u}))>0$ is the constant from the embedding $H^l(\mathbb{T};\mathbb{R}^{n_u})\hookrightarrow W^{1,\infty}(\mathbb{T};\mathbb{R}^{n_u})$ when $l\geq 2$.

Consider the identity function on our time interval
\begin{align*}
    \mathrm{id}: [\hat{t}_1,\hat{t}_2]\rightarrow [\hat{t}_1,\hat{t}_2],\ t\mapsto t,
\end{align*}
which is a $\infty$-times differentiable continuous function. Define the (up to its derivative) periodically extended time variation as 
\begin{align*}
    \tau:\mathbb{T}\rightarrow\mathbb{R},\
    t\mapsto\mathcal{E}\text{id}(t).
\end{align*}
With the properties of the periodic extension operator we obtain
\begin{align*}
    \Vert\tau\Vert_{C^r(\mathbb{T})} \leq C_{\mathcal{E}\text{id}}\Vert \text{id}\Vert_{C^r([\hat{t}_1, \hat{t}_2])}\leq C_{\mathcal{E}\text{id}}\max\{1,\hat{t}_2\}
\end{align*}
for any $r\geq 0$ and where $C_{\mathcal{E}\text{id}}=C_{\mathcal{E}\text{id}}(r, \hat{t}_1,\hat{t}_2)>0$. Further, it is $\tau{|_{{J}}} = \text{id}{|_{{J}}}$. We obtain
\begin{align}\label{tau is smooth}
    \tau\in C^r(\mathbb{T};\mathbb{R})
\end{align}
for any $r\geq0$.
In the same way, we also have
\begin{align*}
    \Vert\tau\Vert_{C(\mathbb{T})} \leq C_1\Vert \tau\Vert_{H^l(\mathbb{T})}
                       \leq C_1C_{\mathcal{E}u}\Vert\text{id}\Vert_{H^l({J})}
                       \leq {C}_u C_{\mathcal{E}u}\sqrt{\frac{\hat{t}_2^3-\hat{t}_1^3}{3}+\hat{t}_2-\hat{t}_1}.
\end{align*}
    We now set
    \begin{align*}
        U^\prime = {C}_u C_{\mathcal{E}u}\Bigg(A + \sqrt{\frac{\hat{t}_2^3-\hat{t}_1^3}{3}+\hat{t}_2-\hat{t}_1}\Bigg)+ 2\pi
    \end{align*}
    to be able to handle the control $(\tau,u)$ simultaneously. Note that we then have
    \begin{align*}
        (\tau,u)\in {(-U^\prime,U^\prime)}^{1+n_u}.
    \end{align*}
    Further notice, that by the notation we introduced, it is $\Vert u\Vert_{C(\mathbb{T};\mathbb{R}^{1+n_u})} \leq U^\prime$.
\end{proof}
In the paper we mentioned that for the interval $I$ we have to at least assume $\mathbb{T}\subset I$ for the Carathéodory assumptions to hold. This is due to the fact that the image of $\tau$ is not necessarily in $\mathbb{T}$. We even can specify a concrete interval $I$. Set
\begin{align}\label{interval I}
I:=(-U^\prime-1,\,U^\prime+1),
\end{align}
which contains $[0,2\pi]$ as well as the ranges of $\tau$ and $\tau_N$, since $U^\prime>2\pi$ and both functions are bounded by $U^\prime$. It coincides with the time interval of the domain $D^\prime$.

%% file: appendixKap3.tex
\section{Fourier analysis}\label{sec: Fourier analysis}
Let $d, d_v, d_y, N\in\mathbb{N}$. Define the grid ${\{t_j\}}_{j\in\mathcal{J}_N}$ and $\mathcal{J}_N$ as in~\cref{Sec:FNOCODE}.
As we use the projection $\dot{P}_N$ in~\cref{our algorithm for CODE}, we make use of following estimation, when showing the convergence of the algorithm.
\begin{lemma}\label{estimation for dot PN}
    Let $s\in\mathbb{N}$. For a function $v\in H^s(\mathbb{T}^d)$ we have $\Vert\dot{P}_N v\Vert_{H^s(\mathbb{T}^d)}\leq\Vert v\Vert_{H^s(\mathbb{T}^d)}$.
\end{lemma}
\begin{proof}
    To see this, we work with the equivalent definition of the $H^s$-norm, see~\cite{Triebel1978InterpolationTF}, that is
    \begin{align*}
        \Vert v\Vert_{H^s(\mathbb{T}^d)}:=\big(\frac{{(2\pi)}^d}{2}\sum_{k\in\mathbb{Z}^d}(1+\vert k\vert^{2s})\vert\hat{v}_k\vert^2\big)^{1/2},
    \end{align*}
    where $\hat{v}_k=\mathcal{F}(v)(k)$ are the Fourier coefficients.
    With this, it easy to verify that
    \begin{align*}
        \Vert\dot{P}_N v\Vert_{H^s(\mathbb{T}^d)} &= \big(\frac{{(2\pi)}^d}{2}\sum_{0<\vert k\vert_\infty\leq N}(1+\vert k\vert^{2s})\vert\hat{v}_k\vert^2\big)^{1/2}\\
                                                    &\leq \big(\frac{{(2\pi)}^d}{2}\sum_{k\in\mathbb{Z}^d}(1+\vert k\vert^{2s})\vert\hat{v}_k\vert^2\big)^{1/2}
                                                    =\Vert v\Vert_{H^s(\mathbb{T}^d)}.
    \end{align*}
\end{proof}

Following spectral estimate for the pseudo-spectral projection is also needed for the proof of convergence for~\cref{our algorithm for CODE}.
\begin{lemma}[Pseudo-spectral approximation estimate]\label{Pseudo-spectral approximation estimate}
    For any $s>\sfrac{d}{2}$ the spectral interpolation operator $\mathcal{I}_N:H^s(\mathbb{T}^d)\rightarrow L^2_N(\mathbb{T}^d)$ is well defined.
    There exists a constant $C=C(s, d)>0$ such that 
    \begin{align*}
        \Vert (1-\mathcal{I}_N)v\Vert_{H^\zeta}\leq CN^{-(s-\zeta)}\Vert v\Vert_{H^s} 
    \end{align*}
    for all $v\in H^s(\mathbb{T}^d)$ and for any $\zeta\in[0,s]$.
\end{lemma}
\begin{proof}
    See Theorem 2.3 in~\cite{ShenSpectral}.
\end{proof}

As a corollary, we state the boundedness of the interpolation operator $\mathcal{I}_N$.
\begin{lemma}\label{boundedness of IN}
 For any $s>\sfrac{d}{2}$ the spectral interpolation operator $\mathcal{I}_N:H^s(\mathbb{T}^d)\rightarrow L^2_N(\mathbb{T}^d)$ is bounded 
    in the following sense
    \begin{align*}
        \Vert\mathcal{I}_N v\Vert_{H^s}\leq C\Vert v\Vert_{H^s}, \ \forall v\in H^s(\mathbb{T}^d)
    \end{align*}
    for $s>\sfrac{1}{2}$ and $C=C(s, d)>0$.
\end{lemma}
\begin{proof}
    For the one-dimensional case one can see~\cite{LiJiyong} but one can also use~\cref{Pseudo-spectral approximation estimate} with $\zeta=s$ to see that
    \begin{align*}
        \Vert\mathcal{I}_N v\Vert_{H^s} \leq&
        \Vert\mathcal{I}_N v - v\Vert_{H^s} + \Vert v\Vert_{H^s}
                                        \leq C\Vert v\Vert_{H^s} + \Vert v\Vert_{H^s}
                                        \leq (C+1)\Vert v\Vert_{H^s}.
    \end{align*}
\end{proof}
We give a spectral estimate for the $L2$-orthogonal projection.
\begin{lemma}[Spectral approximation estimate]\label{Spectral approximation estimate}
    Let $s>0$. There exists a constant $C=C(s, d)>0$ such that for any $v\in H^s(\mathbb{T}^d)$ we have 
    \begin{align*}
        \Vert v-P_N v\Vert_{H^\zeta}\leq CN^{-(s-\zeta)}\Vert v\Vert_{H^s}
    \end{align*}
    for all $v\in H^s(\mathbb{T}^d)$ and for any $\zeta\in[0, s]$.
\end{lemma}
\begin{proof}
    See Theorem 2.1 in~\cite{ShenSpectral}.
\end{proof}
In~\cref{Sec:FNOCODE} we introduced the architecture of FNOs. In the next, we present that we can decompose that architecture.
\begin{defn}[$\sigma$-layer]\label{Psi FNO nonlinear layer}
    A $\sigma$-layer $\mathcal{L}=\mathcal{L}_\sigma$ of a $\Psi$-FNO is a local, non-linear layer of the form $\mathcal{L}_\sigma(v)(x)=\mathcal{I}_N\sigma(W\mathcal{I}_N v(x)+b)$, or, in the numerical implementation 
    \begin{align*}
        {\mathcal{L}_\sigma(v)}_j=\sigma(Wv_j+b_j), \ \forall j\in\mathcal{J}_N,
    \end{align*} 
    with $W\in\mathbb{R}^{d_v\times d_v}$ and $b_j\in\mathbb{R}^{\mathcal{J}_N\times d_v}$ defining an affine mapping.
\end{defn}

\begin{defn}[$\mathcal{F}$-layer]\label{Psi Fno linear layer}
    A $\mathcal{F}$-layer $\mathcal{L}=\mathcal{L}_\mathcal{F}$ of a $\Psi$-FNO is a non-local, linear layer of the form $\mathcal{L}_\mathcal{F}(v)(x)=\mathcal{F}^{-1}(P(k)\cdot\mathcal{F}(\mathcal{I}_N v)(k))(x)$,
    which corresponds in implementation to
    \begin{align*}
        \mathcal{L}_\mathcal{F}{(v)}_j = \mathcal{F}^{-1}_N{\big(P(k)\cdot\mathcal{F}_N(v)(k)\big)}_j, \ \forall j\in\mathcal{J}_N,
    \end{align*}
    where $P:\mathcal{K}_N\rightarrow\mathbb{C}^{d_v\times d_v}$ is a collection of complex weights, 
    with $P(-k)=\overline{P(k)}$ (entry-wise complex conjugation),
    and $\mathcal{F}_N (\mathcal{F}_N^{-1})$ denotes the discrete (inverse) Fourier transform.
\end{defn}
Each $\Psi$-FNO can be decomposed into a finite number of $\sigma$-layers and $\mathcal{F}$-layers. The converse is also true. This is proven in the Appendix D.7 of~\cite{kovachki2021universalapproximationerrorbounds} by proving following lemmata which are also useful for proving our result.

\begin{lemma} [Composition Lemma]\label{Composition lemma}
    For $p\in\mathbb{N}$, let $\mathcal{N}_1,\dots,\mathcal{N}_p$ be $\Psi$-FNOs, defined with respect to the same grid ${\{t_j\}}_{j\in\mathcal{J}_N}$,
    and defining operators $\mathcal{N}_k: L^2_N(\mathbb{T}^d; \mathbb{R}^{d_{k-1}})\rightarrow L^2_N(\mathbb{T}^d; \mathbb{R}^{d_k})$ for $d_0,\dots,d_p\in\mathbb{N}$,
    such that the composition $\mathcal{N}_p\circ\dots\circ\mathcal{N}_1: L^2_N(\mathbb{T}^d; \mathbb{R}^{d_0})\rightarrow L^2_N(\mathbb{T}^d; \mathbb{R}^{d_p})$ is well defined.
    Then there exists a $\Psi$-FNO $\mathcal{N}:L^2_N (\mathbb{T}^d; \mathbb{R}^{d_0})\rightarrow L^2_N(\mathbb{T}^d; \mathbb{R}^{d_p})$, such that 
    \begin{align*}
    \mathcal{N}= \mathcal{N}_p\circ\dots\circ\mathcal{N}_1,
    \end{align*}
    and such that $\textnormal{depth}(\mathcal{N})\leq\sum^p_{k=1}\textnormal{depth}(\mathcal{N}_k), \textnormal{width}(\mathcal{N})\leq \max_{k=1,\dots,p} \textnormal{width}(\mathcal{N}_k)$,
    $\textnormal{lift}(\mathcal{N})\leq \max_{k=1,\dots,p}\textnormal{lift}(\mathcal{N}_k)$.
\end{lemma}
\begin{proof}
    See Lemma D.4 in~\cite{kovachki2021universalapproximationerrorbounds}.
\end{proof}

\begin{lemma}[Linear approximation lemma]\label{Linear approximation lemma}
    Let $\sigma\in C^2$ be a twice continuously differentiable and non-constant activation function. Consider a linear $\Psi$-FNO layer ${\mathcal{L}}$ of form
    \begin{align*}
        {\mathcal{L}}:L^2_N(\mathbb{T}^d;\mathbb{R}^{dv})\rightarrow L^2_N(\mathbb{T}^d;\mathbb{R}^{dv}),\
        v_N \mapsto Wv_N(t_j) + b_j + \mathcal{F}_N^{-1}{(P\mathcal{F}_N v_N)}_j.
    \end{align*}
    Then there exists a constant $C>0$, such that for any $\varepsilon, B>0$, there exists a $\Psi$-FNO $\mathcal{N}:L^2_N(\mathbb{T}^d;\mathbb{R}^{d_v})\rightarrow L^2_N(\mathbb{T};\mathbb{R}^{d_v})$, such that
    \begin{align*}
        \sup_{\Vert v_N\Vert_{L^2}\leq B}\Vert\mathcal{N}(v_N) - {\mathcal{L}}(v_N)\Vert_{L^2}\leq\varepsilon,
    \end{align*}
    and $\textnormal{width}(\mathcal{N})\leq CN^d, \textnormal{depth}(\mathcal{N})\leq C, \textnormal{lift}(\mathcal{N})\leq C$.
\end{lemma}
\begin{proof}
    See Lemma D.3 in~\cite{kovachki2021universalapproximationerrorbounds}.
\end{proof}

\begin{lemma}[Replacement lemma]\label{Replacement lemma}
    Assume that $\sigma$ is Lipschitz continuous. For $p\in\mathbb{N}$, let $\widehat{\mathcal{N}}_1,\dots, \widehat{\mathcal{N}}_p$ be continuous operators $\widehat{\mathcal{N}}_k : L^2_N(\mathbb{T}^d; \mathbb{R}^{d_{k-1}})\rightarrow L^2_N(\mathbb{T}^d;\mathbb{R}^{d_k})$ for $d_0,\dots, d_p\in\mathbb{N}$, such that the composition $\widehat{\mathcal{N}}:=\widehat{\mathcal{N}}_p\circ\dots\circ\widehat{\mathcal{N}}_1: L^2_N(\mathbb{T}^d;\mathbb{R}^{d_0})\rightarrow L^2_N(\mathbb{T}^d;\mathbb{R}^{d_p})$ is well-defined. Assume that there exist constants $D_k, W_k, L_k > 0, k=1,\dots,p,$ such that for any $\varepsilon, M > 0$, there exists a $\Psi$-FNO $\mathcal{N}_k:L^2_N(\mathbb{T}^d;\mathbb{R}^{d_{k-1}})\rightarrow L^2_N(\mathbb{T}^d;\mathbb{R}^{d_k})$, such that  
    \begin{align*}
    \sup_{\Vert u\Vert_{L2}\leq M}\Vert \widehat{\mathcal{N}}_k(u) - \mathcal{N}_k(u)\Vert_{L_N^2}\leq\varepsilon,
    \end{align*}
    and $\textnormal{depth}(\mathcal{N}_k)\leq D_k, \textnormal{width}(\mathcal{N}_k)\leq W_k, \textnormal{lift}(\mathcal{N}_k)\leq L_k$. Then for any $\varepsilon> 0$ and  $M > 0$, there exists a $\Psi$-FNO $\mathcal{N}: L^2_N (\mathbb{T}^d;\mathbb{R}^{d_0})\rightarrow L^2_N(\mathbb{T}^d;\mathbb{R}^{d_p})$, such that  
    \begin{align*}
    \sup_{\Vert u\Vert_{L^2}\leq M}\Vert \widehat{\mathcal{N}}(u)-\mathcal{N}(u)\Vert_{L^2_N}\leq\varepsilon,
    \end{align*}
    and $\textnormal{depth}(\mathcal{N})\leq\sum^p_{k=1} D_k, \textnormal{width}(\mathcal{N})\leq \max_{k=1,...,p} W_k, \textnormal{lift}(\mathcal{N})\leq \max_{k=1,...,p} L_k$.
\end{lemma}
\begin{proof}
    See Lemma D.5 in~\cite{kovachki2021universalapproximationerrorbounds}.
\end{proof}

\begin{lemma}\label{Psi FNO representation of I_N NN}
    Let $\phi:\mathbb{R}^{d_v}\rightarrow \mathbb{R}^{d_y}$ be a (ordinary) neural network with activation function $\sigma$. For any $N\in\mathbb{N}$, the mapping 
    \begin{align*}
    \mathcal{N}: L^2_N (\mathbb{T}^d; \mathbb{R}^{d_v})&\rightarrow L^2_N(\mathbb{T}^d; \mathbb{R}^{d_y}),\\
     v(x)&\mapsto \mathcal{N}(v)(x):=\mathcal{I}_N\phi(v(x))
    \end{align*}
    can be represented by a $\Psi$-FNO, with 
    \begin{align*}
    \textnormal{width}(\mathcal{N}) = \vert \mathcal{J}_N\vert \textnormal{width}(\phi), \textnormal{depth}(\mathcal{N}) = \textnormal{depth}(\phi), \textnormal{lift}(\mathcal{N}) = \textnormal{width}(\phi).
    \end{align*}
\end{lemma}
\begin{proof}
    See Lemma D.2 in~\cite{kovachki2021universalapproximationerrorbounds}.
\end{proof}

\begin{lemma}\label{approximation of multiplication with nn}
    Consider an activation function $\sigma$ with $\sigma\in C^3$. Then there exists a constant $c_\phi=c_\phi(d)>0$, such that for any $\varepsilon>0$ and $B>0$, there exists
    \begin{itemize}
        \item[(1)] a neural network $\phi_1:\mathbb{R}^2\rightarrow\mathbb{R}$, such that
                    \begin{align*}
                        \sup_{\vert a\vert,\vert b\vert\leq B}\vert{\phi_1}(a,b)-ab\vert\leq\varepsilon,
                    \end{align*}
                    and \textnormal{width}$(\phi_1)\leq c_\phi$ and \textnormal{depth}$(\phi_1)\leq c_\phi$,
        \item[(2)] a neural network $\phi_2:\mathbb{R}\times\mathbb{R}^d\rightarrow\mathbb{R}^d$, such that
                    \begin{align*}
                        \sup_{\vert a\vert,\Vert v\Vert_{\mathbb{R}^d}\leq B}\vert{\phi_2}(a,v)-av\vert\leq\varepsilon,
                    \end{align*}
                    and \textnormal{width}$(\phi_2)\leq c_\phi$ and \textnormal{depth}$(\phi_2)\leq c_\phi$.
    \end{itemize}
\end{lemma}
\begin{proof}
    See Lemma E.2. \@in~\cite{kovachki2021universalapproximationerrorbounds}.
\end{proof}

Following results are needed for the emulation of~\cref{our algorithm for CODE} through $\Psi$-FNO layers. First, we clarify that the coefficients of the discrete Fourier transform are in this setting exactly the coefficients of the continuous transform, i.e. the coefficients of a trigonometric polynomial.
\begin{lemma}\label{Discrete Fourier coefficients of trigonometric polynomials}
    For $k\in\mathbb{Z}$ with $\vert k\vert\leq N$ and $v_N\in L^2_N(\mathbb{T})$ with $v_N(t):=\sum_{\vert k\vert_\infty\leq N} c_k e^{i\langle t, k\rangle}$
    we have 
    \begin{align*}    
    \mathcal{F}_N(v_N)(k)=c_k=\mathcal{F}(v_N)(k).
    \end{align*}
\end{lemma}
\begin{proof}
We have
    \begin{align*}
        \mathcal{F}_N(v_N)(k)&=\frac{1}{2N+1}  \sum_{j\in\mathcal{J}_N}{v_N}_j e^{-2\pi i j k / (2N+1)}\\
                            &=\frac{1}{2N+1}  \sum_{j\in\mathcal{J}_N}v_N(t_j)e^{-2\pi i j k / (2N+1)}\\
                            &=\frac{1}{2N+1}  \sum_{j\in\mathcal{J}_N}[\sum_{\vert l\vert\leq N}c_l e^{i t_j l}]e^{-2\pi i j k / (2N+1)}\\
                            &=\frac{1}{2N+1}  \sum_{j\in\mathcal{J}_N}[\sum_{\vert l\vert\leq N}c_l e^{2\pi i j l / (2N+1)}]e^{-2\pi i j k / (2N+1)}\\
                            &=\frac{1}{2N+1}  \sum_{j\in\mathcal{J}_N}\sum_{\vert l\vert\leq N}c_l 
                            e^{2\pi i j(l-k)/(2N+1)}\\
                            &=\frac{1}{2N+1}  \sum_{\vert l\vert\leq N}c_l\sum_{j\in\mathcal{J}_N}e^{2\pi i j(l-k)/(2N+1)}\\
                            &{=}\frac{1}{2N+1}  \sum_{\vert l\vert\leq N}c_l(2N+1)\delta_{lk}
                            = c_k,
    \end{align*}
where $\delta_{lk}$ defines the Kronecker delta and note, that we used the summation property of the root of unity.
\end{proof}

Next result is based on the aforementioned differentiation formula of the Fourier transform. Here, it shows how one can calculate the integral of a trigonometric polynomial with zero mean using the discrete (inverse) Fourier transform.
\begin{lemma}\label{Integral Discrete Fourier Transform Formula}
 Let $v_N\in L^2_N(\mathbb{T})$ with Fourier coefficients $c_k=\mathcal{F}_N(v_N)(k)$, $\vert k\vert\leq N$, let $\delta_0(k):=\frac{1}{ik}$ for $k\neq 0$ and $\delta_0(0):=0$, and define the trigonometric polynomial
 \[
  V_N(t):=\sum_{\vert k\vert\leq N}\delta_0(k)\,c_k\,e^{ikt},\quad t\in\mathbb{T}.
 \]
 Then $V_N\in L^2_N(\mathbb{T})$, $V_N'=\dot{P}_N v_N$, and for all $a,t\in\mathbb{T}$ 
 \[
  \int_a^t(\dot{P}_N v_N)(s)\,ds = V_N(t)-V_N(a).
 \]
 On the grid, $V_N(t_j)=\mathcal{F}_N^{-1}(\delta_0(k)\mathcal{F}_N(v_N)(k))(j)$ for $j\in\mathcal{J}_N$, i.e.\ $V_N$ is obtained from $v_N$ by a Fourier multiplier.
\end{lemma}
\begin{proof}
Term-wise differentiation gives $V_N'(t)=\sum_{0<\vert k\vert\leq N}c_k e^{ikt}=(\dot{P}_N v_N)(t)$. The claim follows from the fundamental theorem of calculus while the statement on the grid is immediately obtained due to~\cref{Discrete Fourier coefficients of trigonometric polynomials} and by the definition of the discrete (inverse) Fourier transform. 
\end{proof}
To calculate the mean over a function inside~\cref{our algorithm for CODE}, we make use of the following result.
\begin{lemma}\label{mean repr using psi fno}
    For the mean value of $v_N\in L^2_N(\mathbb{T})$ it is 
    \begin{align*}
        \fint_{\mathbb{T}} v_N(t)dt = \mathcal{F}_N^{-1}(\delta_1(k)\mathcal{F}_N(v_N)(k))(j)
    \end{align*}
    for all $j\in\mathcal{J}_N$ with 
    \begin{align*}
        \delta_1(k) = \begin{cases}
            0, \qquad k\not=0,\\
            1, \qquad k=0.
        \end{cases}
    \end{align*}
\end{lemma}
\begin{proof}
    It is 
    \begin{align*}
        \fint_\mathbb{T} v_N(t)dt = c_0
        =\sum_{\vert k\vert\leq N}\delta_1(k){c_k}e^{i {t_j} k}
        =\mathcal{F}_N^{-1}(\delta_1(k)\mathcal{F}_N(v_N)(k))(j)
    \end{align*}
    for all $j\in\mathcal{J}_N$ with corresponding $t_j\in\mathbb{T}$.
    \end{proof}

In~\cref{our algorithm for CODE}, it is also required to do a point-wise evaluation with the output function from the FNO layer before. Using the Dirichlet-kernel, which serves then as an additional input to the FNO, allows us to obtain such evaluations inside the FNO's architecture. 
\begin{lemma}\label{Dirichlet kernel property}
    Let $v_N\in L^2_N(\mathbb{T})$ and $a\in\mathbb{T}$, then it is
        \begin{align*}
            v_N(a) = \fint_\mathbb{T} v_N(t) D_N(a-t)dt,
        \end{align*}
    where $D_N$ defines the Dirichlet kernel.
\end{lemma}
\begin{proof}
    It is 
    \begin{align*}
        \fint_\mathbb{T} v_N(t) D_N(a-t)dt&= \sum_{\vert k\vert\leq N} \fint_\mathbb{T} v_N(t)e^{ik(a-t)} dt= \sum_{\vert k\vert\leq N}\mathcal{F}(v_N)(k)e^{ika}= v_N(a).
    \end{align*}
\end{proof}
As a $\Psi$-FNO $\mathcal{N}:L^2_N(\mathbb{T};\mathbb{R}^{n_a})\rightarrow L^2_N(\mathbb{T};\mathbb{R}^{n_a})$ outputs functions of form $\mathcal{I}_Ny$ in each layer and the multiplication of two $L^2_N$ functions lies in $L^2_{2N}$, we need following lemma.
\begin{lemma}\label{discrete and continuous mean align}
    Let $v_{2N}\in L^2_{2N}(\mathbb{T})$, then we have
    \begin{align*}
        \fint_{\mathbb{T}}\mathcal{I}_{N}v_{2N}(s)ds= \fint_\mathbb{T} v_{2N}(s)ds.
    \end{align*}
\end{lemma}
\begin{proof}Consider $v_{2N}\in L^2_{2N}(\mathbb{T})$ with Fourier coefficients $c_k$. Since $\frac{1}{2N+1}\sum_{j=0}^{2N}e^{ikt_j}=0$ for $0<\vert k\vert\leq 2N$, it is
\begin{align*}
    \frac{1}{2N+1}\sum_{j=0}^{2N}v_{2N}(t_j) = \frac{1}{2N+1}\sum_{j=0}^{2N}\sum_{\vert k\vert\leq 2N}c_ke^{ikt_j}= \sum_{\vert k\vert\leq 2N}c_k\frac{1}{2N+1}\sum_{j=0}^{2N}e^{ikt_j}=c_0=\fint_\mathbb{T} v_{2N}(s)ds.
\end{align*}
The same computation applied to the trigonometric polynomial $\mathcal{I}_N v_{2N}$ of degree at most $N$ gives $\fint_\mathbb{T}\mathcal{I}_N v_{2N}(s)ds=\frac{1}{2N+1}\sum_{j=0}^{2N}\mathcal{I}_N v_{2N}(t_j)$. We conclude by noting that $\mathcal{I}_Nv_{2N}(t_j)=v_{2N}(t_j)$ for $j\in\mathcal{J}_N$.
\end{proof}

%% file: appendixKap4.tex
\section{Proof of \texorpdfstring{\cref{main theorem}}{the main theorem}}
\label{sec: well definedness algorithm and main result proof}
Let $n,n_y,N\in\mathbb{N}$.
First we clarify the well-definedness of our proposed~\cref{our algorithm for CODE}. In the algorithm constant values appear which we consider also as trigonometric polynomials due to the next lemma. That this is well-defined shows following statement.
\begin{lemma}\label{constants are L2N}
    For $y\in \mathbb{R}^{n}$ we have that the function give by $y(t):=y$ is in $L^2_N(\mathbb{T};\mathbb{R}^{n})$.
\end{lemma}
\begin{proof}
    It is 
    \begin{align*}
        y(t) = \sum_{\vert k\vert\leq N} c_k e^{ikt} = y
    \end{align*}
    with $c_k=0, k\not=0$ and $c_0=y$.
\end{proof}
To make sense of the expression under the integral, we define 
\begin{align*}
    g_\textnormal{net}(t):=f_\textnormal{net}(y_{N,\textnormal{net}}(t), u_N(t)),
\end{align*}
where $y_{N,\textnormal{net}}\in L_N^2(\mathbb{T};\mathbb{R}^{n_y})$ is any trigonometric polynomial of degree $N$. Then, the expression under the integral is understood as a function
\begin{align*}
    t\mapsto [\dot{P}_N\circ\mathcal{I}_N](g_\textnormal{net}(t)),\ t\in\mathbb{T},
\end{align*}
meaning that the Fourier projection and interpolation are applied with respect to the time variable.
To apply $\mathcal{I}_N$ on $g_\textnormal{net}$, the function $g_\textnormal{net}$ needs to be at least continuous on $\mathbb{T}$.
\begin{bmk}
    We have $L^2_N(\mathbb{T};\mathbb{R}^n)\subset C^\infty(\mathbb{T};\mathbb{R}^n)$.
\end{bmk}
\begin{lemma}
    The function $g_\textnormal{net}$ is a continuous function.
\end{lemma}
\begin{proof}
    Since we know that $f_\textnormal{net}$ is continuous as function 
    \begin{align*}
        \mathbb{R}^{n_y}\times\mathbb{R}^{1+n_u} \rightarrow\mathbb{R}^{n_y},\ (y,u) \mapsto f_\textnormal{net}(y,u),
    \end{align*}
    we know that $g_\textnormal{net}$ is a continuous function as composition of such, as we have 
    \begin{align*}
    u_N\in L^2_N(\mathbb{T};\mathbb{R}^{1+n_u}), l\geq 2,\ \text{and} \
    y_{N,\textnormal{net}}\in L^2_N(\mathbb{T};\mathbb{R}^{n_y}).
    \end{align*}      
\end{proof}
Therefore we can apply $\mathcal{I}_N$ on $g_\textnormal{net}$ and thus apply $\dot{P}_N$ which then solely subtracts the mean of $\mathcal{I}_N f_\textnormal{net}(y,u)$.
Following lemma resolves the integrability of the studied function.
\begin{lemma}
    We have that $t\mapsto\dot{P}_N\mathcal{I}_N g_\textnormal{net}(t),\ t\in\mathbb{T}$, is integrable.
\end{lemma}
\begin{proof}
    After applying $\dot{P}_N$ we end up with a trigonometric polynomial of zero mean. As trigonometric polynomials are continuous, 
    we can integrate the expression $\dot{P}_N\mathcal{I}_N g_\textnormal{net}(t)$ over compact sets.
    Alternatively, we have
    \begin{align*}
        \dot{L}^2_N(\mathbb{T};\mathbb{R}^{n_y})\subset L^2_N(\mathbb{T};\mathbb{R}^{n_y})\subset L^2(\mathbb{T};\mathbb{R}^{n_y})\subset L^1(\mathbb{T};\mathbb{R}^{n_y}).
    \end{align*}    
\end{proof}
It remains the question whether the fixed point formulation of our proposed algorithm is a trigonometric polynomial. 
\begin{lemma}\label{integration is closed in L2N}
    When we have $v_N\in L^2_N(\mathbb{T})$ with $\fint_\mathbb{T} v_N(s)ds=0$, then also 
    \begin{align*}
    w(t):=\int^t_{a} v_N(s)ds\in L^2_N(\mathbb{T}),
    \end{align*}
    where $a\in\mathbb{T}$.
\end{lemma}
\begin{proof}
    We have
    \begin{align*}
        w(t)&=\int^t_{a} v_N(s)ds\\
            &=\int^t_{a} \sum_{\vert k\vert\leq N} c_k e^{iks}ds\\
            &=\sum_{\vert k\vert\leq N}\int^t_{a} c_k e^{iks}ds\\
            &=\int^t_{a} c_0e^{0}ds + \sum_{0<\vert k\vert\leq N}\int^t_{a} c_k e^{iks}ds\\\\
            &=(t-a)c_0+\sum_{0<\vert k\vert\leq N}\frac{c_k}{ik}(e^{ikt}-e^{ika})\\
            &=(t-a)c_0+\sum_{0<\vert k\vert\leq N}\frac{c_k}{ik}e^{ikt}-\sum_{0<\vert k\vert\leq N}\frac{c_k}{ik}e^{ika}\\
            &=\sum_{0<\vert k\vert\leq N}\frac{c_k}{ik}e^{ikt}-\sum_{0<\vert k\vert\leq N}\frac{c_k}{ik}e^{ika},
    \end{align*}
    where we notice that the last summand does not depend on $t$ and is therefore just a constant in $L^2_N$. Further, we know that $c_0=0$ due to the assumption.
    So by setting $d_0:=-\sum_{0<\vert k\vert\leq N}\frac{c_k}{ik}e^{ika}$ and $d_k=\frac{c_k}{ik}, k\not=0$, we have that $w$ is of form
    \begin{align*}
        w(t) = \sum_{{\vert k\vert\leq N}}d_k e^{ikt}\in L^2_N(\mathbb{T}).
    \end{align*}
\end{proof}
We conclude the well-definedness of~\cref{our algorithm for CODE} with following lemma.
\begin{lemma}
    For the iterates $y_{N,\textnormal{net}}^{(k)}, k=0,\dots,K$ it is 
    \begin{align*}
    y_{N,\textnormal{net}}^{(k)}\in L_N^2(\mathbb{T};\mathbb{R}^{n_y})\subset C^\infty(\mathbb{T};\mathbb{R}^{n_y}).
    \end{align*}
\end{lemma}
\begin{proof}
We have that 
    \begin{align*}
    t\mapsto y^0 + \int^t_{\hat{t}_1} \dot{P}_N\mathcal{I}_N g_\textnormal{net}(s)ds + c_0(\tau_N(t)-\hat{t}_1),\ t\in\mathbb{T},
    \end{align*}
    is a function in $L^2_N(\mathbb{T};\mathbb{R}^{n_y})$.
    As the function we integrate is a trigonometric polynomial with zero mean, we again end up with a trigonometric polynomial in $L^2_N(\mathbb{T};\mathbb{R}^{n_y})$,
    see~\cref{integration is closed in L2N}. Since constants are in $L^2_N(\mathbb{T})$, see~\cref{constants are L2N}, and since $\tau_N\in L^2_N(\mathbb{T})$ is also a trigonometric polynomial we are done, due to that $L^2_N(\mathbb{T})$ defines a vector space.
    For the constant function $y_{N,\textnormal{net}}^{(0)}(t)=0$ it is obvious. Inductively, it follows
    \begin{align*}
        y_{N,\textnormal{net}}^{(k)}(t)=y^0+\int^t_{\hat{t}_1} \dot{P}_N\mathcal{I}_N f_\textnormal{net}(y_{N,\textnormal{net}}^{(k-1)}(s), u_N(s))ds+c_{0,k}(\tau_N(t)-\hat{t}_1)\in L^2_N(\mathbb{T};\mathbb{R}^{n_y})
    \end{align*}
    for $k=1,\dots,K$.
\end{proof}
With this, we have proved that the presented algorithm is well defined. We dedicate the rest of this section, to that the proposed algorithm converges as well as that we can emulate that algorithm using layers from the $\Psi$-FNO. We start by choosing certain quantities in a suitable way.

Under~\cref{Final assumptions of ODE system} we consider the following selection of values.
Choose 
\begin{align}\label{choosage of E}
    \Vert f\Vert_{W^{s,\infty}(D^\prime;\mathbb{R}^{n_y})}&\leq E:=\frac{\gamma}{4\widehat{K}\max\{U^\prime,Y^\prime\}}\in\mathbb{R}^+,
\end{align}
where we recall that $\gamma>0$ was chosen such that $Y^\prime =Y+\gamma>1$. We refer to~\cref{choosage of E}
as smallness condition. The constant $\widehat{K}>0$ is defined by
\begin{align*}
    \widehat{K} := (1+U^\prime)2\pi K_1K_2K_3K_4,
\end{align*}
where $K_1>0$ is the constant from the embedding
\begin{align*}
    H^1(\mathbb{T};\mathbb{R}^{n_y})\hookrightarrow L^\infty(\mathbb{T};\mathbb{R}^{n_y}),
\end{align*}
$K_2>0$ is the constant from~\cref{boundedness of IN} with $s=1$, so 
\begin{align*}
        \Vert\mathcal{I}_N v\Vert_{H^1(\mathbb{T})}\leq K_2\Vert v\Vert_{H^1(\mathbb{T})}, \ \forall v\in H^1(\mathbb{T}),
    \end{align*}
$K_3>0$ is the constant from the embedding
\begin{align*}
    W^{1,\infty}(\mathbb{T};\mathbb{R}^{n_y})\hookrightarrow H^1(\mathbb{T};\mathbb{R}^{n_y}),
\end{align*}
see~\cref{Sobolev into continuous functions embedding} and $K_4>0$ is a constant from~\cref{Sobolev composition estimate our case}, which gives us an estimate for the composition of Sobolev functions.
The accuracy levels for the network are determined in the following.
For $N\in\mathbb{N}$ choose $\varepsilon_\textnormal{ReLU}$ with
\begin{align}\label{choosage of epsilon ReLU}
    \varepsilon_\textnormal{ReLU} =N^{-1},
\end{align}
and the number $L\in\mathbb{N}$ of the GELU network as
\begin{align}\label{choosage of epsilon GELU}
    L=\lceil(c_{0,\textnormal{GELU}}N)^{d/(2(s-1))}\rceil\in\mathbb{N},
\end{align}
implying for $\varepsilon_\textnormal{GELU}=c_{0,\textnormal{GELU}}L^{-2(s-1)/d}\leq N^{-1}.$
Recall that the amount of iterations in~\cref{our algorithm for CODE} was chosen as $K=\lceil\log(YN)/\log(1/\kappa)\rceil$.
Together with~\cref{application of disturbance lemma} and~\cref{difference phi and phi N net} we obtain the existence of a constant $C^*=C^*(n_y,n_u,E, \alpha, m, s, l, \mathcal{E}, A, Y^0,\gamma, M_0, J)>0$. Since it can be the case that $3C^*\min\{\sfrac{1}{2},E\}>1$, we set for the main result
\begin{align*}
    \varepsilon\in(0,\min\{\sfrac{1}{2},3C^*\min\{\sfrac{1}{2},E\}\}),
\end{align*}
define $\varepsilon_0:=\min\{\sfrac{1}{2},3C^*\min\{\sfrac{1}{2},E\}\}$ and choose $N\in\mathbb{N}$ by
\begin{align*}
    N=\lceil\frac{3C^*}{\varepsilon}\rceil.
\end{align*}
Otherwise, we would allow for $\log(\varepsilon^{-1})<0$ to happen. Then, for both accuracy levels of the finite neural networks, we have
\begin{align*}
    \varepsilon_\textnormal{net}<\min\{\sfrac{1}{2}, E\},
\end{align*}
which is needed for the existence of the network, as well it is needed for the proof of~\cref{image of iterates}. It is important to note that no circular reasoning is present since $C^*$ is not dependent on the constructed networks $f_\textnormal{net}$. The constant is dependent on the width number $M_0$ that we fix for the GELU network architecture beforehand. The proof of~\cref{image of iterates} could also be done by just requiring $\varepsilon_\textnormal{net}<\min\{\sfrac{1}{2},E\}.$

Following lemma helps us proving the convergence of~\cref{our algorithm for CODE} to the true solution of the CODE by solely utilizing the fact that the Picard operator $\Phi$ already defines a $\kappa$-contraction.
\begin{lemma}[Disturbance lemma]\label{Disturbance lemma}
    Let $\mathcal{B}\subset \mathcal{X}$ be a closed set in a Banach space $\mathcal{X}$. Further, let $\Phi:\mathcal{X}\rightarrow \mathcal{X}$ be a $\kappa$-contraction. Then assume we have an operator $\widehat{\Phi}$ with properties
    \begin{itemize}
        \item $\widehat{\Phi}(\mathcal{B})\subset \mathcal{B}$,
        \item $\sup_{y\in \mathcal{B}}\Vert\Phi(y) - \widehat{\Phi}(y)\Vert_\mathcal{X} \leq\eta$.
    \end{itemize}
    Then, we have 
    \begin{align*}
        \Vert\Phi^{[K]}(y^0) - \widehat{\Phi}^{[K]}(y^0)\Vert_\mathcal{X}\leq \frac{\eta}{1-\kappa}
    \end{align*}
    for any $y^0\in \mathcal{B}$.
\end{lemma}
\begin{proof}
    Define
    \begin{align*}
        e_k := \Vert \Phi^{[k]}(y^0) - \widehat{\Phi}^{[k]}(y^0)\Vert_\mathcal{X}.
    \end{align*}
    Then recursively we obtain
    \begin{align*}
        e_K &= \Vert \Phi^{[K]}(y^0) - \widehat{\Phi}^{[K]}(y^0)\Vert_\mathcal{X}\\
            &= \Vert \Phi(\Phi^{[K-1]}(y^0)) - \widehat{\Phi}(\widehat{\Phi}^{[K-1]}(y^0))\Vert_\mathcal{X}\\
            &\leq \Vert \Phi(\Phi^{[K-1]}(y^0)) - \Phi(\widehat{\Phi}^{[K-1]}(y^0))\Vert_\mathcal{X} + \Vert\Phi(\widehat{\Phi}^{[K-1]}(y^0)) - \widehat{\Phi}(\widehat{\Phi}^{[K-1]}(y^0))\Vert_\mathcal{X}\\
            &\leq \kappa e_{K-1} + \eta\\
            &\leq \kappa(\kappa e_{K-2} + \eta) + \eta\\
            &\leq\dots\\
            &\leq \sum_{k=0}^{K-1} \kappa^k\eta, 
    \end{align*}
    since $e_0=0$. We use the geometric series to obtain $\Vert\Phi^{[K]}(y^0) - \widehat{\Phi}^{[K]}(y^0)\Vert_\mathcal{X} \leq \frac{\eta}{1-\kappa}$.
\end{proof}

To apply the Disturbance lemma, we need to consider a suitable closed subset while also verifying the assumptions regarding the operators $\Phi$ and $\Phi_{N,\textnormal{net}}$. 
We begin with showing that the set $\mathfrak{B}_N$, presented earlier, is a closed subset. We recall that $J\subset\mathbb{T}$ and $Y^\prime\in\mathbb{R}^+$.
\begin{lemma}
For $N\in\mathbb{N}$, the set 
\begin{align*}
    \mathfrak{B}_N=\{z\in L^2_N(J;\mathbb{R}^{n_y}) \ | \ \Vert z\Vert_{W^{1,\infty}(\mathbb{T};\mathbb{R}^{n_y})} \leq Y^\prime \}\subset C(J;\mathbb{R}^{n_y})
\end{align*}
is closed.
\end{lemma}
\begin{proof}
    The operator 
    \begin{align*}
        V:L^2_N(J;\mathbb{R}^{n_y})\rightarrow L^2_N(\mathbb{T};\mathbb{R}^{n_y}),\
        z\mapsto z
    \end{align*}
    is linear and thus continuous due to $L^2_N(J;\mathbb{R}^{n_y})$ and $L^2_N(\mathbb{T};\mathbb{R}^{n_y})$ being finite-dimensional. Also, we have that $L^2_N(\mathbb{T};\mathbb{R}^{n_y})\subset C(\mathbb{T};\mathbb{R}^{n_y})$        
    is closed due to the finite-dimensionality. Further, the set
    \begin{align*}
        A:=\{z\in C(\mathbb{T};\mathbb{R}^{n_y}) \ | \ \Vert z\Vert_{W^{1,\infty}(\mathbb{T};\mathbb{R}^{n_y})}\leq Y^\prime\}\subset C(\mathbb{T};\mathbb{R}^{n_y})
    \end{align*}
    is closed. Thus, we have that
    \begin{align*}
        \mathfrak{B}_0:=L^2_N(\mathbb{T};\mathbb{R}^{n_y})\cap A \subset C(\mathbb{T};\mathbb{R}^{n_y})
    \end{align*}
    is closed. We conclude that
    \begin{align*}
        \mathfrak{B}_N = V^{-1}(\mathfrak{B}_0)\subset C(J;\mathbb{R}^{n_y})
    \end{align*}
    is closed.
\end{proof}
The idea now is to apply the disturbance~\cref{Disturbance lemma} on our setting due to the following.
\begin{lemma}\label{application of disturbance lemma}
    Let~\cref{Final assumptions of ODE system} hold and consider $\Phi$ and $\Phi_{N,\textnormal{net}}$ defined as in~\cref{definition phi} and~\cref{definition phi n net}, while assuming that $\Phi_{N,\textnormal{net}}|_J$ fulfills the assumption in~\cref{Disturbance lemma}. Let $K\in\mathbb{N}$ and $y_{N,\textnormal{net}}^{(0)}\equiv0\in L^2_N(\mathbb{T};\mathbb{R}^{n_y}).$
    Then we have
    \begin{align*}
        \sup_{y^0\in Y^0, u\in\mathcal{U}^{l,J,A}}\Vert \hat{y}(\cdot;u,y^0) - \Phi_{N,\textnormal{net}}^{[K]}(y_{N,\textnormal{net}}^{(0)};u_N, y^0,\tau_N-\hat{t}_1)|_J \Vert_{C(J;\mathbb{R}^{n_y})} \leq e^{\lambda|J|}(\kappa^K Y + 2\eta),
    \end{align*}
    where $\hat{y}(\cdot;u,y^0)$ solves $\dot{y}=f(y,u)$ on $J$ with $y(\hat{t}_1)=y^0, u_N:=P_N(\tau,\mathcal{E}u)\in L^2_N(\mathbb{T};\mathbb{R}^{1+n_u})$ is the projected periodically control with first component $\tau_N=P_N\tau$. The value $e^{\lambda|J|}$ is to calculate between the sup-norm and the Bielecki norm.
\end{lemma}
\begin{proof}
    Let $y_{N,\textnormal{net}}^{(0)}\equiv0\in L^2_N(\mathbb{T};\mathbb{R}^{n_y})$. Define $v:=(u,y^0)$ and $v_N:=(u_N, y^0,\tau_N-\hat{t}_1)$. Note that $\Phi|_J$ defines a $\kappa$-contraction in the Bielecki norm. Thus in the Bielecki norm, denoted by $\Vert\cdot\Vert_X$ we obtain
    \begin{align*}
        &\Vert \hat{y} - \Phi_{N,\textnormal{net}}^{[K]}(y_{N,\textnormal{net}}^{(0)};v_N)|_J\Vert_X\\
        \leq&\Vert \hat{y} -  \Phi^{[K]}(y_{N,\textnormal{net}}^{(0)};v)|_J\Vert_X + \Vert \Phi^{[K]}(y_{N,\textnormal{net}}^{(0)};v)|_J - \Phi_{N,\textnormal{net}}^{[K]}(y_{N,\textnormal{net}}^{(0)};v_N)|_J \Vert_X\\
        \leq&\kappa^K \Vert \hat{y} - y_{N,\textnormal{net}}^{(0)}\Vert_X  + \frac{\eta}{1-\kappa}\\
        \leq&\kappa^K\Vert \hat{y}\Vert_X + 2\eta\\
        \leq&\kappa^K\Vert \hat{y}\Vert_{C(J;\mathbb{R}^{n_y})} + 2\eta \\
        \leq&\kappa^K Y+2\eta,
    \end{align*}
since $\kappa<\sfrac{1}{2}$.
\end{proof}
Before we verify $\Phi_{N,\textnormal{net}}|_J(\mathfrak{B}_N)\subset \mathfrak{B}_N$, we need to establish the bounds for the control functions, when they are not only periodically extended but also projected using
the $L^2$-orthogonal projection.
\begin{lemma}\label{bound of controls}
    For a function $u\in\mathcal{U}^{l,J,A}, l\geq2$ with $u_N:=P_N\mathcal{E} u$ it is 
    \begin{align*}
        \Vert u_N\Vert_{C(\mathbb{T};\mathbb{R}^{1+n_u})}\leq \Vert u_N\Vert_{W^{1,\infty}(\mathbb{T};\mathbb{R}^{1+n_u})}\leq U^\prime.
    \end{align*}
\end{lemma}
\begin{proof}
    It is 
    \begin{align*}
        \Vert u_N\Vert_{C(\mathbb{T};\mathbb{R}^{1+n_u})} &\leq \Vert u_N\Vert_{W^{1,\infty}(\mathbb{T};\mathbb{R}^{1+n_u})}\\
        &= \Vert P_N u\Vert_{W^{1,\infty}(\mathbb{T};\mathbb{R}^{1+n_u})}\\
        &\leq C_2\Vert P_N u\Vert_{H^{l}(\mathbb{T};\mathbb{R}^{1+n_u})}\\
        &\leq C_2\Vert u\Vert_{H^{l}(\mathbb{T};\mathbb{R}^{1+n_u})}\\
        &\leq C_2C_{\mathcal{E}u}\Vert u\Vert_{H^{l}(J;\mathbb{R}^{1+n_u})}\\
        &\leq C_2C_{\mathcal{E}u}\Bigg(A + \sqrt{\frac{\hat{t}_2^3-\hat{t}_1^3}{3}+\hat{t}_2-\hat{t}_1}\Bigg)\\
        &\leq C_u C_{\mathcal{E}u}\Bigg(A + \sqrt{\frac{\hat{t}_2^3-\hat{t}_1^3}{3}+\hat{t}_2-\hat{t}_1}\Bigg)+ 2\pi
        = U^\prime,
    \end{align*}
    where the first constant $C_2=C_2(W^{1,\infty}(\mathbb{T};\mathbb{R}^{1+n_u}),l)>0$ comes from the Sobolev embedding~\cref{Sobolev into continuous functions embedding}
    and the second constant $C_{\mathcal{E}u}>0$ comes from the continuous, linear periodic extension operator, see~\cref{periodic extension operator}. 
    Recalling~\cref{Definition of Cu} and~\cref{bound of periodicially extended controls} and we see that we are done. Note that we used that $P_N$ defines a projection.
\end{proof}
To show the aforementioned invariance of $\Phi_{N,\textnormal{net}}|_J$, we need to use an estimate for compositions of Sobolev functions, presented as next.
\begin{cor}\label{Sobolev composition estimate our case}
    Let $u_N$ with $u_N:=P_N u$ for $u\in \mathcal{U}^{l, \mathbb{T}, A} ,l\geq 2$ and $y_{N,\textnormal{net}}\in L^2_N(\mathbb{T};\mathbb{R}^{n_y})$ such that
    \begin{align*}
        \textnormal{Im}(y_{N,\textnormal{net}})\subset {(-Y^\prime,Y^\prime)}^{n_y}.
    \end{align*}
    Consider $f_\textnormal{net}:D^\prime\subset\mathbb{R}^d\rightarrow\mathbb{R}^{n_y}, d=1+n_y+n_u$, to be a ReLU or GELU neural network, i.e.\ $f_\textnormal{net}\in W^{1,\infty}(D^\prime;\mathbb{R}^{n_y})$. Then, there exists a constant, which we denote by $K_4=K_4(d,\mathbb{T})>0$, such that:
    \begin{align*}
    f_\textnormal{net}(y_{N,\textnormal{net}},u_N)\in W^{1,\infty}(\mathbb{T};\mathbb{R}^{n_y})\subset H^1(\mathbb{T};\mathbb{R}^{n_y})
    \end{align*}
    and 
    \begin{align*}
        \Vert f_\textnormal{net}(y_{N,\textnormal{net}},u_N)\Vert_{W^{1,\infty}(\mathbb{T};\mathbb{R}^{n_y})}
        \leq K_4\max\{\Vert f_\textnormal{net}\Vert_{L^\infty({D^\prime};\mathbb{R}^{n_y})}, \Vert f_\textnormal{net}\Vert_{W^{1,\infty}(D^\prime;\mathbb{R}^{n_y})}\Vert h\Vert_{W^{1,\infty}(\mathbb{T};\mathbb{R}^d)}\},
    \end{align*}
    where 
    \begin{align*}
        h:\mathbb{T}\rightarrow\mathbb{R}^d,\ t\mapsto(y_{N,\textnormal{net}}(t),u_N(t)).
    \end{align*}
\end{cor}
\begin{proof}
    Use Corollary B.5. from~\cite{guehring2019errorboundsapproximationsdeep} and consider there the mentioned domain $\Omega_1$ as $\Omega_1=\mathbb{T}\backslash\partial\mathbb{T}$. Then, the norms coincide
    \begin{align*}
        \Vert f_\textnormal{net}(y_{N,\textnormal{net}},u_N)\Vert_{W^{1,\infty}(\mathbb{T}\backslash\partial\mathbb{T};\mathbb{R}^{n_y})} 
      = \Vert f_\textnormal{net}(y_{N,\textnormal{net}},u_N)\Vert_{W^{1,\infty}(\mathbb{T};\mathbb{R}^{n_y})}.
    \end{align*}
    The inner functions are in $W^{1,\infty}(\mathbb{T};\mathbb{R}^d)$ since they are in $L^2_N(\mathbb{T};\mathbb{R}^d)\subset C^\infty(\mathbb{T};\mathbb{R}^d)$. Then use that the Sobolev semi-norm is smaller than the Sobolev norm.
\end{proof}
Last lemma allows us to give an estimate on the bound of the iterates, which then will lead automatically to the invariance of the operator $\Phi_{N,\textnormal{net}}|_J$.
\begin{lemma}\label{image of iterates}
    Consider $u\in\mathcal{U}^{l,\mathbb{T},A}, u_N:=P_N u, y^0\in Y^0, y_{N,\textnormal{net}}^{(0)}\in\mathfrak{B}_N$ and $f:D^\prime\subset\mathbb{R}^d\rightarrow\mathbb{R}^{n_y}, d=1+n_y+n_u, f\in W^{s,\infty}(D^\prime;\mathbb{R}^{n_y}), \Vert f\Vert_{W^{s,\infty}(D^\prime;\mathbb{R}^{n_y})}\leq E, s\geq 2$, with $f_\textnormal{net},\varepsilon_\textnormal{net}$ as in~\cref{choosage of fnet}.
    For the iterates of~\cref{our algorithm for CODE} it holds that 
    \begin{align*}
        \Vert y_{N,\textnormal{net}}^{(k)}\Vert_{W^{1,\infty}(\mathbb{T};\mathbb{R}^{n_y})} \leq Y^\prime,
    \end{align*}
    and thus
    \begin{align*}
        y_{N,\textnormal{net}}^{(k)}(t)\in {(-Y^\prime,Y^\prime)}^{n_y}
    \end{align*}
    for all $t\in\mathbb{T}$ and $k=0,\dots,K$. Recall that $Y^\prime = Y + \gamma>1$.
\end{lemma}
\begin{proof}
    It is $y_{N,\textnormal{net}}^{(0)}\in\mathfrak{B}_N$ and thus $\Vert y_{N,\textnormal{net}}^{(0)}\Vert_{W^{1,\infty}(\mathbb{T};\mathbb{R}^{n_y})} \leq Y^\prime$.
    So by induction we assume that $\Vert y_{N,\textnormal{net}}^{(k-1)}\Vert_{W^{1,\infty}(\mathbb{T};\mathbb{R}^{n_y})} \leq Y^\prime$, which implies
    \begin{align*}
         y_{N,\textnormal{net}}^{(k-1)}(t) \leq (-Y^\prime,Y^\prime)^{n_y}
    \end{align*}
    for all $t\in\mathbb{T}$.
    Therefore we can apply~\cref{Sobolev composition estimate our case} and obtain
    \begin{align*}
        &\Vert y^{(k)}_{N,\textnormal{net}}\Vert_{W^{1,\infty}(\mathbb{T};\mathbb{R}^{n_y})}
        \leq \Vert y^0\Vert_{\mathbb{R}^{n_y}} + \Vert \int^{\cdot}_{\hat{t}_1}\dot{P}_N\mathcal{I}_N f_\textnormal{net}({y}^{(k-1)}_{N,\textnormal{net}}, u_N)ds\Vert_{W^{1,\infty}(\mathbb{T};\mathbb{R}^{n_y})} + \vert c_{0,k}\vert\cdot\Vert \tau_N-\hat{t}_1\Vert_{W^{1,\infty}(\mathbb{T})},
    \end{align*}
    where the norm is taking over the function $\Big[t\mapsto \int^{t}_{\hat{t}_1}\dot{P}_N\mathcal{I}_N f_\textnormal{net}({y}^{(k-1)}_{N,\textnormal{net}}, u_N)ds\Big]$.
    We can further estimate
    \begin{align*}
        \Vert \int^{\cdot}_{\hat{t}_1}\dot{P}_N\mathcal{I}_N f_\textnormal{net}({y}^{(k-1)}_{N,\textnormal{net}}, u_N)ds\Vert_{W^{1,\infty}(\mathbb{T};\mathbb{R}^{n_y})}
        &\leq 2\pi\Vert\dot{P}_N\mathcal{I}_N f_\textnormal{net}(y_{N,\textnormal{net}}^{(k-1)},u_N)\Vert_{L^\infty(\mathbb{T};\mathbb{R}^{n_y})}\\
        &\leq 2\pi K_1\Vert\dot{P}_N\mathcal{I}_N f_\textnormal{net}(y_{N,\textnormal{net}}^{(k-1)},u_N)\Vert_{H^1(\mathbb{T};\mathbb{R}^{n_y})}\\
        &\leq 2\pi K_1\Vert\mathcal{I}_N f_\textnormal{net}(y_{N,\textnormal{net}}^{(k-1)},u_N)\Vert_{H^1(\mathbb{T};\mathbb{R}^{n_y})},
    \end{align*}
    and 
    \begin{align*}
        \vert c_{0,k}\vert\cdot\Vert \tau_N-\hat{t}_1\Vert_{W^{1,\infty}(\mathbb{T})}
        &\leq  U^\prime\fint_\mathbb{T}\vert\mathcal{I}_N f_\textnormal{net}(y_{N,\textnormal{net}}^{(k-1)},u_N)\vert ds\\
        &\leq  U^\prime 2\pi\Vert\mathcal{I}_N f_\textnormal{net}(y_{N,\textnormal{net}}^{(k-1)},u_N)\Vert_{L^\infty(\mathbb{T};\mathbb{R}^{n_y})}\\ 
        &\leq  U^\prime 2\pi K_1\Vert\mathcal{I}_N f_\textnormal{net}(y_{N,\textnormal{net}}^{(k-1)},u_N)\Vert_{H^1(\mathbb{T};\mathbb{R}^{n_y})},
    \end{align*}
    where we used that we can estimate $\Vert \tau_N-\hat{t}_1\Vert_{W^{1,\infty}(\mathbb{T})}$ also with $U^\prime\in\mathbb{R}^+$. Combining the estimates, we obtain
    \begin{align*}    
        \Vert y^{(k)}_{N,\textnormal{net}}\Vert_{W^{1,\infty}(\mathbb{T};\mathbb{R}^{n_y})}
        &\leq \Vert y^0\Vert_{\mathbb{R}^{n_y}}+(1+U^\prime)2\pi K_1\Vert\mathcal{I}_N f_\textnormal{net}(y_{N,\textnormal{net}}^{(k-1)},u_N)\Vert_{H^1(\mathbb{T};\mathbb{R}^{n_y})}\\
        &\leq \Vert y^0\Vert_{\mathbb{R}^{n_y}}+2\pi (1+U^\prime)K_1K_2\Vert f_\textnormal{net}(y_{N,\textnormal{net}}^{(k-1)},u_N)\Vert_{H^1(\mathbb{T};\mathbb{R}^{n_y})}\\
        &\leq \Vert y^0\Vert_{\mathbb{R}^{n_y}}+2\pi (1+U^\prime)K_1K_2K_3\Vert f_\textnormal{net}(y_{N,\textnormal{net}}^{(k-1)},u_N)\Vert_{W^{1,\infty}(\mathbb{T};\mathbb{R}^{n_y})}.
    \end{align*}
    We do a case distinction due to two possible upper bounds in~\cref{Sobolev composition estimate our case}. For the first case we have
    \begin{align*}
        &2\pi (1+U^\prime)K_1K_2K_3\Vert f_\textnormal{net}(y_{N,\textnormal{net}}^{(k-1)},u_N)\Vert_{W^{1,\infty}(\mathbb{T};\mathbb{R}^{n_y})}\\
        \leq&2\pi (1+U^\prime)K_1K_2K_3K_4\Vert f_\textnormal{net}\Vert_{W^{1,\infty}(D^\prime;\mathbb{R}^{n_y})}\Vert h\Vert_{W^{1,\infty}(\mathbb{T};\mathbb{R}^{d})}\\
        =&\widehat{K}\Vert f_\textnormal{net}\Vert_{W^{1,\infty}(D^\prime;\mathbb{R}^{n_y})}\Vert h\Vert_{W^{1,\infty}(\mathbb{T};\mathbb{R}^{d})}\\
        \leq&\widehat{K}\Vert f_\textnormal{net}\Vert_{W^{1,\infty}(D^\prime;\mathbb{R}^{n_y})}2\max\{U^\prime,Y^\prime\}\\
        \leq&\widehat{K}(\Vert f_\textnormal{net}-f\Vert_{W^{1,\infty}(D^\prime;\mathbb{R}^{n_y})} + \Vert f \Vert_{W^{1,\infty}(D^\prime);\mathbb{R}^{n_y}})2\max\{U^\prime,Y^\prime\}\\
        \leq&\widehat{K}(\Vert f_\textnormal{net}-f\Vert_{W^{1,\infty}(D^\prime;\mathbb{R}^{n_y})} + \Vert f \Vert_{W^{s,\infty}(D^\prime);\mathbb{R}^{n_y}})2\max\{U^\prime,Y^\prime\}\\
        \leq&\widehat{K}(\Vert f_\textnormal{net}-f\Vert_{W^{1,\infty}(D^\prime;\mathbb{R}^{n_y})} + E)2\max\{U^\prime,Y^\prime\}\\
        \leq&\widehat{K}(\frac{\gamma}{4\widehat{K}\max\{U^\prime,Y^\prime\}} + \frac{\gamma}{4\widehat{K}\max\{U^\prime,Y^\prime\}})2\max\{U^\prime,Y^\prime\}\\
        =&\gamma,
    \end{align*}
    where $h$ is given again by $h:\mathbb{T}\rightarrow\mathbb{R}^d,\
    t\mapsto(y_{N,\textnormal{net}}(t),u_N(t))$.
    For the other case we have
    \begin{align*}
        &2\pi (1+U^\prime)K_1K_2K_3\Vert f_\textnormal{net}(y_{N,\textnormal{net}}^{(k-1)},u_N)\Vert_{W^{1,\infty}(\mathbb{T};\mathbb{R}^{n_y})}\\
        \leq& 2\pi (1+U^\prime)K_1K_2K_3K_4\Vert f_\textnormal{net}\Vert_{L^{\infty}(D^\prime;\mathbb{R}^{n_y})}\\
        =&\widehat{K}\Vert f_\textnormal{net}\Vert_{L^{\infty}(D^\prime;\mathbb{R}^{n_y})}\\
        \leq&\widehat{K}\Vert f_\textnormal{net}\Vert_{W^{1,\infty}(D^\prime;\mathbb{R}^{n_y})}\\
        \leq&\widehat{K}(\Vert f_\textnormal{net}-f\Vert_{W^{1,\infty}(D^\prime;\mathbb{R}^{n_y})} + \Vert f \Vert_{W^{1,\infty}(D^\prime;\mathbb{R}^{n_y})})\\
        \leq&\widehat{K}(\Vert f_\textnormal{net}-f\Vert_{W^{1,\infty}(D^\prime;\mathbb{R}^{n_y})} + E)\\
        \leq&\frac{\gamma}{2\max\{U^\prime,Y^\prime\}}\\
        \leq&\frac{\gamma}{Y^\prime}\\
        \leq&{\gamma},
    \end{align*}
    since we choose $\gamma$ exactly like that, see \cref{choosage of gamma}.
    Besides, we used that $\dot{P}_N$ is a projection in the sense of~\cref{estimation for dot PN}.
    We also used that $\mathcal{I}_N$ is bounded, see~\cref{boundedness of IN}. Further we used the different Sobolev embeddings. 
    Now since $y^0\in{(-Y,Y)}^{n_y}$ and $Y^\prime:=Y+\gamma$, we are done.
\end{proof}
In the proof, we also showed the boundedness for the composition $f_\textnormal{net}(y, u_N)$ regarding $\Vert\cdot\Vert_{H^1(\mathbb{T};\mathbb{R}^{n_y})}$-norm, which is stated next as a corollary.
\begin{cor}\label{boundedness of f ReLU}
For $f_\textnormal{net}, y$ and $u_N$ as defined in~\cref{image of iterates}, we have that $f_\textnormal{net}(y, u_N)$ is bounded with respect to the $\Vert\cdot\Vert_{H^1(\mathbb{T};\mathbb{R}^{n_y})}$-norm for every
$y\in\mathfrak{B}_N$ and $u\in \mathcal{U}^{l,\mathbb{T},A}, u_N:=P_N u$. The upper bound, denoted by $\widehat{D}>0$, is dependent on $n_y,n_u,E, \alpha, m, s, l, \mathcal{E}, A, Y^0, \gamma$ and on the GELU network width $M_0$. 
\end{cor}
\begin{proof}
    See proof of~\cref{image of iterates}. There, one can see that the upper bound depends on $K_4$, which depends on $d=1+n_y+n_u$, and $f_\textnormal{net}, E$ and $h$. The network depends solely on $f$ while $h$ is dependent on $U^\prime$ and $Y^\prime$. The remaining dependence of $f$ is $\alpha, m$ and $s$, see~\cref{Final assumptions of ODE system}. The domain $U^\prime$ is dependent on $l, \mathcal{E}$ and $A$, see~\cref{bound of controls} for this dependence. The remaining dependence of the domain $Y^\prime$ is on $Y^0$ and $\gamma$, see~\cref{choosage of gamma}. When $f_\textnormal{net}=f_\textnormal{GELU}$, then we also have dependence on the fixed width $M_0$. Thus, we choose the constant $\widehat{D}>0$ as the maximum over both constants obtained either by $f_\textnormal{ReLU}$ or $f_\textnormal{GELU}$.
\end{proof}
We can now show that $\Phi_{N,\textnormal{net}}|_J$ verifies the first assumption from~\cref{Disturbance lemma}.
\begin{cor}
    Consider the closed subset
    \begin{align*}
        \mathfrak{B}_N=\{z\in L^2_N(J;\mathbb{R}^{n_y}) \ | \ \Vert z\Vert_{W^{1,\infty}(\mathbb{T};\mathbb{R}^{n_y})} \leq Y^\prime \}\subset C(J;\mathbb{R}^{n_y})
    \end{align*}
    and $\Phi_{N,\textnormal{net}}$ defined as in~\cref{definition phi n net} under the setting of~\cref{image of iterates}, we then have
    \begin{align*}
        \Phi_{N,\textnormal{net}}|_J(\mathfrak{B}_N)\subset \mathfrak{B}_N.
    \end{align*}
\end{cor}
\begin{proof}
Let $z\in \mathfrak{B}_N$, then $\Phi_{N,\textnormal{net}}|_J(z)$ is a trigonometric polynomial on $J$ with norm on $\mathbb{T}$ such that
    \begin{align*}
        \Vert \Phi_{N,\textnormal{net}}(z)\Vert_{W^{1,\infty}(\mathbb{T};\mathbb{R}^{n_y})}
        \leq Y^\prime,
    \end{align*}
    due to~\cref{image of iterates}.
\end{proof}
Following useful lemma can be found in the appendix of~\cite{guehring2019errorboundsapproximationsdeep} and will help us to show the second assumption in~\cref{Disturbance lemma}.
\begin{lemma}[Lipschitz characterization of Sobolev functions~\cite{guehring2019errorboundsapproximationsdeep}]\label{Lipschitz characterization of Sobolev functions}
    Let $D\subset\mathbb{R}^d, d\in\mathbb{N}$ be open and convex. If $f\in W^{1,\infty}(D)$, then $f$ is $\Vert\vert\nabla f(x)\vert\Vert_{L^\infty(D)}$-Lipschitz. 
    Conversely, if $f$ is bounded and $L$-Lipschitz on $D$, then $f\in W^{1,\infty}(D)$. 
\end{lemma}
\begin{proof}
    See~\cite{guehring2019errorboundsapproximationsdeep}, Theorem B.3.
\end{proof}
\begin{bmk}
    The characterization above is also true for vector-valued $W^{1,\infty}$ functions.
\end{bmk}
\begin{cor}
    Consider $f$ under~\cref{Final assumptions of ODE system}.
    There exists a Lipschitz-constant $L_u$ such that
    \begin{align*}
        \Vert f(y,u) - f(y,v)\Vert_{\mathbb{R}^{n_y}} \leq L_u\Vert u-v\Vert_{\mathbb{R}^{1+n_u}}
    \end{align*}
for all $u,v\in {(-U^\prime, U^\prime)}^{1+n_u}$.
\end{cor}
\begin{proof}
    This is simply due to the fact that now, due to the finite neural network approximation, we assume $f\in W^{s,\infty}(D^\prime;\mathbb{R}^{n_y})$. 
    Then using~\cref{Lipschitz characterization of Sobolev functions} gives us the statement.
\end{proof}
\begin{bmk}
Notice that the composition $f_\textnormal{net}(y_{N,\textnormal{net}}^{(k-1)}, u_N)$ is a periodic function by construction. 
Even further we can show that the composition is periodic up to its derivative.
This is due to the paper by~\cite{Murat2003} which states an analogue to the chain rule for the case when the outer function is a piecewise $C^1$ function. 
A ReLU neural network is a piecewise affine-linear function.
In the analogue chain rule, we see that the periodicity of the derivative only relies on the derivative of the inner functions, which are again, due to construction, periodic. 
Thus the composition is periodic up to its derivative. For the GELU function it is more immediate since can one use the standard chain rule.
\end{bmk}
We are left with proving the second assumption.
\begin{lemma}\label{difference phi and phi N net}
Consider the closed subset
    \begin{align*}
        \mathfrak{B}_N=\{z\in L^2_N(J;\mathbb{R}^{n_y}) \ | \ \Vert z\Vert_{W^{1,\infty}(\mathbb{T};\mathbb{R}^{n_y})} \leq Y^\prime \}\subset C(J;\mathbb{R}^{n_y}).
    \end{align*}
    Let $J=(\hat{t}_1,\hat{t}_2)\subset\mathbb{T}, N\in\mathbb{N}$ and $f_\textnormal{net},\varepsilon_\textnormal{net}$ as in~\cref{choosage of fnet}.
    For $\Phi$ and $\Phi_{N,\textnormal{net}}$ defined as in~\cref{definition phi} and~\cref{definition phi n net} there exists a constant $C^\prime=C^\prime(n_y,n_u,E, \alpha, m, s, l, \mathcal{E}, A, Y^0,\gamma, M_0, J)>0$ such that
    \begin{align*}
        \sup_{z\in \mathfrak{B}_N}\Vert \Phi(z;u,y^0)|_J - \Phi_{N,\textnormal{net}}(z;u_N,y^0,\tau_N-\hat{t}_1)|_J\Vert_{C(J;\mathbb{R}^{n_y})}\leq C^\prime (N^{-1} + \varepsilon_\textnormal{net}),
    \end{align*}
    where $y^0\in Y^0, u\in\mathcal{U}^{l,J,A}, u_N:=P_N(\tau,\mathcal{E}u)\in L^2_N(\mathbb{T};\mathbb{R}^{1+n_u})$ is the projected periodically extended control with first component $\tau_N=P_N\tau$.
    \begin{proof}
        Define $v:=(u,y^0)$ and $v_N:=(u_N,y^0,\tau_N-\hat{t}_1)$.
        For $z\in\mathfrak{B}_N$ and $t\in J=(\hat{t}_1,\hat{t}_2)\subset\mathbb{T}$ it is
        \begin{align*}
            &\Vert \Phi(z;v)[t] - \Phi_{N,\textnormal{net}}(z;v_N)[t] \Vert_{\mathbb{R}^{n_y}}\\
            =&\Vert \int^t_{\hat{t}_1} \Big[f(z,u) - \dot{P}_N\mathcal{I}_N f_\textnormal{net}(z,u_N)\Big]ds - c_0(\tau_N(t)-\hat{t}_1)\Vert_{\mathbb{R}^{n_y}}\\
            =&\Vert \int^t_{\hat{t}_1} \Big[f(z,u) - \mathcal{I}_N f_\textnormal{net}(z,u_N) + \big(\fint_\mathbb{T}\mathcal{I}_N f_\textnormal{net}(z,u_N)ds_2\big)\Big]ds_1 - c_0(\tau_N(t)-\hat{t}_1)\Vert_{\mathbb{R}^{n_y}}\\
            \leq&\int^t_{\hat{t}_1} \Vert f(z,u) - \mathcal{I}_N f_\textnormal{net}(z,u_N)\Vert_{\mathbb{R}^{n_y}} ds + \Vert (t-\hat{t}_1)\fint_\mathbb{T}\mathcal{I}_N f_\textnormal{net}(z,u_N)ds - c_0(\tau_N(t)-\hat{t}_1)\Vert_{\mathbb{R}^{n_y}}\\
            \leq&\underbrace{\int^t_{\hat{t}_1} \Vert f(z,u) - \mathcal{I}_N f_\textnormal{net}(z,u_N)\Vert_{\mathbb{R}^{n_y}} ds}_{=:
            (\mathrm{I})} + \vert c_0\vert\underbrace{\vert t-\tau_N(t)\vert}_{=:
            (\mathrm{II})},
        \end{align*}
        where we already know that $c_0$ is bounded independent of $z\in \mathfrak{B}_N$ and $u_N$ as seen in the proof of~\cref{image of iterates}.
        For
        $(\mathrm{I})$ we get
        \begin{align*}
            &\int^t_{\hat{t}_1} \Vert f(z,u) - \mathcal{I}_N f_\textnormal{net}(z,u_N)\Vert_{\mathbb{R}^{n_y}} ds\\
            \leq&\int_{0}^{2\pi}\Vert f(z,u) - f(z,u_N)\Vert_{\mathbb{R}^{n_y}} ds+\int_{\hat{t}_1}^{\hat{t}_2}\Vert f(z,u_N) - f_\textnormal{net}(z,u_N)\Vert_{\mathbb{R}^{n_y}} ds + \int_{0}^{2\pi}\Vert f_\textnormal{net}(z,u_N) - \mathcal{I}_N f_\textnormal{net}(z,u_N)\Vert_{\mathbb{R}^{n_y}} ds\\
            \leq& \int_{0}^{2\pi} L_u\Vert u - u_N\Vert_{\mathbb{R}^{n_u}} ds + |J|\cdot\Vert f-f_\textnormal{net}\Vert_{L^\infty(D^\prime)} + \sqrt{2\pi}\Vert f_\textnormal{net}(z,u_N) - \mathcal{I}_N f_\textnormal{net}(z,u_N)\Vert_{L^2(\mathbb{T};\mathbb{R}^{n_y})}\\
            \leq&\sqrt{2\pi}L_u \Vert u-u_N\Vert_{L^2(\mathbb{T};\mathbb{R}^{1+n_u})} + |J| \varepsilon_\textnormal{net} + \sqrt{2\pi}\Vert f_\textnormal{net}(z,u_N) - \mathcal{I}_N f_\textnormal{net}(z,u_N)\Vert_{L^2(\mathbb{T};\mathbb{R}^{n_y})}\\
            \leq&\sqrt{2\pi}L_u\widetilde{C}N^{-l}\Vert u\Vert_{H^l(\mathbb{T};\mathbb{R}^{1+n_u})} + |J| \varepsilon_\textnormal{net} + \sqrt{2\pi}\Vert f_\textnormal{net}(z,u_N) - \mathcal{I}_N f_\textnormal{net}(z,u_N)\Vert_{L^2(\mathbb{T};\mathbb{R}^{n_y})}\\
            \leq&\sqrt{2\pi}L_u \widetilde{C}N^{-l}C_{\tau,u} + |J| \varepsilon_\textnormal{net} + \sqrt{2\pi}\Vert f_\textnormal{net}(z,u_N) - \mathcal{I}_N f_\textnormal{net}(z,u_N)\Vert_{L^2(\mathbb{T};\mathbb{R}^{n_y})},
        \end{align*}
        with constant $\widetilde{C}=\widetilde{C}(l)>0$ coming from~\cref{Spectral approximation estimate} and where $C_{\tau,u}>0$ is the bound of the periodically extended controls with respect to the included time variation $\tau$, for which showed boundedness in $\Vert \cdot\Vert_{H^l(\mathbb{T})}$ at the end of~\cref{sec: Proof of well-posedness of CODE}. For the further estimation we do 
        \begin{align*}
            &\sqrt{2\pi}L_u \widetilde{C}N^{-l}C_{\tau,u} + |J| \varepsilon_\textnormal{net} + \sqrt{2\pi}\Vert f_\textnormal{net}(z,u_N) - \mathcal{I}_N f_\textnormal{net}(z,u_N)\Vert_{L^2(\mathbb{T};\mathbb{R}^{n_y})}\\
            \leq&\sqrt{2\pi}L_u \widetilde{C}N^{-l}C_{\tau,u} + |J| \varepsilon_\textnormal{net} + \sqrt{2\pi}\widetilde{D}N^{-1}\Vert f_\textnormal{net}(z,u_N)\Vert_{H^1(\mathbb{T};\mathbb{R}^{n_y})}\\
            \leq&\sqrt{2\pi}L_u \widetilde{C}N^{-l}C_{\tau,u} + |J| \varepsilon_\textnormal{net} + \sqrt{2\pi}\widetilde{D}\widehat{D} N^{-1},
        \end{align*}
        where the constant $\widetilde{D}>0$ comes from the~\cref{Pseudo-spectral approximation estimate} and $\widehat{D}=\widehat{D}(n_y,n_u,E, \alpha, m, s, l, \mathcal{E}, A, Y^0,\gamma, M_0)>0$ comes from the fact that the composition $f_\textnormal{net}(z,u_N)$ is bounded as stated in~\cref{boundedness of f ReLU}.
        For 
        $(\mathrm{II})$, we get for $t\in J$ and any $r\geq 1$ that
        \begin{align*}
            \vert t-\tau_N(t)\vert &\leq \Vert \text{id}-\tau_N\Vert_{C(J)}\\
                                   &\leq \Vert\tau -\tau_N\Vert_{C(\mathbb{T})}\\
                                   &\leq C_1\Vert\tau -\tau_N\Vert_{H^1(\mathbb{T})}\\
                                   &\leq C_1\Vert\tau -P_N\tau\Vert_{H^1(\mathbb{T})}\\
                                   &\leq C_1 D_1N^{-(r-1)}\Vert \tau\Vert_{H^r(\mathbb{T})},
        \end{align*}
        where $D_1=D_1(r)>0$ is the constant from~\cref{Spectral approximation estimate} and the last term $\Vert \tau\Vert_{H^r(\mathbb{T})}$ is also bounded as shown in~\cref{tau is smooth}.
        Note that the Lipschitz constant $L_u$, in the same way as $\widehat{D}$, is dependent on $f$, thus we can drop the explicit dependency on $L_u$, when stating the existence of a constant $C^\prime>0$. Since we can choose $r$ arbitrary, we also drop that dependency.
        In total, we can conclude that there exists a constant $C^\prime=C^\prime(n_y,n_u,E, \alpha, m, s, l, \mathcal{E}, A, Y^0,\gamma, M_0, J)>0$ such that
        \begin{align*}
            \sup_{z\in \mathfrak{B}_N}\Vert \Phi(z;u,y^0)|_J - \Phi_{N,\textnormal{net}}(z;u_N,y^0,\tau_N-\hat{t}_1)|_J\Vert_{C(J;\mathbb{R}^{n_y})}\leq C^\prime (N^{-1} + \varepsilon_\textnormal{net}).
        \end{align*}
    \end{proof}
\end{lemma}
Combining the last lemma with~\cref{application of disturbance lemma} gives us a convergence result for our proposed algorithm.
\begin{thm}\label{approximation by our algorithm}
Let~\cref{Final assumptions of ODE system} hold and consider $y_{N,\textnormal{net}}^{(0)}\equiv0\in\mathfrak{B}_N$. For $\Phi_{N,\textnormal{net}}$ defined as in~\cref{definition phi n net} with $K = \lceil\log(YN)/\log(1/\kappa)\rceil$ there are constants $C^*>0$ and $\widehat{C}>0$ depending on $n_y,n_u,E, \alpha, m, s, l, A, Y^0, \gamma, J$, on the fixed extension operator $\mathcal{E}$ and on the GELU width parameter $M_0$, but not on $\varepsilon$, such that for every $\varepsilon\in(0,\min\{\sfrac{1}{2},3C^*\min\{\sfrac{1}{2},E\}\})$ there exists a $N=N(\varepsilon)\leq \widehat{C}\varepsilon^{-1}$ with
    \begin{align*}
        \sup_{y^0\in Y^0, u\in\mathcal{U}^{l,J,A}}\Vert \hat{y}(\cdot;u,y^0) - \Phi_{N,\textnormal{net}}^{[K]}(y_{N,\textnormal{net}}^{(0)};u_N,y^0,\tau_N-\hat{t}_1)|_J \Vert_{C(J;\mathbb{R}^{n_y})} \leq \varepsilon,
    \end{align*}
    where $\hat{y}(\cdot;u,y^0)$ solves $\dot{y}=f(y,u)$ on $J$ with $y(\hat{t}_1)=y^0, u_N:=P_N(\tau,\mathcal{E}u)\in L^2_N(\mathbb{T};\mathbb{R}^{1+n_u})$ is the projected periodically extended control with first component $\tau_N=P_N\tau$. Moreover, for the obtained networks in this case, it holds
    \begin{align*}
        \textnormal{depth}(f_\textnormal{ReLU})\leq \widehat{C}\log({\varepsilon}^{-s/(s-1)}),\
        &\textnormal{width}(f_\textnormal{ReLU}) \leq \widehat{C}{\varepsilon}^{-d/(s-1)}\log({\varepsilon}^{-s/(s-1)}),\\
        \textnormal{depth}(f_\textnormal{GELU})\leq &\widehat{C} {\varepsilon}^{-d/(2(s-1))}\log({\varepsilon}^{-d/(2(s-1))}),
    \end{align*}
    where $s$ is the regularity of $f\in W^{s,\infty}(D^\prime;\mathbb{R}^{n_y})$ and $d=1+n_y+n_u$.
\end{thm}
\begin{proof}
    In~\cref{our algorithm for CODE} we choose 
    \begin{align*}
        K &= \lceil\log(YN)/\log(1/\kappa)\rceil,
    \end{align*}
    which implicates $\kappa^K\leq\sfrac{1}{YN}$ for $0<\kappa<1$. Recall that we chose $\varepsilon_\textnormal{net}$ (see~\cref{choosage of epsilon ReLU,choosage of epsilon GELU}) such that $\varepsilon_\textnormal{net}\leq N^{-1}$. The~\cref{application of disturbance lemma,difference phi and phi N net} immediately give us the existence of a constant $C^*=C^*(n_y,n_u,E, \alpha, m, s, l, \mathcal{E}, A, Y^0,\gamma, M_0, J)>0$ by defining $C^*=e^{\lambda \vert J\vert}\max\{1,2C^\prime\}$, such that
    \begin{align*}
        \sup_{y^0\in Y^0, u\in\mathcal{U}^{l,J,A}}\Vert \hat{y}(\cdot;u,y^0) - \Phi_{N,\textnormal{net}}^{[K]}(y_{N,\textnormal{net}}^{(0)};u_N,y^0,\tau_N-\hat{t}_1)|_J \Vert_{C(J;\mathbb{R}^{n_y})}\leq C^*(\kappa^K Y + N^{-1} + \varepsilon_\textnormal{net}).
    \end{align*}
    Now choosing $N = \lceil\frac{3C^*}{\varepsilon}\rceil$ implicates $\frac{3C^*}{\varepsilon}\leq N\leq \frac{3C^*}{\varepsilon}+1$ and
    \begin{align*}
            C^*(\kappa^K Y + N^{-1} + \varepsilon_\textnormal{net})\leq C^*(\frac{Y}{Y N} + N^{-1} + N^{-1})\leq C^*3N^{-1}\leq\varepsilon.
    \end{align*}
    For the size of the ReLU network with $\varepsilon_\textnormal{ReLU}=N^{-1}$,
    \begin{align*}
    \textnormal{depth}(f_\textnormal{ReLU})\leq c_\textnormal{ReLU}\cdot\log(\varepsilon_\textnormal{ReLU}^{-s/(s-1)}),\
    \textnormal{width}(f_\textnormal{ReLU}) \leq c_\textnormal{ReLU}\cdot\varepsilon_\textnormal{ReLU}^{-d/(s-1)}\cdot\log(\varepsilon_\textnormal{ReLU}^{-s/(s-1)}),
    \end{align*}
    it therefore follows for a constant $D_1$ depending on $C^*$ that
    \begin{align*}
        \textnormal{depth}(f_\textnormal{ReLU})\leq D_1\cdot\log({\varepsilon}^{-s/(s-1)}),\
        \textnormal{width}(f_\textnormal{ReLU}) \leq D_1\cdot{\varepsilon}^{-d/(s-1)}\cdot\log({\varepsilon}^{-s/(s-1)}).
    \end{align*}
    Recall that for the GELU network we set $L=\lceil(c_{0,\textnormal{GELU}}N)^{d/(2(s-1))}\rceil$ implying
    \begin{align*}
        L\leq(c_{0,\textnormal{GELU}}N)^{d/(2(s-1))}+1\leq(c_{0,\textnormal{GELU}}(\frac{3C^*}{\varepsilon}+1))^{{d/(2(s-1))}+1,}
    \end{align*}
    so that we can find a constant $D_2$ depending on $C^*$, and $c_{1,\textnormal{GELU}}$, such that
    \begin{align*}
        \textnormal{depth}(f_\textnormal{GELU}) \leq D_2 \cdot {\varepsilon}^{-d/(2(s-1))}\log({\varepsilon}^{-d/(2(s-1))}),
    \end{align*}
    since we chose the GELU neural network with constant width $M_0$. To conclude, we can now choose a constant $\widehat{C}=\widehat{C}(n_y,n_u,E, \alpha, m, s, l, \mathcal{E}, A, Y^0,\gamma, M_0, J)>0$ such that $\widehat{C}\geq\max\{D_1,D_2\}$ and also $N\leq \widehat{C}\varepsilon^{-1}$.
\end{proof}
We define the constant $\varepsilon_0:=\min\{\sfrac{1}{2},3C^*\min\{\sfrac{1}{2},E\}\}$ where $\varepsilon_0$ has the same dependencies as $C^*>0.$
\begin{bmk}\label{linear approximation lemma explanation}
    To show~\cref{main theorem}, we use the linear approximation~\cref{Linear approximation lemma} which is stated in the $L^2$-norm, but calculated for trigonometric polynomials. Due to the finite dimensionality of the space of trigonometric polynomials, we can reformulate the statement using the sup-norm but with a factor $D(N)>0$ depending on $N$ and of order $D(N)=\mathcal{O}(N^{1/2})$. The lemmas from~\cite{kovachki2021universalapproximationerrorbounds} that we use require smooth activation functions: \cref{Linear approximation lemma} assumes $\sigma\in C^2$ and~\cref{approximation of multiplication with nn} assumes $\sigma\in C^3$. We therefore use the GELU activation function in all layers that do not belong to $f_\textnormal{net}$. In particular, if $f$ is approximated by a GELU network, the whole $\Psi$-FNO uses a single activation function.
\end{bmk}
It remains to show, that we can emulate~\cref{our algorithm for CODE} using $\Psi$-FNO layers, which ultimately leads to proving~\cref{main theorem}.
\begin{proof}
    Consider $u_N:=P_N\mathcal{E}u$ for $u\in\mathcal{U}^{l,J,A}, l\geq 2, y^0\in Y^0$ and the projected periodically extended time variable $\tau_N=P_N\tau$. Then, under~\cref{Final assumptions of ODE system} we have the existence of a solution $\hat{y}(\cdot;u,y^0)$ which solves $\dot{y}=f(y,u)$ on $J$ with $\hat{y}(\hat{t}_1)=y^0$. 
    For the approximation by our proposed~\cref{our algorithm for CODE}, we choose $\varepsilon/2>0$ and $y_{N,\textnormal{net}}^{(0)}:=0\in\mathfrak{B}_N$ to apply~\cref{approximation by our algorithm}, which results in the existence of constants $C^*$ and $\widehat{C}$ depending on $n_y,n_u,E, \alpha, m, s, l, \mathcal{E}, A, Y^0,\gamma, M_0$ and $J$ while $\lceil\frac{6C^*}{\varepsilon}\rceil=N=N(\varepsilon/2)\leq \widehat{C}(\frac{\varepsilon}{2})^{-1}$ with
    \begin{align*}
        \sup_{y^0\in Y^0,u\in\mathcal{U}^{l,J,A}}\Vert \hat{y}(\cdot;u,y^0) - \Phi_{N,\textnormal{net}}^{[K]}(y_{N,\textnormal{net}}^{(0)};u_N,y^0,\tau_N-\hat{t}_1)|_J \Vert_{C(J;\mathbb{R}^{n_y})} \leq \varepsilon/2,
    \end{align*}
    and \begin{align*}
        \textnormal{depth}(f_\textnormal{ReLU})\leq \widehat{C}\log(({\frac{\varepsilon}{2}})^{-s/(s-1)}),\
        &\textnormal{width}(f_\textnormal{ReLU}) \leq \widehat{C}{({\frac{\varepsilon}{2}})}^{-d/(s-1)}\log({({\frac{\varepsilon}{2}})}^{-s/(s-1)}),\\
        \textnormal{depth}(f_\textnormal{GELU})\leq &\widehat{C} {({\frac{\varepsilon}{2}})}^{-d/(2(s-1))}\log({({\frac{\varepsilon}{2}})}^{-d/(2(s-1))}).
    \end{align*}
    By triangle inequality, we want to prove the existence of a $\Psi$-FNO $\widehat{\mathcal{N}}:L^2_N(\mathbb{T};\mathbb{R}^{(1+n_u)+n_y+2})\rightarrow L^2_N(\mathbb{T};\mathbb{R}^{n_y})$ such that
    \begin{align*}
        \sup_{y^0\in Y^0,u\in\mathcal{U}^{l,J,A}}\Vert \Phi_{N,\textnormal{net}}^{[K]}(y_{N,\textnormal{net}}^{(0)};u_N,y^0,\tau_N-\hat{t}_1) -\widehat{\mathcal{N}}\big(u_N,y^0,\tau_N-\hat{t}_1,D_N(\hat{t}_1-\cdot)\big)\Vert_{C(\mathbb{T};\mathbb{R}^{n_y})} \leq \varepsilon/2.    
    \end{align*}
    Note that we show the estimation on $\mathbb{T}$. To prove this estimation, we replace the steps of~\cref{our algorithm for CODE} with $\Psi$-FNO layers using either exact representation or the linear approximation~\cref{Linear approximation lemma}. Due to~\cref{linear approximation lemma explanation}, we set $\varepsilon_{\widehat{\mathcal{N}}}:=3C^*/D(N)N>0$, when using~\cref{Linear approximation lemma}.
    Define the lifting operator
    \begin{align*}
        {\mathcal{R}}:L^2_N(\mathbb{T};\mathbb{R}^{(1+n_u)+ n_y+ 2})&\rightarrow L^2_N(\mathbb{T}; \mathbb{R}^{(1+n_u)+ 2n_y+ 2})\\
        \big(u_N,y^0, \tau_N-\hat{t}_1, D_N(\hat{t}_1-\cdot)\big) &\mapsto \big(u_N, y^0, 0, \tau_N-\hat{t}_1, D_N(\hat{t}_1-\cdot)\big),
    \end{align*}
    which corresponds to the input step and step (3) in~\cref{our algorithm for CODE}. Step (1) is already fulfilled due to the assumption and step (2) will be resolved by following construction.
    Define the operators
    \begin{align*}
         {\mathcal{N}}_1:L^2_N(\mathbb{T};\mathbb{R}^{(1+n_u)+ 2n_y + 2})&\rightarrow L^2_N(\mathbb{T};\mathbb{R}^{(1+n_u)+ 2n_y + 2})\\
        \big(u_N, y^0, y, \tau_N-\hat{t}_1, D_N(\hat{t}_1-\cdot)\big) &\mapsto \big(u_N, y^0, \mathcal{I}_N f_\textnormal{net}(y, u_N), \tau_N-\hat{t}_1, D_N(\hat{t}_1-\cdot)\big),
    \end{align*}
    \begin{align*}
         {\mathcal{N}}_2:L^2_N(\mathbb{T};\mathbb{R}^{(1+n_u)+ 2n_y + 2})&\rightarrow L^2_N(\mathbb{T};\mathbb{R}^{(1+n_u)+ 3n_y + 2})\\
        \big(u_N, y^0, y, \tau_N-\hat{t}_1, D_N(\hat{t}_1-\cdot)\big) &\mapsto \big(u_N, y^0, \sum_{\vert k\vert\leq N}\delta_0(k)\mathcal{F}_N(y)(k)e^{ik\cdot}, \fint_\mathbb{T} y(s)ds, \tau_N-\hat{t}_1, D_N(\hat{t}_1-\cdot)\big),
    \end{align*}
    \begin{align*}
         {\mathcal{N}}_3:L^2_N(\mathbb{T};\mathbb{R}^{(1+n_u)+ 3n_y + 2})&\rightarrow L^2_N(\mathbb{T};\mathbb{R}^{(1+n_u)+ 4n_y + 2})\\
        \big(u_N, y^0, y, c_0, \tau_N-\hat{t}_1, D_N(\hat{t}_1-\cdot)\big) &\mapsto \big(u_N,y^0,y,\mathcal{I}_N\big(yD_N(\hat{t}_1-\cdot)\big), c_0 (\tau_N-\hat{t}_1),
        \tau_N-\hat{t}_1, D_N(\hat{t}_1-\cdot)\big),
    \end{align*} 
    \begin{align*}
        {\mathcal{N}}_4:L^2_N(\mathbb{T};\mathbb{R}^{(1+n_u)+ 4n_y + 2})&\rightarrow L^2_N(\mathbb{T};\mathbb{R}^{(1+n_u)+ 4n_y+ 2})\\
        \big(u_N, y^0, y, y_1, c_{0t}, \tau_N-\hat{t}_1, D_N(\hat{t}_1-\cdot)\big) &\mapsto \big(u_N,y^0,y, \fint_\mathbb{T} y_1(s)ds, c_{0t}, \tau_N-\hat{t}_1, D_N(\hat{t}_1-\cdot)\big),
    \end{align*} 
    and
    \begin{align*}
         {\mathcal{N}}_5:L^2_N(\mathbb{T};\mathbb{R}^{(1+n_u)+ 4n_y + 2})&\rightarrow L^2_N(\mathbb{T};\mathbb{R}^{(1+n_u)+ 2n_y + 2})\\
        \big(u_N, y^0, y, y_1, c_{0t}, \tau_N-\hat{t}_1, D_N(\hat{t}_1-\cdot)\big) &\mapsto \big(u_N, y^0, y^0 + y - y_1 + c_{0t}, \tau_N-\hat{t}_1, D_N(\hat{t}_1-\cdot)\big).
    \end{align*}
    Define the projection operator
    \begin{align*}
         {\mathcal{Q}}:L^2_N(\mathbb{T};\mathbb{R}^{(1+n_u)+ 2n_y + 2})&\rightarrow L^2_N(\mathbb{T};\mathbb{R}^{n_y})\\
        (u_N, y^0, y, \tau_N-\hat{t}_1, D_N)&\mapsto y,
    \end{align*}
    to see that our proposed~\cref{our algorithm for CODE} can now be written as
    \begin{align*}
    \mathcal{N}:L^2_N(\mathbb{T};\mathbb{R}^{(1+n_u)+ n_y+2})&\rightarrow L^2_N(\mathbb{T};\mathbb{R}^{n_y})\\
        v_N&\mapsto  \bigg[{\mathcal{Q}}\circ \underbrace{ \mathcal{N}^\prime\circ\dots\circ \mathcal{N}^\prime}_{K\text{-fold composition}}\circ {\mathcal{R}}\bigg]\big(u_N, y^0, \tau_N-\hat{t}_1, D_N(\hat{t}_1-\cdot)\big),
    \end{align*}
    where \begin{align*}
        \mathcal{N}^\prime := {\mathcal{N}}_5\circ{\mathcal{N}}_4\circ {\mathcal{N}}_3\circ {\mathcal{N}}_2\circ {\mathcal{N}}_1,\
        K &= \lceil\log(YN)/\log(1/\kappa)\rceil.
    \end{align*}
    For each operator $\mathcal{N}_i, i=1,\dots,5$, we argue that we can represent it exactly or can at least find an approximation by a $\Psi$-FNO layer. Note that when approximating one of the operators $\mathcal{N}_i, i=1,\dots,5$ we introduce an error which propagates through the iteration. This causes no issue, as we can simply apply the Replacement~\cref{Replacement lemma}.
    According to~\cref{Psi FNO representation of I_N NN} the operator ${\mathcal{N}}_1$ can be represented by a $\Psi$-FNO $\widehat{\mathcal{N}}_1$ with
    \begin{align*}
    \textnormal{width}(\widehat{\mathcal{N}}_1) = \vert \mathcal{J}_N\vert \textnormal{width}(f_\textnormal{net}), \textnormal{depth}(\widehat{\mathcal{N}}_1) = \textnormal{depth}(f_\textnormal{net}), \textnormal{lift}(\widehat{\mathcal{N}}_1) = \textnormal{width}(f_\textnormal{net}).
    \end{align*}
The channels $u_N, y^0, \tau_N-\hat{t}_1$ and $D_N(\hat{t}_1-\cdot)$ are passed through the $\textnormal{depth}(f_\textnormal{net})$ layers of $\widehat{\mathcal{N}}_1$ exactly, since $\sigma(x)-\sigma(-x)=x$ holds for the ReLU and the (standard) GELU activation function. This changes width and lift of $\widehat{\mathcal{N}}_1$ only by a constant factor.

To see that we can represent $ {\mathcal{N}}_2$ using $\Psi$-FNO structure, recall that~\cref{Integral Discrete Fourier Transform Formula} gives us
    \begin{align*}
        \sum_{\vert k\vert\leq N}\delta_0(k)\mathcal{F}_N(y)(k) e^{it_jk} = \mathcal{F}_N^{-1}(\delta_0(k)\mathcal{F}_N(y)(k))(j),\ \forall j\in\mathcal{J}_N,
    \end{align*}
and, implied by~\cref{mean repr using psi fno}, that
    \begin{align*}
        \fint^{2\pi}_{0} y(s)ds = \mathcal{F}_N^{-1}(\delta_1(k)\mathcal{F}_N(y)(k))(j),\ \forall j\in\mathcal{J}_N,
    \end{align*}
where
        \begin{align*}
            \delta_0(k):=\begin{cases}
                            \frac{1}{ik},\quad k\not=0,\\
                            0,\quad k=0,
            \end{cases} \ \textnormal{and} \ \delta_1(k):=\begin{cases}
                            0,\quad k\not=0,\\
                            1,\quad k=0.
            \end{cases}
    \end{align*}
Both right-hand side already have a similar structure of a linear $\Psi$-FNO layer, or also called $\mathcal{F}$-layer, as defined in~\cref{Psi Fno linear layer}.
We now introduce the precise weight matrix, bias and Fourier multiplier matrix to show that we can write $\mathcal{N}_2$ as a single linear $\mathcal{F}$-layer.
Set the weight matrix, bias and Fourier multiplier matrix as 
    \begin{align*}
    W_2:=\begin{pmatrix}
        1 &0 &0 &0 &0\\
        0 &1 &0 &0 &0\\
        0 &0 &0 &0 &0\\
        0 &0 &0 &0 &0\\
        0 &0 &0 &1 &0\\
        0 &0 &0 &0 &1
    \end{pmatrix}, {b_2}_j:=\begin{pmatrix}0\\0\\0\\0\\0\\0\end{pmatrix}, P_2(k):=\begin{pmatrix}
        0 &0 &0 &0 &0\\
        0 &0 &0 &0 &0\\
        0 &0 &\delta_0(k) &0 &0\\
        0 &0 &\delta_1(k) &0 &0\\
        0 &0 &0 &0 &0\\
        0 &0 &0 &0 &0\\
    \end{pmatrix},
    \end{align*}
then
    \begin{align*}
        \mathcal{N}_2{(v_N)}_j = W_2{v_N}_j + {b_2}_j + \mathcal{F}_N^{-1}{(P_2\mathcal{F}_N v)}_j,\ \forall j\in\mathcal{J}_N,
    \end{align*}
where $v_N:=\big(u_N, y^0, y, \tau_N-\hat{t}_1, D_N(\hat{t}_1-\cdot)\big)\in L^2_N(\mathbb{T};\mathbb{R}^{(1+n_u)+ 2n_y+ 2})$.
To apply the linear approximation~\cref{Linear approximation lemma}, we choose $\varepsilon_{\widehat{\mathcal{N}}}>0$ as mentioned and we choose as activation function the GELU activation function, see~\cref{linear approximation lemma explanation}. Further to apply the lemma, we need to show that $\Vert v_N\Vert_{L^2}\leq B$ for some $B>0$.
This follows easily since $u_N$ is bounded in the $L^2$-norm due to~\cref{bound of controls}. Further we have  $y^0\in Y^0, y\in\mathfrak{B}_N, \tau_N-\hat{t}_1\in C(\mathbb{T})$ and $D_N\in L^2_N(\mathbb{T})$.
Thus, we can find a common $L^2$-bound $B>0$. By the linear approximation lemma, one obtains a constant $C_1>0$ and $\Psi$-FNO $\widehat{\mathcal{N}}_2$ with
    \begin{align*}
        \textnormal{width}(\widehat{\mathcal{N}}_2)&\leq C_1N,\
        \textnormal{depth}(\widehat{\mathcal{N}}_2)\leq C_1,\
        \textnormal{lift}(\widehat{\mathcal{N}}_2)\leq  C_1
    \end{align*}
and
    \begin{align*}
        \sup_{\Vert v_N\Vert_{L^2}\leq B}\Vert\mathcal{N}_2(v_N)-\widehat{\mathcal{N}}_2(v_N)\Vert_{L^2}\leq\varepsilon_{\widehat{\mathcal{N}}}.
    \end{align*}
For $\mathcal{N}_3$ we use that the multiplication can be approximated by finite neural networks due to~\cref{approximation of multiplication with nn}. As an activation function,
we again choose the GELU and set $\varepsilon_{\widehat{\mathcal{N}}}/\sqrt{4\pi}>0$. Notice that we need to apply the lemma two times. Set $z(t):= \sum_{\vert k\vert\leq N}\delta_0(k)\mathcal{F}_N(y)(k)e^{ikt}, t\in\mathbb{T}$, and $c_0:=\fint_\mathbb{T} y(s)ds\in\mathbb{R}^{n_y}$. Due to the finite dimensionality of the space of trigonometric polynomials it is 
    \begin{align*}
        \sup_{t\in\mathbb{T}}\vert D_N(\hat{t}_1-t)\vert &\leq B^\prime \Vert D_N(\hat{t}_1-\cdot)\Vert_{L^2}\leq B^{\prime\prime},
    \end{align*}
where $B^\prime=B^\prime(N)>0$ and $B^{\prime\prime}=B^{\prime\prime}(N)>0$. Then for $t\in\mathbb{T}$ we have
    \begin{align*}
        \Vert z(t)\Vert_{\mathbb{R}^{n_y}}&\leq Y^\prime \sum_{0<\vert k\vert\leq N}\frac{\vert e^{ikt}\vert}{\vert k\vert} \leq 2Y^\prime(1+\log(N))\\
        \Rightarrow\sup_{t\in\mathbb{T}}\Vert z(t)\Vert_{\mathbb{R}^{n_y}}&\leq 2Y^\prime(1+\log(N)).
    \end{align*}
It also is
    \begin{align*}
        \sup_{t\in\mathbb{T}}\vert\tau_N(t)-\hat{t}_1 \vert\leq B^{\prime\prime\prime},\
        \sup_{t\in\mathbb{T}}\Vert c_0\Vert_{\mathbb{R}^{n_y}}\leq 2\pi Y^\prime,
    \end{align*}
where $B^{\prime\prime\prime}=B^{\prime\prime\prime}(N)>0$. By choosing the maximum out of all bounds we have a bound $\widehat{B}=\widehat{B}(N)>0$ such that we obtain the existence of constants $C_2=C_2(n_y)>0$ and $C_3=C_3(n_y)>0$ and finite GELU neural networks $\phi_1, \phi_2:\mathbb{R}\times\mathbb{R}^{n_y}\rightarrow\mathbb{R}^{n_y}$ for the multiplication of a scalar and a vector, see~\cref{approximation of multiplication with nn}~(2), such that
    \begin{align*}
        \sup_{\vert D_N\vert, \Vert z\Vert_{\mathbb{R}^{n_y}}\leq \widehat{B}}\Vert
        \phi_1(D_N, z) - D_N z\Vert_{\mathbb{R}^{n_y}}&\leq{\varepsilon_{\widehat{\mathcal{N}}}}/{\sqrt{4\pi}}\\
        \sup_{\vert \tau_N-\hat{t}_1\vert, \Vert c_0\Vert_{\mathbb{R}^{n_y}}\leq \widehat{B}}\Vert 
        \phi_2(\tau_N-\hat{t}_1, c_0) - c_0(\tau_N-\hat{t}_1)\Vert_{\mathbb{R}^{n_y}}&\leq{\varepsilon_{\widehat{\mathcal{N}}}}/{\sqrt{4\pi}}
    \end{align*}
    where the size is constant with $\textnormal{width}(\phi_1), \textnormal{depth}(\phi_1)\leq C_2$ and $\textnormal{width}(\phi_2), \textnormal{depth}(\phi_2)\leq C_3$. 
    Then, we again use~\cref{Psi FNO representation of I_N NN} just like for the operator $\mathcal{N}_1$ to obtain the existence
    of $\Psi$-FNOs $\widehat{\mathcal{N}}_{3,i},i=1,2$, of constant size with
    \begin{align*}
&\max_{j\in\mathcal{J}_N}\Vert \widehat{\mathcal{N}}_{3,1}\big(z, D_N(\hat{t}_1-\cdot)\big)(t_j) - z(t_j)D_N(\hat{t}_1-t_j)\Vert_{\mathbb{R}^{n_y}}     \leq{\varepsilon_{\widehat{\mathcal{N}}}}/{\sqrt{4\pi}},\\
&\max_{j\in\mathcal{J}_N}\Vert \widehat{\mathcal{N}}_{3,2}(c_0, \tau_N-\hat{t}_1)(t_j) - c_0(\tau_N(t_j)-\hat{t}_1)\Vert_{\mathbb{R}^{n_y}}
\leq{\varepsilon_{\widehat{\mathcal{N}}}}/{\sqrt{4\pi}}.
    \end{align*}
To obtain one single $\Psi$-FNO $\widehat{\mathcal{N}}_3$ we just choose again a weight matrix, bias and Fourier multiplier matrix appropriately. 
Both estimates hold at the grid points. Since $\widehat{\mathcal{N}}_{3,1}(z,D_N(\hat{t}_1-\cdot))$ and $\mathcal{I}_N\big(zD_N(\hat{t}_1-\cdot)\big)$, as well as $\widehat{\mathcal{N}}_{3,2}(c_0,\tau_N-\hat{t}_1)$ and $c_0(\tau_N-\hat{t}_1)$, are trigonometric polynomials of degree at most $N$, the discrete Parseval identity $\Vert p\Vert_{L^2(\mathbb{T})}^2=\frac{2\pi}{2N+1}\sum_{j\in\mathcal{J}_N}\Vert p(t_j)\Vert^2$ yields $L^2$-estimates with upper bound $\varepsilon_{\widehat{\mathcal{N}}}/\sqrt{2}$ for each of the two output channels and hence with upper bound $\varepsilon_{\widehat{\mathcal{N}}}$ for $\widehat{\mathcal{N}}_3$, as needed for the Replacement lemma.
    The bounds $B$ and $\widehat{B}$ used for $\mathcal{N}_2$, $\mathcal{N}_3$ and $\mathcal{N}_4$ are not restricted to the exact Picard iterates: for every $M>0$ and every input $v\in L^2_N(\mathbb{T};\mathbb{R}^m)$ with $\Vert v\Vert_{L^2}\leq M$, the inverse estimate $\Vert v\Vert_{L^\infty}\leq\sqrt{(2N+1)/(2\pi)}\,\Vert v\Vert_{L^2}$ yields point-wise bounds depending only on $M$ and $N$. Since the sizes of the networks in~\cref{Linear approximation lemma} and~\cref{approximation of multiplication with nn} do not depend on these bounds, every operator $\mathcal{N}_i$ can be approximated uniformly on $L^2$-bounded sets by $\Psi$-FNOs of fixed size, as required in~\cref{Replacement lemma}. In particular, this covers inputs that are perturbed by the approximation errors of preceding layers.
    The operator $\mathcal{N}_4$ is resolved exactly the same way as operator $\mathcal{N}_2$. Note that the product of $z$ and $D_N(\hat{t}_1-\cdot)$ lies in the space $L^2_{2N}$ while the $\Psi$-FNO calculates $\fint_\mathbb{T}\mathcal{I}_N(z(s) D_N(\hat{t}_1-s))ds$. This causes no issue as we have $\fint_\mathbb{T}\mathcal{I}_N(z(s) D_N(\hat{t}_1-s))ds = \fint_\mathbb{T}z(s) D_N(\hat{t}_1-s)ds$, due to~\cref{discrete and continuous mean align}. Thus, the result from~\cref{Dirichlet kernel property} implies
    \begin{align*}
        z(\hat{t}_1) = \fint_\mathbb{T} z(s)D_N(\hat{t}_1-s)ds,
    \end{align*}
    where here it is $z(\hat{t}_1)=\sum_{\vert k\vert\leq N}\delta_0(k)\mathcal{F}_N(y)(k)e^{ik\hat{t}_1}$. To apply again the linear approximation lemma,
    we need to show that $\mathcal{I}_N\big(zD_N(\hat{t}_1-\cdot)\big)$ and $c_0(\tau_N-\hat{t}_1)$ are bounded in $L^2$-norm, which is immediate due to the previous arguments. Also, the size of the obtained neural operator $\widehat{\mathcal{N}}_4$ only differs by a positive constant compared to the size of $\widehat{\mathcal{N}}_2$. The mapping $\mathcal{N}_5$ can be written as a linear $\Psi$-FNO layer with weight matrix, bias and Fourier multiplier
    \begin{align*}
        W_5:=\begin{pmatrix}
        1 &0 &0 &0 &0 &0 &0\\
        0 &1 &0 &0 &0 &0 &0\\
        0 &1 &1 &-1 &1 &0 &0\\
        0 &0 &0 &0 &0 &1 &0\\
        0 &0 &0 &0 &0 &0 &1\\
    \end{pmatrix}, {b_5}_j:=\begin{pmatrix}0\\0\\0\\0\\0\end{pmatrix}, P_5(k):= \underline{0}.
    \end{align*}
    Remark that the calculated function $y-y_1$ in $\mathcal{N}_5$ is $\sum_{\vert k\vert\leq N}\delta_0(k)\mathcal{F}_N(\tilde{y})(k)(e^{ikt}
    -e^{ik\hat{t}_1})$, where it is $\tilde{y}:=\mathcal{I}_Nf_\textnormal{net}(0,u_N)$ for the first iteration. Therefore, due to~\cref{Integral Discrete Fourier Transform Formula} we calculate $\int^t_{\hat{t}_1}(\dot{P}_N\mathcal{I}_Nf_\textnormal{net}(0,u_N))(s)ds$ as required for the Picard type iteration.
    Again, to apply~\cref{Linear approximation lemma} we argue that the input is bounded which is immediate. The obtained neural operator $\widehat{\mathcal{N}}_5$ is of same size up to a constant as $\widehat{\mathcal{N}}_2$. Define
    \begin{align*}
        \mathcal{N}^\dagger = \underbrace{\mathcal{N}_5\circ\mathcal{N}_4\circ\mathcal{N}_3\circ\mathcal{N}_2\circ\mathcal{N}_1\circ\dots\circ\mathcal{N}_5\circ\mathcal{N}_4\circ\mathcal{N}_3\circ\mathcal{N}_2\circ\mathcal{N}_1}_{K\text{-fold composition}}
    \end{align*}
    to perform $K$-times iteration as defined in~\cref{our algorithm for CODE}.
    While not all operators $\mathcal{N}_i, i=1,\dots,5$ can be represented exactly by a $\Psi$-FNO but rather rely on approximation, we apply the Replacement~\cref{Replacement lemma} on $\mathcal{N}^\dagger$ to obtain a single $\Psi$-FNO $\widehat{\mathcal{N}}^\dagger$ with
    \begin{align*}
        \sup_{\Vert v_N\Vert_{L^2}\leq B^\dagger}\Vert \mathcal{N}^\dagger(v_N) - \widehat{\mathcal{N}}^\dagger(v_N)\Vert_{L^2}\leq \varepsilon_{\widehat{\mathcal{N}}},
    \end{align*}
    where $B^\dagger$ is chosen as the maximum out of all occurring $L^2$-input bounds. For the size of $\widehat{\mathcal{N}}^\dagger$ it is $\textnormal{depth}(\widehat{\mathcal{N}}^\dagger)\leq 5K\cdot \max_{i=1,\dots,5}\{\textnormal{depth}(\widehat{\mathcal{N}}_i)\}, \textnormal{width}(\widehat{\mathcal{N}}^\dagger)\leq\max_{i=1,\dots,5}\{\textnormal{width}(\widehat{\mathcal{N}}_i)\}$ and $\textnormal{lift}(\widehat{\mathcal{N}}^\dagger)\leq \max_{i=1,\dots,5}\{\textnormal{lift}(\widehat{\mathcal{N}}_i)\}$. 
    We define
    \begin{align*}
        \widehat{\mathcal{N}}(v_N) := (\mathcal{Q}\circ\widehat{\mathcal{N}}^\dagger\circ\mathcal{R})(v_N)  
    \end{align*}
    and since
    \begin{align*}
        \Phi_{N,\textnormal{net}}^{[K]}(y_{N,\textnormal{net}}^{(0)};u_N, y^0,\tau_N-\hat{t}_1)=\mathcal{N}(v_N)=(\mathcal{Q}\circ\mathcal{N}^\dagger\circ\mathcal{R})(v_N),
    \end{align*}
    we conclude 
    \begin{align*}
        \sup_{\Vert v_N\Vert_{L^2}\leq B^\dagger}\Vert \Phi_{N,\textnormal{net}}^{[K]}(y_{N,\textnormal{net}}^{(0)};u_N,y^0,\tau_N-\hat{t}_1) - \widehat{\mathcal{N}}(v_N)\Vert_{L^2(\mathbb{T};\mathbb{R}^{n_y})}\leq\varepsilon_{\widehat{\mathcal{N}}},
    \end{align*}
    where we recall that $v_N:=\big(u_N, y^0, \tau_N-\hat{t}_1, D_N(\hat{t}_1-\cdot)\big)$.
    Due to the finite dimensionality of $L^2_N$, see also~\cref{linear approximation lemma explanation}, we obtain the existence of a constant $D=D(N)>0$ such that 
    \begin{align*}
        \sup_{\Vert v_N\Vert_{L^2}\leq B^\dagger}\Vert \Phi_{N,\textnormal{net}}^{[K]}(y_{N,\textnormal{net}}^{(0)};u_N,y^0,\tau_N-\hat{t}_1) - \widehat{\mathcal{N}}(v_N)\Vert_{C(\mathbb{T};\mathbb{R}^{n_y})}\leq D(N)\varepsilon_{\widehat{\mathcal{N}}} =D(N)\frac{3C^*}{D(N)N} = 3C^* N^{-1} \leq \varepsilon/2,
    \end{align*}
    where $N=\lceil\frac{6C^*}{\varepsilon}\rceil$ and $\varepsilon_{\widehat{\mathcal{N}}}=3C^*/D(N)N$. Since our input always fulfills $\Vert v_N\Vert_{L^2}\leq B^\dagger$, we have shown the estimate formulated in the statement.
    For the size of $\widehat{\mathcal{N}}$ we do following estimations first
    \begin{align*}
        &\textnormal{depth}(\widehat{\mathcal{N}})\leq 5K\cdot\max_{i=1,\dots,5}\{\textnormal{depth}(\widehat{\mathcal{N}}_i)\}\lesssim K\cdot\textnormal{depth}(f_\textnormal{net}),\\
        &\textnormal{width}(\widehat{\mathcal{N}})\leq \max_{i=1,\dots,5}\{\textnormal{width}(\widehat{\mathcal{N}}_i)\},\\
        &\textnormal{lift}(\widehat{\mathcal{N}})\leq \max_{i=1,\dots,5}\{\textnormal{lift}(\widehat{\mathcal{N}}_i)\}\lesssim \textnormal{width}(f_\textnormal{net}).
    \end{align*}
    In case of $f_\textnormal{net}=f_\textnormal{ReLU}$, it is then
    \begin{align*}
        &\textnormal{depth}(\widehat{\mathcal{N}}_\textnormal{ReLU})\lesssim \log(\frac{1}{\varepsilon})\log({\varepsilon}^{-s/(s-1)}),\\
        &\textnormal{lift}(\widehat{\mathcal{N}}_\textnormal{ReLU})\leq \max_{i=1,\dots,5}\{\textnormal{lift}(\widehat{\mathcal{N}}_i)\}\lesssim \textnormal{width}(f_\textnormal{ReLU})\lesssim {\varepsilon}^{-d/(s-1)}\cdot\log({\varepsilon}^{-s/(s-1)}),\\
        &\textnormal{width}(\widehat{\mathcal{N}}_\textnormal{ReLU})\lesssim \max\{\textnormal{width}(\widehat{\mathcal{N}}_1), \textnormal{width}(\widehat{\mathcal{N}}_2)\}\\ 
        &\qquad\leq\max\{\vert\mathcal{J}_N\vert \textnormal{width}(f_\textnormal{ReLU}), C_1N\}
        \lesssim \varepsilon^{-1}{\varepsilon}^{-d/(s-1)}\cdot\log({\varepsilon}^{-s/(s-1)}),
    \end{align*}
    and in case of $f_\textnormal{net}=f_\textnormal{GELU}$, we have
    \begin{align*}
        &\textnormal{depth}(\widehat{\mathcal{N}}_\textnormal{GELU})\lesssim \log(\frac{1}{\varepsilon}){\varepsilon}^{-d/(2(s-1))}\log({\varepsilon}^{-d/(2(s-1))}),\\
        &\textnormal{width}(\widehat{\mathcal{N}}_\textnormal{GELU})\lesssim N\lesssim {(\frac{1}{\varepsilon})},\\
        &\textnormal{lift}(\widehat{\mathcal{N}}_\textnormal{GELU})\leq \max_{i=1,\dots,5}\{\textnormal{lift}(\widehat{\mathcal{N}}_i)\}\lesssim \textnormal{width}(f_\textnormal{GELU}),
    \end{align*}
    where $\textnormal{width}(f_\textnormal{GELU})=c_{2,\textnormal{GELU}}M_0\log M_0$ is constant, since the width parameter $M_0\in\mathbb{N}$ was fixed beforehand.
    Recalling~\cref{estimate for size of psi fno} and choosing a constant $C>0$ in a similar fashion as in the proof of~\cref{approximation by our algorithm} we obtain
    \begin{align*}
        \textnormal{size}(\widehat{\mathcal{N}}_\textnormal{ReLU}) \leq C\cdot {(\frac{1}{\varepsilon})}^{\frac{2d}{s-1}+1}{\big(\log(\frac{1}{\varepsilon})\big)}^4,\
        \textnormal{size}(\widehat{\mathcal{N}}_\textnormal{GELU}) \leq C\cdot{(\frac{1}{\varepsilon})}^{\frac{d}{2(s-1)}+1} {\big(\log(\frac{1}{\varepsilon})\big)}^2,
    \end{align*}
    and $N\leq C\varepsilon^{-1}$. Since $C$ depends on $C^*$ and $\widehat{C}$, it has the same dependencies and thus we have proven our main result.
\end{proof}

%% file: references.bib
@article{RATHABA200389,
title = {Identifiability of an {Electric Arc Furnace Model}},
journal = {IFAC Proceedings Volumes},
volume = {36},
number = {24},
pages = {89-94},
year = {2003},
note = {IFAC Workshop on New Technologies for Automation of the Metallurgical Industry 2003, Shanghai, PR China, 11-13 October 2003},
issn = {1474-6670},
doi = {10.1016/S1474-6670(17)37609-7},
url = {https://www.sciencedirect.com/science/article/pii/S1474667017376097},
author = {P.L. Rathaba and I.K. Craig and X. Xia and P.C. Pistorius}
}

@article{MeierThomas,
author = {Meier, Thomas and Gandt, Karima and Hay, Thomas and Echterhof, Thomas},
title = {{Process Modeling and Simulation of the Radiation in the Electric Arc Furnace}},
journal = {steel research international},
volume = {89},
number = {4},
pages = {1700487},
doi = {10.1002/srin.201700487},
url = {https://onlinelibrary.wiley.com/doi/abs/10.1002/srin.201700487},
year = {2018}
}

@article{MarulandaDurango,
author = {Marulanda-Durango, Jesser and Zuluaga-Ríos, C.D.},
year = {2022},
month = {12},
pages = {100850},
title = {{A Meta-heuristic Optimization-based Method for Parameter Estimation of an Electric Arc Furnace Model}},
volume = {17},
journal = {Results in Engineering},
doi = {10.1016/j.rineng.2022.100850}
}

@article{identification_electric_arc_furnace,
  title={Modeling and Simulation of an Electric Arc Furnace Process},
  author={Johannes Gerhardt Bekker and Ian Keith Craig and Petrus Christiaan Pistorius},
  journal={ISIJ International},
  volume={39},
  number={1},
  pages={23-32},
  year={1999},
  doi={10.2355/isijinternational.39.23}
}

@article{Lpez2014CALIBRATIONOP,
author = {Álvarez, Mauricio and Henao Baena, Carlos and Marulanda-Durango, Jesser},
year = {2014},
month = {12},
volume = {11},
number = {22},
pages = {35--44},
title = {Calibration of parameters for electric arc furnace model using simulation and neural networks},
journal = {Revista EIA}
}

@article{Beldzinski,
author = {Bełdziński, Michał and Galewski, Marek and Kossowski, Igor},
title = {Dependence on parameters for nonlinear equations—Abstract principles and applications},
journal = {Mathematical Methods in the Applied Sciences},
volume = {45},
number = {3},
pages = {1668-1686},
doi = {10.1002/mma.7882},
url = {https://onlinelibrary.wiley.com/doi/abs/10.1002/mma.7882},
year = {2022}
}

@Inbook{book_sontag,
author="Sontag, Eduardo D.",
title="Systems",
bookTitle="Mathematical Control Theory: Deterministic Finite Dimensional Systems",
year="1998",
publisher="Springer New York",
address="New York, NY",
pages="25--80",
isbn="978-1-4612-0577-7",
doi="10.1007/978-1-4612-0577-7_2",
url="https://doi.org/10.1007/978-1-4612-0577-7_2"
}

@misc{decarvalhoneto2025generalversioncaratheodorysexistence,
      title={{A General Version of Carathéodory's Existence and Uniqueness Theorem}}, 
      author={Paulo M. de Carvalho-Neto and Cícero L. Frota and Pedro G. P. Torelli},
      year={2025},
      eprint={2505.24516},
      archivePrefix={arXiv},
      primaryClass={math.CA},
      url={https://arxiv.org/abs/2505.24516}, 
}

@book{o2013existence,
  title={{Existence Theory for Nonlinear Ordinary Differential Equations}},
  author={O'Regan, D.},
  isbn={9789401715171},
  lccn={97008146},
  series={Mathematics and Its Applications},
  url={https://books.google.de/books?id=ijjpCAAAQBAJ},
  year={2013},
  publisher={Springer Netherlands}
}

@Inbook{roubicek,
author="Roubíček, Tomáš",
title="Preliminary general material",
bookTitle="Nonlinear Partial Differential Equations with Applications",
year="2013",
publisher="Springer Basel",
address="Basel",
pages="1--27",
isbn="978-3-0348-0513-1",
doi="10.1007/978-3-0348-0513-1_1",
url="https://doi.org/10.1007/978-3-0348-0513-1_1"
}

@book{alt2006lineare,
  title={Linear Functional Analysis: An Application-Oriented Introduction},
  author={Alt, H.W. and Nürnberg, R.},
  isbn={9781447172802},
  lccn={2016944464},
  series={Universitext},
  url={https://books.google.de/books?id=mzilDAAAQBAJ},
  year={2016},
  publisher={Springer-Verlag Berlin Heidelberg}
}

@article{ern2016mollificationstronglylipschitzdomains,
url = {https://doi.org/10.1515/cmam-2015-0034},
title = {Mollification in Strongly Lipschitz Domains with Application to Continuous and Discrete De Rham Complexes},
author = {Alexandre Ern and Jean-Luc Guermond},
pages = {51--75},
volume = {16},
number = {1},
journal = {Computational Methods in Applied Mathematics},
doi = {10.1515/cmam-2015-0034},
year = {2016},
lastchecked = {2026-11-09}
}

@article{Murat2003,
author = {Murat, François and Trombetti, Cristina},
year = {2003},
month = {10},
pages = {},
title = {{A chain rule formula for the composition of a vector-valued function by a piecewise smooth function}},
volume = {3},
journal = {Bollettino della Unione Matematica Italiana. Serie VIII. Sezione B. Articoli di Ricerca Matematica}
}

@book{Stein,
 ISBN = {9780691080796},
 URL = {http://www.jstor.org/stable/j.ctt1bpmb07},
 author = {Elias M. Stein},
 publisher = {Princeton University Press},
 title = {{Singular Integrals and Differentiability Properties of Functions (PMS-30)}},
 urldate = {2026-10-31},
 year = {1970}
}

@book{Brzis2010FunctionalAS,
  title={Functional Analysis, Sobolev Spaces and Partial Differential Equations},
  author={Brezis, H.},
  isbn={9780387709147},
  lccn={2010938382},
  series={Universitext},
  url={https://books.google.de/books?id=1bFCAAAAQBAJ},
  year={2010},
  publisher={Springer New York}
}

@book{Triebel1978InterpolationTF,
  title={Interpolation Theory, Function Spaces, Differential Operators},
  author={Triebel, H.},
  isbn={9783335004202},
  lccn={95202085},
  url={https://books.google.de/books?id=ZzvvAAAAMAAJ},
  year={1995},
  publisher={Huethig \& Wepf}
}

@article{guehring2019errorboundsapproximationsdeep,
author = {Gühring, Ingo and Kutyniok, Gitta and Petersen, Philipp},
title = {Error bounds for approximations with deep ReLU neural networks in Ws,p norms},
journal = {Analysis and Applications},
volume = {18},
number = {05},
pages = {803-859},
year = {2020},
doi = {10.1142/S0219530519410021},
URL = { 
        https://doi.org/10.1142/S0219530519410021
},
eprint = { 
        https://doi.org/10.1142/S0219530519410021
}
}

@misc{yang2026deepneuralnetworksgeneral,
      title={{Deep Neural Networks with General Activations: Super-Convergence in Sobolev Norms}}, 
      author={Yahong Yang and Juncai He},
      year={2026},
      eprint={2508.05141},
      archivePrefix={arXiv},
      primaryClass={cs.LG},
      url={https://arxiv.org/abs/2508.05141}, 
}

@misc{hendrycks2023gaussianerrorlinearunits,
      title={Gaussian Error Linear Units (GELUs)}, 
      author={Dan Hendrycks and Kevin Gimpel},
      year={2023},
      eprint={1606.08415},
      archivePrefix={arXiv},
      primaryClass={cs.LG},
      url={https://arxiv.org/abs/1606.08415}, 
}

@article{LiJiyong,
author = {Li, Jiyong},
year = {2025},
month = {07},
pages = {},
title = {{Uniform error bounds of a nested Picard iterative integrator for the Klein-Gordon-Zakharov system in the subsonic limit regime}},
volume = {51},
journal = {Advances in Computational Mathematics},
doi = {10.1007/s10444-025-10251-x}
}

@book{ShenSpectral,
author = {Shen, Jie and Tang, Tao and Wang, Li-Lian},
year = {2011},
month = {01},
publisher = {Springer},
title = {{Spectral Methods: Algorithms, Analysis and Applications}},
volume = {41},
isbn = {978-3-540-71040-0},
series = {Springer Series in Computational Mathematics},
doi = {10.1007/978-3-540-71041-7}
}

@article{kovachki2021neural,
author = {Kovachki, Nikola and Li, Zongyi and Liu, Burigede and Azizzadenesheli, Kamyar and Bhattacharya, Kaushik and Stuart, Andrew and Anandkumar, Anima},
title = {Neural operator: learning maps between function spaces with applications to PDEs},
year = {2023},
issue_date = {January 2023},
publisher = {JMLR.org},
volume = {24},
number = {1},
issn = {1532-4435},
journal = {J. Mach. Learn. Res.},
month = jan,
articleno = {89},
numpages = {97}
}

@inproceedings{
li2020neuraloperatorgraphkernel,
title={Neural Operator: Graph Kernel Network for Partial Differential Equations},
author={Anima Anandkumar and Kamyar Azizzadenesheli and Kaushik Bhattacharya and Nikola Kovachki and Zongyi Li and Burigede Liu and Andrew Stuart},
booktitle={ICLR 2020 Workshop on Integration of Deep Neural Models and Differential Equations},
year={2019},
url={https://openreview.net/forum?id=fg2ZFmXFO3}
}

@article{lanthaler2024nonlocalitynonlinearityimpliesuniversality,
  title={Nonlocality and nonlinearity implies universality in operator learning},
  author={Lanthaler, Samuel and Li, Zongyi and Stuart, Andrew M},
  journal={Constructive Approximation},
  volume={62},
  number={2},
  pages={261--303},
  year={2025},
  publisher={Springer}
}

@article{Nganyu_Tanyu_2023,
   title={{Deep learning methods for partial differential equations and related parameter identification problems}},
   volume={39},
   ISSN={1361-6420},
   url={http://dx.doi.org/10.1088/1361-6420/ace9d4},
   DOI={10.1088/1361-6420/ace9d4},
   number={10},
   journal={Inverse Problems},
   publisher={IOP Publishing},
   author={Nganyu Tanyu, Derick and Ning, Jianfeng and Freudenberg, Tom and Heilenkötter, Nick and Rademacher, Andreas and Iben, Uwe and Maass, Peter},
   year={2023},
   month=Aug, pages={103001} }

@article{kovachki2024operatorlearningalgorithmsanalysis,
  title={Operator learning: Algorithms and analysis},
  author={Kovachki, Nikola B and Lanthaler, Samuel and Stuart, Andrew M},
  journal={Handbook of Numerical Analysis},
  volume={25},
  pages={419--467},
  year={2024},
  publisher={Elsevier}
}

@article{azizzadenesheli2024neuraloperatorsacceleratingscientific,
  title={Neural operators for accelerating scientific simulations and design},
  author={Azizzadenesheli, Kamyar and Kovachki, Nikola and Li, Zongyi and Liu-Schiaffini, Miguel and Kossaifi, Jean and Anandkumar, Anima},
  journal={Nature Reviews Physics},
  volume={6},
  number={5},
  pages={320--328},
  year={2024},
  publisher={Nature Publishing Group UK, London}
}

@article{lanthaler2023operatorlearningpcanetupper,
  title={Operator learning with PCA-Net: upper and lower complexity bounds},
  author={Lanthaler, Samuel},
  journal={Journal of Machine Learning Research},
  volume={24},
  number={318},
  pages={1--67},
  year={2023}
}

@article{raonic2023convolutionalneuraloperatorsrobust,
  title={Convolutional neural operators for robust and accurate learning of pdes},
  author={Raonić, Bogdan and Molinaro, Roberto and De Ryck, Tim and Rohner, Tobias and Bartolucci, Francesca and Alaifari, Rima and Mishra, Siddhartha and De Bézenac, Emmanuel},
  journal={Advances in Neural Information Processing Systems},
  volume={36},
  pages={77187--77200},
  year={2023}
}

@article{cao2023lnolaplaceneuraloperator,
  title={Laplace neural operator for solving differential equations},
  author={Cao, Qianying and Goswami, Somdatta and Karniadakis, George Em},
  journal={Nature Machine Intelligence},
  volume={6},
  number={6},
  pages={631--640},
  year={2024},
  publisher={Nature Publishing Group UK, London}
}

@article{liu2023renderuntonumericsorthogonal,
  title={Render unto numerics: Orthogonal polynomial neural operator for PDEs with nonperiodic boundary conditions},
  author={Liu, Ziyuan and Wang, Haifeng and Zhang, Hong and Bao, Kaijun and Qian, Xu and Song, Songhe},
  journal={SIAM Journal on Scientific Computing},
  volume={46},
  number={4},
  pages={C323--C348},
  year={2024},
  publisher={SIAM}
}

@article{tripura2022waveletneuraloperatorneural,
author = {Tripura, Tapas and Chakraborty, Souvik},
year = {2023},
month = {02},
pages = {115783},
title = {Wavelet Neural Operator for solving parametric partial differential equations in computational mechanics problems},
volume = {404},
journal = {Computer Methods in Applied Mechanics and Engineering},
doi = {10.1016/j.cma.2022.115783}
}

@inproceedings{
wagner2024neural,
title={{Neural Operators as Fast Surrogate Models for the Transmission Loss of Parameterized Sonic Crystals}},
author={Jakob Elias Wagner and Samuel Burbulla and Miguel de Benito Delgado and Johannes D. Schmid},
booktitle={NeurIPS 2024 Workshop on Data-driven and Differentiable Simulations, Surrogates, and Solvers},
year={2024},
url={https://openreview.net/forum?id=ltXAQvjhjt}
}

@inproceedings{Yin_2022,
    author = {Yin, Ziyi and Siahkoohi, Ali and Louboutin, Mathias and Herrmann, Felix J.},
    title = {Learned coupled inversion for carbon sequestration monitoring and forecasting with Fourier neural operators},
    booktitle = {Second International Meeting for Applied Geoscience and Energy},
    series = {SEG Technical Program Expanded Abstracts},
    pages = {467-472},
    year = {2022},
    month = {08},
    doi = {10.1190/image2022-3722848.1},
    url = {https://doi.org/10.1190/image2022-3722848.1},
}

@inproceedings{
behroozi2025sensitivityconstrainedfourierneuraloperators,
title={Sensitivity-Constrained Fourier Neural Operators for Forward and Inverse Problems in Parametric Differential Equations},
author={Abdolmehdi Behroozi and Chaopeng Shen and Daniel Kifer},
booktitle={The Thirteenth International Conference on Learning Representations},
year={2025},
url={https://openreview.net/forum?id=DPzQ5n3mNm}
}

@article{Yin_2023,
   title={{Solving multiphysics-based inverse problems with learned surrogates and constraints}},
   volume={10},
   ISSN={2213-7467},
   url={http://dx.doi.org/10.1186/s40323-023-00252-0},
   DOI={10.1186/s40323-023-00252-0},
   number={1},
   journal={Advanced Modeling and Simulation in Engineering Sciences},
   publisher={Springer Science and Business Media LLC},
   author={Yin, Ziyi and Orozco, Rafael and Louboutin, Mathias and Herrmann, Felix J.},
   year={2023},
   month=Oct }

@article{Cao_2023,
   title={{Residual-based error correction for neural operator accelerated infinite-dimensional Bayesian inverse problems}},
   volume={486},
   ISSN={0021-9991},
   url={http://dx.doi.org/10.1016/j.jcp.2023.112104},
   DOI={10.1016/j.jcp.2023.112104},
   journal={Journal of Computational Physics},
   publisher={Elsevier BV},
   author={Cao, Lianghao and O’Leary-Roseberry, Thomas and Jha, Prashant K. and Oden, J. Tinsley and Ghattas, Omar},
   year={2023},
   month=Aug, pages={112104} }

@article{scherzer2026neuraloperatorssolvingnonlinear,
  title={Neural operators for solving nonlinear inverse problems},
  author={Scherzer, Otmar and Vu, Thi Lan Nhi and Yan, Jikai},
  journal={Numerical Functional Analysis and Optimization},
  volume={47},
  number={11-14},
  pages={315--356},
  year={2026},
  publisher={Taylor \& Francis}
}

@article{wang2024latentneuraloperatorsolving,
  title={Latent neural operator for solving forward and inverse pde problems},
  author={Wang, Tian and Wang, Chuang},
  journal={Advances in Neural Information Processing Systems},
  volume={37},
  pages={33085--33107},
  year={2024}
}

@inproceedings{molinaro2023neuralinverseoperatorssolving,
  title={Neural Inverse Operators for Solving PDE Inverse Problems},
  author={Roberto Molinaro and Yu-Nan Yang and Bjorn Engquist and Siddhartha Mishra},
  booktitle={International Conference on Machine Learning},
  year={2023},
  url={https://api.semanticscholar.org/CorpusID:256274656}
}

@article{Lu_2021,
   title={{Learning nonlinear operators via DeepONet based on the universal approximation theorem of operators}},
   volume={3},
   ISSN={2522-5839},
   url={http://dx.doi.org/10.1038/s42256-021-00302-5},
   DOI={10.1038/s42256-021-00302-5},
   number={3},
   journal={Nature Machine Intelligence},
   publisher={Springer Science and Business Media LLC},
   author={Lu, Lu and Jin, Pengzhan and Pang, Guofei and Zhang, Zhongqiang and Karniadakis, George Em},
   year={2021},
   month=Mar, pages={218–229} }

@article{jin2022mionetlearningmultipleinputoperators,
  title={MIONet: Learning multiple-input operators via tensor product},
  author={Jin, Pengzhan and Meng, Shuai and Lu, Lu},
  journal={SIAM Journal on Scientific Computing},
  volume={44},
  number={6},
  pages={A3490--A3514},
  year={2022},
  publisher={SIAM}
}

@misc{liu2026inverseneuraloperatorode,
      title={{Inverse Neural Operator for ODE Parameter Optimization}}, 
      author={Zhi-Song Liu and Wenqing Peng and Helmi Toropainen and Ammar Kheder and Andreas Rupp and Holger Froning and Xiaojie Lin and Michael Boy},
      year={2026},
      eprint={2603.11854},
      archivePrefix={arXiv},
      primaryClass={cs.LG},
      url={https://arxiv.org/abs/2603.11854}, 
}

@inproceedings{li2021fourierneuraloperatorparametric,
  author       = {Zongyi Li and
                  Nikola Borislavov Kovachki and
                  Kamyar Azizzadenesheli and
                  Burigede Liu and
                  Kaushik Bhattacharya and
                  Andrew M. Stuart and
                  Anima Anandkumar},
  title        = {Fourier Neural Operator for Parametric Partial Differential Equations},
  booktitle    = {9th International Conference on Learning Representations, {ICLR} 2021,
                  Virtual Event, Austria, May 3-7, 2021},
  publisher    = {OpenReview.net},
  year         = {2021},
  url          = {https://openreview.net/forum?id=c8P9NQVtmnO},
  bibsource    = {dblp computer science bibliography, https://dblp.org}
}

@article{kovachki2021universalapproximationerrorbounds,
  title={On universal approximation and error bounds for Fourier neural operators},
  author={Kovachki, Nikola and Lanthaler, Samuel and Mishra, Siddhartha},
  journal={Journal of Machine Learning Research},
  volume={22},
  number={290},
  pages={1--76},
  year={2021}
}

@inproceedings{furuya2024quantitativeapproximationneuraloperators,
  title={Quantitative approximation for neural operators in nonlinear parabolic equations},
  author={Furuya, Takashi and Taniguchi, Koichi and Okuda, Satoshi},
  booktitle={International Conference on Learning Representations},
  pages={85528--85556},
  year={2025}
}

@article{lanthaler2025parametriccomplexityoperatorlearning,
  title={The parametric complexity of operator learning},
  author={Lanthaler, Samuel and Stuart, Andrew M},
  journal={IMA Journal of Numerical Analysis},
  volume={46},
  number={2},
  pages={647--712},
  year={2026},
  publisher={Oxford University Press}
}

@article{lanthaler2022errorestimatesdeeponetsdeep,
  title={Error estimates for deeponets: A deep learning framework in infinite dimensions},
  author={Lanthaler, Samuel and Mishra, Siddhartha and Karniadakis, George E},
  journal={Transactions of Mathematics and its Applications},
  volume={6},
  number={1},
  pages={tnac001},
  year={2022},
  publisher={Oxford University Press}
}

@misc{kow2026quantitativeapproximationtheoremneural,
      title={{Quantitative approximation theorem of neural operators for the Dirichlet-to-Neumann map}}, 
      author={Pu-Yun Kow and Andrea Aspri and Elisa Francini and Jenn-Nan Wang},
      year={2026},
      eprint={2610.35404},
      archivePrefix={arXiv},
      primaryClass={math.AP},
      url={https://arxiv.org/abs/2610.35404}, 
}

@misc{furuya2026approximationtheorylaplacianbasedneural,
      title={{Approximation Theory of Laplacian-Based Neural Operators for Reaction-Diffusion System}}, 
      author={Takashi Furuya and Ryo Ozawa and Jenn-Nan Wang},
      year={2026},
      eprint={2605.12025},
      archivePrefix={arXiv},
      primaryClass={cs.LG},
      url={https://arxiv.org/abs/2605.12025}, 
}

@misc{clason2021regularizationinverseproblems,
      title={{Regularization of Inverse Problems}}, 
      author={Christian Clason},
      year={2021},
      eprint={2001.00617},
      archivePrefix={arXiv},
      primaryClass={math.FA},
      url={https://arxiv.org/abs/2001.00617}, 
}

@book{Tarantola2004InversePT,
author = {Tarantola, Albert},
title = {Inverse Problem Theory and Methods for Model Parameter Estimation},
publisher = {Society for Industrial and Applied Mathematics},
year = {2005},
doi = {10.1137/1.9780898717921},
address = {},
edition   = {},
URL = {https://epubs.siam.org/doi/abs/10.1137/1.9780898717921},
}

@article{Rieder2003TheAI,
 ISSN = {00255718, 10886842},
 URL = {http://www.jstor.org/stable/4099840},
 author = {Andreas Rieder and Thomas Schuster},
 journal = {Mathematics of Computation},
 number = {243},
 pages = {1399--1415},
 publisher = {American Mathematical Society},
 title = {The Approximate Inverse in Action II: Convergence and Stability},
 urldate = {2026-11-09},
 volume = {72},
 year = {2003}
}

@misc{rothermel2019parameterestimationtemperaturedependent,
      title={{Parameter estimation of temperature dependent material parameters in the cooling process of TMCP steel plates}}, 
      author={Dimitri Rothermel and Thomas Schuster and Roland Schorr and Martin Peglow},
      year={2019},
      eprint={1912.12166},
      archivePrefix={arXiv},
      primaryClass={math.NA},
      url={https://arxiv.org/abs/1912.12166}, 
}

@article{Rothermel2020SolvingAI,
  title={Solving an inverse heat convection problem with an implicit forward operator by using a projected quasi-Newton method},
  author={Dimitri Rothermel and Thomas Schuster},
  journal={Inverse Problems},
  year={2020},
  volume={37},
  url={https://api.semanticscholar.org/CorpusID:222133406}
}

@incollection{klein2019sequentialsubspaceoptimizationrecovering,
  title={Sequential subspace optimization for recovering stored energy functions in hyperelastic materials from time-dependent data},
  author={Klein, Rebecca and Schuster, Thomas and Wald, Anne},
  booktitle={Time-dependent problems in imaging and parameter identification},
  pages={165--190},
  year={2021},
  publisher={Springer}
}

@article{Leit_o_2016,
   title={{On projective Landweber–Kaczmarz methods for solving systems of nonlinear ill-posed equations}},
   volume={32},
   ISSN={1361-6420},
   url={http://dx.doi.org/10.1088/0266-5611/32/2/025004},
   DOI={10.1088/0266-5611/32/2/025004},
   number={2},
   journal={Inverse Problems},
   publisher={IOP Publishing},
   author={Leitão, A and Svaiter, B F},
   year={2016},
   month=jan, pages={025004}
}

@article{Kowar2002ConvergenceAO,
author = {Richard, Kowar and Scherzer, O.},
year = {2002},
month = {01},
pages = {},
title = {Convergence analysis of a Landweber-Kaczmarz method for solving nonlinear ill-posed problems},
volume = {23},
journal = {Ill-Posed and Inverse Problems}
}

@article{NanYu2026InverseProblemsandImaging,
title = {{Landweber–Kaczmarz method for inverse problems using multiple repeated measurements}},
journal = {Inverse Problems and Imaging},
volume = {20},
number = {0},
pages = {1-24},
year = {2026},
issn = {1930-8337},
doi = {10.3934/ipi.2025016},
url = {https://www.aimsciences.org/article/id/67f8db1e45e64536d7529a5f},
author = {Nan Yu and Yuxin Xia and Zhenwu Fu and Yong Chen}
}

@article{Li_2018,
   title={{The Averaged Kaczmarz Iteration for Solving Inverse Problems}},
   volume={11},
   ISSN={1936-4954},
   url={http://dx.doi.org/10.1137/17M1146178},
   DOI={10.1137/17m1146178},
   number={1},
   journal={SIAM Journal on Imaging Sciences},
   publisher={Society for Industrial & Applied Mathematics (SIAM)},
   author={Li, Housen and Haltmeier, Markus},
   year={2018},
   month=Jan, pages={618–642} }

@article{iglesias2025deeponetacceleratedbayesianinversionmoving,
  title={DeepONet-accelerated Bayesian inversion for moving boundary problems},
  author={Iglesias, Marco A and Causon, Michael E and Matveev, Mikhail Y and Endruweit, Andreas and Tretyakov, Michael V},
  journal={Computer Methods in Applied Mechanics and Engineering},
  volume={457},
  pages={118988},
  year={2026},
  publisher={Elsevier}
}

@article{zhou2025parameterestimationstructuraldynamics,
  title={Parameter estimation of structural dynamics with neural operators enabled surrogate modeling},
  author={Zhou, Mingyuan and Song, Haoze and Ye, Wenjing and Wang, Wei and Lai, Zhilu},
  journal={Mechanical Systems and Signal Processing},
  volume={237},
  pages={112914},
  year={2025},
  publisher={Elsevier}
}

@article{Iglesias_2013,
   title={{Ensemble Kalman methods for inverse problems}},
   volume={29},
   ISSN={1361-6420},
   url={http://dx.doi.org/10.1088/0266-5611/29/4/045001},
   DOI={10.1088/0266-5611/29/4/045001},
   number={4},
   journal={Inverse Problems},
   publisher={IOP Publishing},
   author={Iglesias, Marco A and Law, Kody J H and Stuart, Andrew M},
   year={2013},
   month=Mar, pages={045001} }

@article{Iglesias_2016,
title={{A regularizing iterative ensemble Kalman method for PDE-constrained inverse problems}},
volume={32},
ISSN={1361-6420},
url={http://dx.doi.org/10.1088/0266-5611/32/2/025002},
DOI={10.1088/0266-5611/32/2/025002},
number={2},
journal={Inverse Problems},
publisher={IOP Publishing},
author={Iglesias, Marco A},
year={2016},
month=Jan, pages={025002} }

@misc{iglesias2020adaptiveregularisationensemblekalman,
      title={{Adaptive regularisation for ensemble Kalman inversion}}, 
      author={Marco Iglesias and Yuchen Yang},
      year={2020},
      eprint={2006.14980},
      archivePrefix={arXiv},
      primaryClass={math.NA},
      url={https://arxiv.org/abs/2006.14980}, 
}

@article{Causon2024RealtimeBI,
title = {Real-time Bayesian inversion in resin transfer moulding using neural surrogates},
journal = {Composites Part A: Applied Science and Manufacturing},
volume = {185},
pages = {108355},
year = {2024},
issn = {1359-835X},
doi = {https://doi.org/10.1016/j.compositesa.2024.108355},
url = {https://www.sciencedirect.com/science/article/pii/S1359835X2400352X},
author = {M.E. Causon and M.A. Iglesias and M.Y. Matveev and A. Endruweit and M.V. Tretyakov}
}

@article{weissmann2022gradient,
  title={Gradient flow structure and convergence analysis of the ensemble Kalman inversion for nonlinear forward models},
  author={Weissmann, Simon},
  journal={Inverse Problems},
  volume={38},
  number={10},
  pages={105011},
  year={2022},
  publisher={IOP Publishing}
}

@article{Kennedy2001BayesianCO,
  title={{Bayesian calibration of computer models}},
  author={Marc C. Kennedy and Anthony O'Hagan},
  journal={Journal of the Royal Statistical Society: Series B (Statistical Methodology)},
  year={2001},
  volume={63},
  url={https://api.semanticscholar.org/CorpusID:119562136}
}

@article{Anthony,
author = {Brynjarsdóttir, Jenný and {O'Hagan}, Anthony},
year = {2014},
month = {10},
pages = {},
title = {{Learning about physical parameters: the importance of model discrepancy}},
volume = {30},
journal = {Inverse Problems},
doi = {10.1088/0266-5611/30/11/114007}
}

@article{schillings2017analysis,
  title={Analysis of the ensemble {K}alman filter for inverse problems},
  author={Schillings, Claudia and Stuart, Andrew M.},
  journal={SIAM Journal on Numerical Analysis},
  volume={55},
  number={3},
  pages={1264--1290},
  year={2017},
  doi={10.1137/16M105959X}
}

@article{laurent2000adaptive,
author = {B. Laurent and P. Massart},
title = {{Adaptive estimation of a quadratic functional by model selection}},
volume = {28},
journal = {The Annals of Statistics},
number = {5},
publisher = {Institute of Mathematical Statistics},
pages = {1302 -- 1338},
year = {2000},
doi = {10.1214/aos/1015957395},
URL = {https://doi.org/10.1214/aos/1015957395}
}
